\documentclass[a4paper,oneside,11pt]{article}

\usepackage[a4paper]{geometry}
\usepackage{aeguill}                 
\usepackage{graphicx}
\usepackage{amsmath,amsfonts,amssymb,amsthm}
\usepackage{pstricks} 
\usepackage{epstopdf}
\usepackage{srcltx}
\usepackage[utf8]{inputenc} 
\usepackage[T1]{fontenc} 
\usepackage{hyperref} 
\usepackage[english]{babel} 
\usepackage{comment}

\title{Boundary amenability and measure equivalence rigidity among two-dimensional Artin groups of hyperbolic type}
\author{Camille Horbez and Jingyin Huang}

\usepackage[all]{xy}
\usepackage{bbm}
\begin{document}
\maketitle
\newtheorem{de}{Definition} [section]
\newtheorem{theo}[de]{Theorem} 
\newtheorem{prop}[de]{Proposition}
\newtheorem{lemma}[de]{Lemma}
\newtheorem{cor}[de]{Corollary}
\newtheorem{propd}[de]{Proposition-Definition}
\newtheorem{conj}[de]{Conjecture}
\newtheorem{claim}{Claim}
\newtheorem*{claim2}{Claim}

\newtheorem{theointro}{Theorem}
\newtheorem*{defintro}{Definition}
\newtheorem{corintro}[theointro]{Corollary}

\theoremstyle{remark}
\newtheorem{rk}[de]{Remark}
\newtheorem{ex}[de]{Example}
\newtheorem{question}[de]{Question}

\newtheorem*{rkintro}{Remark}

\normalsize

\newcommand{\Aut}{\mathrm{Aut}}
\newcommand{\Out}{\mathrm{Out}}
\newcommand{\Inn}{\mathrm{Inn}}
\newcommand{\stab}{\operatorname{Stab}}
\newcommand{\pstab}{\operatorname{Pstab}}
\newcommand{\dunion}{\sqcup}
\newcommand{\eps}{\varepsilon}
\renewcommand{\epsilon}{\varepsilon}
\newcommand{\calf}{\mathcal{F}}
\newcommand{\cali}{\mathcal{I}}
\newcommand{\caly}{\mathcal{Y}}
\newcommand{\calx}{\mathcal{X}}
\newcommand{\calz}{\mathcal{Z}}
\newcommand{\calo}{\mathcal{O}}
\newcommand{\calb}{\mathcal{B}}
\newcommand{\calq}{\mathcal{Q}}
\newcommand{\calu}{\mathcal{U}}
\newcommand{\call}{\mathcal{L}}
\newcommand{\bbR}{\mathbb{R}}
\newcommand{\bbZ}{\mathbb{Z}}
\newcommand{\bbD}{\mathbb{D}}
\newcommand{\NT}{\mathrm{NT}}
\newcommand{\cat}{\mathrm{CAT}(-1)}
\newcommand{\CAT}{\mathrm{CAT}(0)}
\newcommand{\actson}{\curvearrowright}
\newcommand{\caln}{\mathcal{N}}
\newcommand{\calg}{\mathcal{G}}
\newcommand{\Prob}{\mathrm{Prob}}
\newcommand{\calt}{\mathcal{T}}
\newcommand{\calc}{\mathcal{C}}
\newcommand{\adm}{\mathrm{adm}}
\newcommand{\cala}{\mathcal{A}}
\newcommand{\cals}{\mathcal{S}}
\newcommand{\calh}{\mathcal{H}}
\newcommand{\Stab}{\mathrm{Stab}}
\newcommand{\Isom}{\mathrm{Isom}}
\newcommand{\cod}{\mathbb{D}^*}
\newcommand{\horo}{\overline{\cod}^h}
\newcommand{\horox}{\overline{X}^h}
\newcommand{\bdd}{\mathrm{bdd}}
\newcommand{\calp}{\mathcal{P}}
\newcommand{\flex}{\mathrm{flex}}
\newcommand{\bv}{\mathcal{B}V}
\newcommand{\Fix}{\mathrm{Fix}}
\newcommand{\bad}{\mathrm{bad}}
\newcommand{\todo}{\textbf{TO DO:~}}
\newcommand{\Conv}{\mathrm{Conv}}
\newcommand{\cent}{c}
\newcommand{\cald}{\mathcal{D}}
\newcommand{\rad}{r}
\newcommand{\good}{\mathrm{good}}
\newcommand{\Faces}{\mathrm{Faces}}
\newcommand{\del}{\mathbb{D}_\Gamma}
\newcommand{\Mod}{\mathrm{Mod}}
\newcommand{\Comm}{\mathrm{Comm}}
\newcommand{\si}{\sigma}
\newcommand{\st}{\mathrm{st}}
\newcommand{\PML}{\mathrm{PML}}
\newcommand{\Inv}{\mathrm{SInv}}
\newcommand{\Homeo}{\mathrm{Homeo}}
\newcommand{\Cay}{\mathfrak{C}}
\newcommand{\lk}{\mathrm{lk}}

\makeatletter
\edef\@tempa#1#2{\def#1{\mathaccent\string"\noexpand\accentclass@#2 }}
\@tempa\rond{017}
\makeatother

\begin{abstract}
We study $2$-dimensional Artin groups of hyperbolic type from the viewpoint of measure equivalence, and establish rigidity theorems.

We first prove that they are boundary amenable. So is every group acting discretely by simplicial isometries on a connected piecewise hyperbolic $\cat$ simplicial complex with countably many simplices in finitely many isometry types, assuming that vertex stabilizers are boundary amenable. Consequently, they satisfy the Novikov conjecture.

We then show that measure equivalent $2$-dimensional Artin groups of hyperbolic type have isomorphic fixed set graphs -- an analogue of the curve graph, introduced by Crisp. This yields classification results.

We obtain strong rigidity theorems. Let $G=G_\Gamma$ be a $2$-dimensional Artin group of hyperbolic type, with $\Out(G)$ finite. When the automorphism groups of the fixed set graph and of the Cayley complex $\mathfrak{C}$ coincide, every countable group $H$ which is measure equivalent to $G$, is commensurable to a lattice in $\Aut(\mathfrak{C})$. This happens whenever $\Gamma$ is triangle-free with all labels at least $3$ -- unless $G$ is commensurable to the direct sum of $\mathbb{Z}$ and a free group. When $\Gamma$ satisfies an additional star-rigidity condition, then $\Aut(\mathfrak{C})$ is countable, and $H$ is almost isomorphic to $G$. 

This has applications to orbit equivalence rigidity, and rigidity results for von Neumann algebras associated to ergodic actions of Artin groups. We also derive a rigidity statement regarding possible lattice envelopes of certain Artin groups, and a cocycle superrigidity theorem from higher-rank lattices to $2$-dimensional Artin groups of hyperbolic type. 
\end{abstract}

\section{Introduction} 

The present paper is concerned with measure equivalence rigidity and classification of $2$-dimensional Artin groups of hyperbolic type. This has several applications, for instance to a complete description of possible lattice embeddings of these groups in locally compact second countable groups. We also obtain strong $W^*$-rigidity results for cross-product von Neumann algebras associated to ergodic actions of a large class of Artin groups. Along the way, we establish that all $2$-dimensional Artin groups of hyperbolic type are boundary amenable -- as a byproduct, they satisfy the Novikov conjecture on higher signatures. 

\subsection{Artin groups} 

Let $\Gamma$ be a finite simple labeled graph, where every edge is labeled by an integer at least equal to $2$. The \emph{Artin group with defining graph $\Gamma$}, denoted $G_{\Gamma}$, is the group with the following presentation: its generators are the vertices of $\Gamma$, and its relators are given by $$\underbrace{aba\cdots}_{m}=\underbrace{bab\cdots}_{m}$$ whenever $a$ and $b$ span an edge with label $m$. Artin groups emerged in the study of braids, configurations spaces, hyperplane complements and Coxeter groups (the latter being obtained by adding the extra relation that every generator squares to $1$). 
Most Artin groups are still far from being understood 
 -- we refer to \cite{mccammond2017mysterious} for perspectives.

In the present paper, we will be interested in $2$-dimensional Artin groups of hyperbolic type, i.e.\ Artin groups $G_\Gamma$ with cohomological dimension $2$ and whose associated Coxeter group is Gromov-hyperbolic. We refer to Section~\ref{sec:background} for a combinatorial characterization in terms of the defining graph, and simply mention that this class properly contains all Artin groups such that all labels of the defining graph are at least $4$. 

Two-dimensional Artin groups of hyperbolic type are not Gromov-hyperbolic, but have interesting actions on Gromov-hyperbolic spaces. One example is the \emph{modified Deligne complex} introduced in Charney and Davis \cite{CharneyDavis}, which will play an important role in this paper.

\subsection{Boundary amenability}

Our first goal is to establish a general criterion for boundary amenability of groups acting on $\cat$ complexes, which applies in particular to 2-dimensional Artin groups of hyperbolic type. We recall that a countable group $G$ is \emph{boundary amenable} if it admits a \emph{Borel amenable} action on a nonempty compact space $X$, i.e.\ there exists a sequence of Borel maps $\nu_n:X\to\Prob(G)$ such that for every $g\in G$, one has $$||\nu_n(gx)-g\cdot\nu_n(x)||_1\to 0$$ as $n$ goes to $+\infty$. Here the topology on the group $\Prob(G)$ is the subspace topology from $\ell^1(G)$. The usual definition requires the maps $\nu_n$ to be continuous and the convergence to be uniform in $x$, but this in fact yields an equivalent definition \cite{Ren,Oza}.

One motivation for studying boundary amenability is that every countable group which is boundary amenable satisfies the Novikov conjecture on higher signatures, as follows from works of Yu \cite{Yu}, Higson--Roe \cite{HR} and Higson \cite{Hig}. Another motivation comes from operator algebras: boundary amenability of a countable group is equivalent to the exactness of its reduced $C^*$-algebra \cite{Ana,Oza2}. A third motivation -- which is the main motivation of the present paper -- comes from applications to measured group theory.

Boundary amenability has already been established for many classes of groups, among which linear groups \cite{GHW}, hyperbolic groups \cite{Ada2} and more generally relatively hyperbolic groups with boundary amenable parabolic subgroups \cite{Oza},  locally compact groups acting properly by automorphisms on buildings \cite{Lec}, mapping class groups of finite-type surfaces \cite{Kid-memoir,Ham}, outer automorphism groups of finitely generated free groups (and more generally of right-angled Artin groups and of torsion-free hyperbolic groups or toral relatively hyperbolic groups) \cite{BGH}. It has also been established for groups acting properly and cocompactly on $\CAT$ cube complexes \cite{CN}, and more generally for countable groups acting on finite-dimensional $\CAT$ cube complexes with boundary amenable vertex stabilizers \cite{GN}; in particular all Artin groups of type FC are boundary amenable \cite{GN}. On the other hand, Gromov's monster containing a properly embedded expander fails to be boundary amenable \cite{Gro2,AD,Osa}.

\begin{theointro}[Theorems~\ref{theo:exactness} and~\ref{theo:amenability-boundary}]\label{theo:boundary-amenable-intro}
Let $X$ be a connected piecewise hyperbolic $\cat$ simplicial complex with countably many simplices that belong to finitely many isometry types. Let $G$ be a countable group acting discretely on $X$ by simplicial isometries.
\begin{enumerate}
\item If the $G$-stabilizer of every vertex of $X$ is boundary amenable, then $G$ is boundary amenable.
\item If the $G$-stabilizer of every edge of $X$ is amenable, then the $G$-action on the visual boundary $\partial_\infty X$ is Borel amenable.
\end{enumerate}
\end{theointro}

Notice that none of the two statements directly implies the other; in particular, the visual boundary $\partial_\infty X$ need not be compact, as $X$ is not assumed to be proper.

Using the modified Deligne complex constructed by Charney and Davis \cite{CharneyDavis}, Martin and Przytycki showed in \cite{MP} that every $2$-dimensional Artin group of hyperbolic type admits an action as in Theorem~\ref{theo:boundary-amenable-intro}, where every vertex stabilizer splits as a direct sum of a finitely generated free group and $\mathbb{Z}$. We therefore attain the following corollary (using the fact that boundary amenability is stable under passing to subgroups, as follows from \cite[Proposition~11]{Oza}).

\begin{corintro}[Theorem~\ref{theo:exactness-Artin}]\label{corintro:novikov}
Every group which embeds in a $2$-dimensional Artin group of hyperbolic type is boundary amenable, and therefore satisfies the Novikov conjecture. 
\end{corintro}

We mention that the Novikov conjecture was already known for Artin groups of type FC (those being boundary amenable) and for Artin groups of large type (i.e.\ with all labels of $\Gamma$ at least $3$) by work of Osajda and the second named author \cite{HO2}; Corollary~\ref{corintro:novikov} includes many new cases. 

Let us say a word about our proof of Theorem~\ref{theo:boundary-amenable-intro}. The model case for proofs of boundary amenability is the action of a finitely generated free group $F_N$ on its boundary: the measures $\nu_n$ in this case associate to every point $\xi\in\partial_\infty F_N$ the uniform probability measure on the first $n$ vertices on the geodesic ray from $e$ to $\xi$ in a fixed Cayley tree. The key geometric feature ensuring that these maps are asymptotically equivariant is that any two rays targeting the same point at infinity eventually coincide -- see e.g.\ \cite[Example~2.2]{Oza-icm}.

Assume now that a countable group $G$ acts isometrically on a $\cat$ space $X$. Geodesic rays in $X$ converging to the same boundary point also get arbitrary close to one another, but one needs a bit more structure on $X$ to define the maps $\nu_n$. Kaimanovich established the boundary amenability of $G$ under a  `bounded geometry' assumption (which includes the case when $X$ is proper), or when the group action has a finite critical exponent  \cite{Kai}. In our case, we use the simplicial structure as follows: one first associates to every geodesic ray $r$ the collection of all simplices traversed by $r$ before time $n$. The problem is that when two rays $r$ and $r'$ converge to the same boundary point, it is possible that $r$ follows an edge $e$ very closely while staying in a simplex $\Delta$ that contains $e$ as a face, while $r'$ follows that same edge $e$ very closely while staying in a different simplex $\Delta'$ that also contains $e$ as a face. We tackle this difficulty via a system of \emph{weights} associated to each simplex: in the above case, the probability measure associated to the ray $r$ will give a very small weight to $\Delta$ and a big weight to $e$ -- intuitively, the weight measures how close to the `center' of the simplex the ray passes. Our hypothesis that there are finitely many types of simplices enables us to choose the weights in an appropriate way. In this way we get maps from $\partial_\infty X$ to the set of simplices of $X$, which is why we need the extra assumption that stabilizers of simplices are boundary amenable to conclude. But another problem is that $\partial_\infty X$ need not be compact. This is tackled by working instead in the horofunction compactification of $X$. A key observation coming from work of Maher and Tiozzo \cite{MT} is that every horofunction on $X$ which is unbounded from below naturally corresponds to a point in $\partial_\infty X$, while to every horofunction which is bounded from below, one can naturally associate a point in $X$, defined as the circumcenter of the minset of the horofunction.

\subsection{Measure equivalence rigidity} 

The notion of measure equivalence between countable groups was coined by Gromov in \cite[{Definition~0.5.$\text{E}_1$}]{Gro} as a measurable analogue of the notion of quasi-isometry between finitely generated groups. 

\begin{defintro}[Measure equivalence]
Let $G_1$ and $G_2$ be two countable groups. A \emph{measure equivalence coupling} between $G_1$ and $G_2$ is a standard non-zero measure space $\Sigma$ equipped with a measure-preserving action of $G_1\times G_2$ by Borel automorphisms such that for every $i\in\{1,2\}$, the $G_i$-action on $\Sigma$ is free and has a finite measure fundamental domain. 

Two countable groups $G_1$ and $G_2$ are \emph{measure equivalent} if there exists a measure equivalence coupling between them.
\end{defintro}

For example, two lattices in the same locally compact second countable group are measure equivalent -- as a special case, passing to a finite-index subgroup does not change the measure equivalence class. Taking a quotient by a finite normal subgroup does not change the measure equivalence class either.  

A central question in measured group theory is to determine, given a countable group $G$, to what extent the measure equivalence class of $G$ determines $G$ (up to finite-index subgroups and finite kernels).  
In general it does not: for example, a celebrated theorem of Ornstein and Weiss \cite{OW} asserts that all countably infinite amenable groups are measure equivalent to one another; also, Gaboriau constructed in \cite{Gab} many examples of non-free groups that are measure equivalent to free groups (see also \cite{BTW} for more examples). 

In contrast, a strong rigidity theorem was established by Furman \cite{Fur-me}, building on earlier work of Zimmer \cite{Zim2,Zim3}: if $\Lambda$ is a lattice in a higher rank simple connected Lie group $G$ with finite center, and if $\Lambda'$ is a countable group that is measure equivalent to $\Lambda$, then $\Lambda'$ is commensurable to a lattice in $G$ up to a finite kernel. An even stronger rigidity theorem was later obtained by Kida \cite{Kid} for mapping class groups of finite-type surfaces: with the exception of a few low-complexity cases, if $S$ is a connected oriented hyperbolic surface of finite type, then every countable group that is measure equivalent to the mapping class group $\Mod(S)$, is in fact commensurable to $\Mod(S)$ up to finite kernels. Very recently, Guirardel and the first named author proved in \cite{GH} that for every integer $N\ge 3$, the group $\Out(F_N)$ of outer automorphisms of a free group of rank $N$ satisfies the same strong form of measure equivalence rigidity as in Kida's theorem.

\paragraph*{Measure equivalence rigidity among Artin groups: a special case.} 
 We establish measure equivalence rigidity theorems for $2$-dimensional Artin groups of hyperbolic type, thereby providing a measure-theoretic counterpart to the quasi-isometry rigidity theorems proved by Osajda and the second-named author in \cite{HO}. One important ingredient in Kida's work on mapping class groups \cite{Kid} is Ivanov's theorem stating that every automorphism of the curve graph of $S$ is induced by a homeomorphism of $S$ \cite{Iva}. In our setting, work of Crisp \cite{crisp} provides an analogue of the curve graph (the \emph{fixed set graph} $\Theta_\Gamma$\footnote{A possible definition is the following: vertices are maximal cyclic subgroups of $G_\Gamma$ with non-abelian centralizer; two distinct vertices are joined by an edge when the corresponding subgroups commute. See Section~\ref{sec:theta} for more information.}). By exploiting analogies between Artin groups and mapping class groups, we transfer the measure equivalence rigidity problem for Artin groups to a combinatorial rigidity statement about the fixed set graph, from which we derive superrigidity under suitable assumptions on $\Gamma$. Before explaining how this is done in general, let us first state a concrete measure equivalence rigidity theorem (which ) that we reach under natural assumptions on the defining graph $\Gamma$: for its statement, we say that two countable groups $G_1,G_2$ are \emph{almost isomorphic} if there are finite-index subgroups $G'_i\subseteq G_i$ and finite normal subgroups $F_i\unlhd G'_i$ such that $G'_1/F_1$ and $G'_2/F_2$ are isomorphic. 

\begin{theointro}[Corollary~\ref{cor:main-large}]\label{theointro:sample}
Let $G=G_\Gamma$ be a large type Artin group (i.e.\ all labels of $\Gamma$ are at least $3$) with $\Out(G_\Gamma)$ finite. Assume that $\Gamma$ is triangle-free (i.e. there are no 3-cycles in $\Gamma$) and has more than 2 vertices.
\begin{enumerate}
\item If a countable group $H$ is measure equivalent to $G_\Gamma$, then $H$ is commensurable to a lattice in $\Aut(\Cay_\Gamma)$, the automorphism group of the Cayley complex of $G_\Gamma$ (i.e.\ the universal cover of its presentation complex for the standard Artin presentation).\footnote{Notice that $G_\Gamma$ acts geometrically on $\Cay_\Gamma$ and is therefore itself a cocompact lattice in $\Aut(\Cay_\Gamma)$.}
\item Assume additionally that $\Gamma$ is \emph{star-rigid}, i.e.\ every label-preserving automorphism of $\Gamma$ that fixes the star of a vertex pointwise is the identity. If a countable group $H$ is measure equivalent to $G_\Gamma$, then $H$ is almost isomorphic to $G_\Gamma$.
\end{enumerate} 
\end{theointro}

The assumptions made on $\Gamma$ are natural and can be checked on the defining graph (we suspect however that the triangle-freeness assumption might be removable but more work would be required). When $\Gamma$ has at most $2$ vertices, the associated Artin group is either a free group or commensurable to the direct sum of $\mathbb Z$ with a free group, which is not rigid.

In particular, finiteness of $\Out(G_\Gamma)$ is equivalent to $\Gamma$ being connected, with no separating vertex or edge \cite{crisp}. This assumption is motivated by that in many other cases when measure equivalence rigidity are obtained, i.e. higher-rank lattices, mapping class groups or $\Out(F_N)$, the underlying groups have finite outer automorphism group (and it is natural to require that the group has few algebraic symmetries when looking for rigidity phenomena). Notice however that the groups we work with do not have property (T) (another property often related to rigidity), as they all surject onto $\mathbb{Z}$.

Imposing that $G_\Gamma$ is of large type is also natural, and one should stay aware that a strong form of rigidity cannot hold for all $2$-dimensional Artin groups of hyperbolic type. For instance, if all labels of $\Gamma$ are equal to $2$ (i.e.\ $G_\Gamma$ is a \emph{right-angled Artin group}), then $G_\Gamma$ is not rigid, as every graph product of countably infinite amenable groups over the same defining graph is measure equivalent to $G_\Gamma$, see \cite[Theorem~4.1]{HH}. In the next paragraph of this introduction, we will provide more general statements that actually apply even beyond the large type case. 

We also mention that the \emph{star-rigidity} condition is not merely a technical assumption: from the work of Crisp \cite[p.~1436]{crisp}, one can deduce that this is a necessary and sufficient condition to ensure that $\Aut(\Cay_\Gamma)$ is countable, and actually contains $G_\Gamma$ as a finite-index subgroup. 

As a consequence of our work, we also obtain measure equivalence classification theorems among $2$-dimensional Artin groups of hyperbolic type (see also Corollary~\ref{cor:me-raag} for an application to right-angled Artin groups, which has since been improved in \cite{HH}). In the following statement, the measure equivalence classification coincides with the quasi-isometry classification \cite[Theorem~1.6]{HO}.

\begin{theointro}[Corollary~\ref{cor:me-clttf}]
	\label{theo:1}
Let $G_1$ and $G_2$ be two large type Artin groups whose defining graphs are triangle-free and have more than $2$ vertices. Assume that $\Out(G_1)$ and $\Out(G_2)$ are finite.

Then $G_1$ and $G_2$ are measure equivalent if and only if they are isomorphic.
\end{theointro}

\paragraph*{Measure equivalence rigidity of Artin groups: general statements.}

As mentioned above, our proof of Theorem~\ref{theointro:sample} reduces the measure equivalence rigidity question to a combinatorial rigidity statement about the fixed set graph $\Theta_\Gamma$. This is recorded in three theorems which provide a more general framework:
\begin{enumerate}
\item (Theorem~\ref{theo:main}) Measure equivalent $2$-dimensional Artin group of hyperbolic type must have isomorphic fixed set graphs.
\item (Theorem~\ref{theointro:main-furman}) If every automorphism of $\Theta_\Gamma$ is induced by an automorphism of $\Cay_\Gamma$ (a form of rigidiy analogous to Tits' theorem for buildings \cite{tits2009buildings}, as will be explained in a more detailed way below), then every group which is measure equivalent to $G_\Gamma$ is commensurable to a lattice in $\Aut(\Cay_\Gamma)$ (as in Furman's theorem \cite{Fur-me}).
\item (Theorem~\ref{theointro:main-kida}) If (up to finite-index) every automorphism of $\Theta_\Gamma$ comes from the action of an element of $G_\Gamma$ (as in Ivanov's theorem about the curve graph \cite{Iva}), then $G_\Gamma$ satisfies the strongest form of measure equivalence rigidity (as in Kida's theorem \cite{Kid}).   
\end{enumerate}

Our first theorem is a measure group theoretic analogue of a theorem of Osajda and the second-named author, stating that every quasi-isometry between two $2$-dimensional Artin groups of hyperbolic type induces an isomorphism between their fixed set graphs \cite[Theorem~1.3]{HO}. Given two graphs $\Theta_1$ and $\Theta_2$, the set $\Isom(\Theta_1\to\Theta_2)$ of all graph isomorphisms from $\Theta_1$ to $\Theta_2$ is equipped with the compact-open topology. Work of Alvarez and Gaboriau regarding the behavior of measure equivalence with respect to free products \cite{AG} allows us to restrict the study to the case where $G_\Gamma$ is freely indecomposable, i.e.\ $\Gamma$ is connected.

\begin{theointro}[Theorem~\ref{theo:main-artin}]\label{theo:main}
Let $G_1=G_{\Gamma_1}$ and $G_2=G_{\Gamma_2}$ be two $2$-dimensional Artin groups of hyperbolic type with connected defining graphs. If $G_1$ and $G_2$ are measure equivalent, then $\Theta_{\Gamma_1}$ and $\Theta_{\Gamma_2}$ are isomorphic. In addition, for every measure equivalence coupling $\Sigma$ between $G_1$ and $G_2$, there is an almost $(G_1\times G_2)$-equivariant Borel map $$\Sigma\to\Isom(\Theta_{\Gamma_1}\to\Theta_{\Gamma_2}),$$ where the $(G_1\times G_2)$-action on $\Isom(\Theta_{\Gamma_1}\to\Theta_{\Gamma_2})$ is via $(g_1,g_2)\cdot f (v):=g_2 f(g_1^{-1}v)$.
\end{theointro}    

Using standard techniques in measured group theory \cite{Fur-me,MS,Kid,Kid-amalgams,BFS}, strong rigidity results are often deduced from a rigidity statement about self couplings -- in analogy, studying the group of self quasi-isometries is often a crucial step to establish quasi-isometry rigidity statements. In this way, Theorem~\ref{theo:main} can be considered as the main theorem of the present work, from which all other theorems are deduced by specifying to the case where $\Gamma_1=\Gamma_2$, and exploiting various forms of rigidity of the fixed set graph that build on the work of Crisp \cite{crisp}. 

To motivate the statement of our next theorem, recall that when $X$ is a higher rank irreducible symmetric space of non-compact type, measurable automorphisms of the Tits boundary of $X$ are induced by isometries of $X$  \cite{tits2009buildings}. In our case, the Cayley complex $\Cay_\Gamma$ will play the role of $X$. Just like the Tits boundary encodes the intersection of flats in the associated symmetric space, the graph $\Theta_\Gamma$ encodes the intersection pattern of a canonical collection of flat planes in $\Cay_\Gamma$. It natural to ask whether automorphisms of $\Theta_\Gamma$ are induced by automorphisms of $\Cay_\Gamma$. For this purpose, one can always define a canonical \emph{comparison map} $\Aut(\Cay_\Gamma)\to\Aut(\Theta_\Gamma)$ (see Lemma~\ref{lemma:cay-to-theta}), which is actually an isomorphism in many cases \cite{crisp}. 

\begin{theointro}[Theorem~\ref{theo:main-furman}]\label{theointro:main-furman}
Let $G_\Gamma$ be a $2$-dimensional Artin group of hyperbolic type. Assume that $\Gamma$ is connected and has no valence $1$ vertex, and that the comparison map $\Aut(\Cay_\Gamma)\to\Aut(\Theta_\Gamma)$ is an isomorphism. Let $H$ be a countable group that is measure equivalent to $G_\Gamma$.

Then there is a homomorphism $H\to\Aut(\Cay_\Gamma)$ with finite kernel, and whose image is a lattice.
\end{theointro}

The assumptions of Theorem~\ref{theointro:main-furman} hold for all groups appearing in Theorem~\ref{theointro:sample}. With a bit more work, one could actually prove that our assumptions that $\Gamma$ is connected with no valence $1$ vertex are redundant with the comparison map being an isomorphism. On the other hand, these assumptions are natural to make, and exclude situations where $G_\Gamma$ splits as a free product or as an amalgam over $\mathbb{Z}$.

\begin{theointro}[Theorem~\ref{theo:main-kida}]\label{theointro:main-kida}
Let $G_\Gamma$ be a $2$-dimensional Artin group of hyperbolic type. Assume that $\Gamma$ is connected and that the natural map $G_\Gamma\to\Aut(\Theta_\Gamma)$ is injective and has finite-index image. Let $H$ be a countable group which is measure equivalent to $G_\Gamma$.

Then $H$ is almost isomorphic to $G_\Gamma$. 
\end{theointro}

The injectivity assumption is satisfied as long as $G_\Gamma$ is not commensurable to the direct sum of $\mathbb{Z}$ and a free group (i.e.\ as long as $\Gamma$ is neither an edge, nor a star with all edges labelled $2$). The finite-index image assumption is the most crucial. The assumptions of Theorem~\ref{theointro:main-kida} hold for all groups appearing in Theorem~\ref{theointro:sample}(2). 

\begin{rk} Theorems~\ref{theointro:sample} and \ref{theointro:main-furman} naturally raise the question of describing all lattices in  $\Aut(\Cay_\Gamma)$. In particular, are there situations where every lattice is in fact commensurable to $G_\Gamma$, without $\Aut(\Cay_\Gamma)$ being countable? A positive answer to this question would provide a new pathway towards the strongest possible form of rigidity. We can actually say a bit more if we additionally assume that $H$ is finitely presented. Indeed, as $G_\Gamma$ is $2$-dimensional, its Cayley complex is contractible \cite{CharneyDavis}, so a theorem of Gandini \cite[Theorem~3.2]{Gan} ensures that there is no finitely presented non-uniform lattice in $\Aut(\Cay_\Gamma)$. So if $G_\Gamma$ satisfies the assumptions of Theorem~\ref{theointro:main-furman} (which is the case, for instance, for all the groups appearing in Theorem~\ref{theointro:sample}), and if $H$ is finitely presented and measure equivalent to $G_\Gamma$, then $H$ must be commensurable to a uniform lattice in $\Aut(\Cay_\Gamma)$. In other words, every finitely presented group $H$ which is measure equivalent to an Artin group $G_\Gamma$ as in Theorem~\ref{theointro:main-furman}, is also quasi-isometric to $G_\Gamma$, and understanding these groups amounts to describing uniform lattices in $\Aut(\Cay_\Gamma)$.
\end{rk}

\paragraph*{Proof ingredients.} 

The central point of the paper is to prove Theorem~\ref{theo:main}: all other statements follow using standard techniques in measured group theory and combinatorial rigidity results regarding the fixed set graph.

The broad outline of our proof follows Kida's strategy in the mapping class group setting. In fact, it turns out that there are many striking analogies between mapping class groups and $2$-dimensional Artin groups of hyperbolic type. These analogies motivated us to provide in Section~\ref{sec:axioms} an axiomatic framework that covers both settings. In this introduction, we will focus on Artin groups, and highlight the most striking analogies with mapping class groups along the way. 

It is helpful to have in mind the group-theoretic analogue of our theorem, which consists in building a graph isomorphism between $\Theta_{\Gamma_1}$ and $\Theta_{\Gamma_2}$ from a group isomorphism (or perhaps a commensuration) between $G_{\Gamma_1}$ and $G_{\Gamma_2}$. Following a celebrated strategy of Ivanov \cite{Iva}, this can be done by providing a purely algebraic characterization of stabilizers of vertices and of adjacency in the fixed set graph. Our strategy should be understood as a groupoid-theoretic version of the above. Namely, starting only from a measure equivalence coupling between $G_{\Gamma_1}$ and $G_{\Gamma_2}$, a standard construction originating in work of Furman \cite{Fur-oe} yields a probability space equipped with measure-preserving actions of $G_{\Gamma_1}$ and $G_{\Gamma_2}$, which are orbit equivalent in restriction to some positive measure subset. In other words, we have a measured groupoid with two cocycles towards $G_{\Gamma_1}$ and $G_{\Gamma_2}$. Our task is then to give a groupoid-theoretic characterization of subgroupoids that ``fix'' (in an appropriate sense) a vertex of $\Theta_\Gamma$, and a groupoid-theoretic characterization of adjacency. In order to keep this introduction not too technical, we will only provide the group-theoretic analogues of these statements. 

The group-theoretic version of our characterization of vertex stabilizers is the following (here $G=G_\Gamma$ is a $2$-dimensional Artin group of hyperbolic type). 
\medskip

\emph{A subgroup $H\subseteq G$ is equal to the stabilizer in $G$ of some vertex of $\Theta_\Gamma$ if and only if $H$ is a maximal (for inclusion) subgroup of $G$ satisfying the following property: $H$ is nonamenable and has a normal infinite amenable subgroup.}
\medskip

Interestingly, the exact same characterization holds for stabilizers of essential simple closed curves in mapping class groups of finite-type surfaces, the infinite normal amenable subgroup being the subgroup generated by the Dehn twist about the curve. The group-theoretic version of our characterization of adjacency is the following -- notice again that the same characterization holds for the curve graph. 
\medskip

\emph{Two distinct vertices $v_1$ and $v_2$ of $\Theta_\Gamma$ are adjacent if and only if the intersection of their stabilizers in $G$ does not fix any third vertex of $\Theta_\Gamma$.}
\medskip

Among the geometric ingredients used in the proof of Theorem~\ref{theo:main}, we mention the amenability of the $G$-action on the boundary $\partial_\infty \cod_\Gamma$ of its coned-off Deligne complex $\cod_\Gamma$ (guaranteed by Theorem~\ref{theo:boundary-amenable-intro}), the partition of the horocompactification of $\cod_\Gamma$ already mentioned in the proof sketch of boundary amenability, and the existence of a barycenter map -- which easily follows from $\cat$ geometry -- that associates a point in $\cod_\Gamma$ to every triple of pairwise distinct points in $\partial_\infty \cod_\Gamma$. Also, it is important to have an analogue to the notion of a \emph{canonical reduction system} from the mapping class group setting: this comes in the form of a map that canonically associates to every (possibly infinite) collection of vertices of $\Theta_\Gamma$ with infinite elementwise stabilizer, a unique vertex of $\Theta_\Gamma$ -- this construction, and its relationship to canonical reduction systems, is explained in detail in Section~\ref{sec:canonical-vertices}. A reader familiar with Kida's proof of measure equivalence rigidity in \cite{Kid} will have recognized analogues of the essential ingredients from there.

\subsection{Orbit equivalence rigidity, $W^*$-rigidity, lattice embeddings}

\paragraph*{Orbit equivalence rigidity.} Measure equivalence rigidity is intimately related to orbit equivalence rigidity of ergodic probability measure-preserving actions. We record its consequences when we have the strongest form of measure equivalence rigidity, as in Theorem~\ref{theointro:main-kida}. We refer to Section~\ref{sec:oe} for the relevant terminology and variations. The following theorem applies in particular to all Artin groups as in Theorem~\ref{theointro:sample}(2).

\begin{theointro}[Theorem~\ref{theo:oe-superrigid}]\label{theo:intro-oe}
Let $G=G_\Gamma$ be as in Theorem~\ref{theointro:main-kida}. Let $H$ be a countable group. Suppose 
$G\actson (X,\mu)$ and $H\actson (Y,\nu)$ are two free ergodic measure-preserving  actions on standard probability spaces by Borel automorphisms.

If the actions $G\actson (X,\mu)$ and $H\actson (Y,\nu)$ are stably orbit equivalent, then they are virtually conjugate.
\end{theointro}

\paragraph*{$W^*$-rigidity.} We also obtain strong rigidity results for the cross-product von Neumann algebras associated to actions as above. Recall that a group $G$ is \emph{Cartan rigid} \cite{PV2} if for every free, ergodic, probability measure-preserving action of $G$ on a standard probability space $X$, the cross-product von Neumann algebra $L^\infty(X)\rtimes G$ (obtained by the construction of Murray and von Neumann \cite{MvN}) contains $L^\infty(X)$ as its unique \emph{Cartan subalgebra}, i.e.\ maximal abelian subalgebra whose normalizer generates $L^\infty(X)\rtimes G$, up to unitary conjugacy. Via a theorem of Singer \cite{Sin}, Cartan rigidity implies that the von Neumann algebra $L^\infty(X)\rtimes G$ recovers the action $G\actson X$ up to orbit equivalence. When combined with orbit equivalence rigidity as in Theorem~\ref{theo:intro-oe}, this yields a strong rigidity result for $L^\infty(X)\rtimes G$. The following theorem yields a new infinite class of examples of countable groups, all of whose free ergodic probability measure-preserving actions satisfy a strong form of $W^*$-superrigidity, after the groundbreaking work of Popa and Vaes \cite{PV2} (see also \cite{HPV}), and work of Chifan, Ioana and Kida \cite{CIK} in mapping class group settings.

\begin{theointro}[Theorem~\ref{theo:cartan-rigid} and Corollary~\ref{cor:von-neumann-strong}]\label{theo:intro-strong-von-neumann}
Let $G=G_\Gamma$ be a $2$-dimensional Artin group of hyperbolic type. 
\begin{enumerate}
\item Assume that $\Gamma$ is not a clique and that $G_\Gamma$ is not isomorphic to the direct sum of $\mathbb{Z}$ and a free group. Then $G_\Gamma$ is Cartan rigid. 
\item Assume additionally that $G_\Gamma$ is as in Theorem~3(2), or more generally as in Theorem~\ref{theointro:main-kida}. Let $H$ be a countable group. Let $G\actson X$ and $H\actson Y$ be free, ergodic, probability measure-preserving actions on standard probability spaces. If $L^\infty(X)\rtimes G$ is isomorphic to $L^\infty(Y)\rtimes H$, then the actions $G\actson X$ and $H\actson Y$ are virtually conjugate.
\end{enumerate}
\end{theointro}
 
Cartan rigidity is derived from a theorem of Ioana \cite[Theorem~1.1]{Ioa2}, combined with the existence of splittings of $G_\Gamma$ as an amalgamated free product that are readable from the defining graph. For the second part, when the $G$-action on $X$ is additionally assumed \emph{aperiodic} (i.e.\ every finite index subgroup of $G$ acts ergodically on $X$), then we get the stronger conclusion that the actions $G\actson X$ and $H\actson Y$ are actually conjugate.

\paragraph*{Lattice embeddings.} Our next theorem describes all possible ways in which a $2$-dimensional Artin group of hyperbolic type that satisfies one the above measure equivalence rigidity statements can embed as a lattice in a locally compact second countable group.

\begin{theointro}[Theorems~\ref{theo:lattice-embeddings} and~\ref{theo:lattice-embeddings-2}]\label{theointro:lattice-embeddings}
	Let $f:G_\Gamma\to \mathsf{H}$ be a lattice embedding from an Artin group in a locally compact second countable group $\mathsf{H}$. Then 
	\begin{enumerate}
		\item if $G_\Gamma$ satisfies the assumptions of Theorem~\ref{theointro:main-furman}, then there is a continuous homomorphism $g:\mathsf{H}\to\Aut(\Cay_\Gamma)$ with compact kernel such that $g\circ f$ coincides with the natural map $G_\Gamma\to\Aut(\Cay_\Gamma)$;
		\item if $G_\Gamma$ satisfies the assumptions of Theorem~\ref{theointro:main-kida}, then there is a compact normal subgroup $K\unlhd \mathsf{H}$ such that $f$ induces a finite-index inclusion $G_\Gamma\to \mathsf{H}/K$. 
	\end{enumerate}
\end{theointro}

\paragraph*{Cocycle superrigidity.} Our last theorem is a cocycle rigidity statement: this is not a direct consequence of measure equivalence rigidity, but it follows from the geometric techniques developed in the present work together with work of Guirardel, Lécureux and the first named author \cite{GHL}. It extends the fact that there is no nontrivial homomorphism from a higher-rank lattice to a $2$-dimensional Artin group of hyperbolic type.

\begin{theointro}[Theorem~\ref{theo:cocycle}]\label{theointro:cocycle}
Let $G$ be a $2$-dimensional Artin group of hyperbolic type. Let $\mathsf{H}$ be a product of connected higher rank simple algebraic groups over local fields. Let $\Lambda$ be either $\mathsf{H}$, or a lattice in $\mathsf{H}$. Let $\Omega$ be a standard probability space equipped with an ergodic measure-preserving $\Lambda$-action by Borel automorphisms.

Then every Borel cocycle $\Lambda\times\Omega\to G$ is cohomologically trivial.
\end{theointro}

\subsection{Organization of the paper} 

Section~\ref{sec:background} gives general background on Artin groups and their actions on certain geometric spaces, in particular the modified Deligne complex and its coned-off version. 

In Section~\ref{sec:geometry-cat}, we establish several geometric features of $\cat$ spaces of crucial importance in the paper: in particular, we study the horofunction compactification of such spaces, and we describe a barycenter map that associates a center to every triple of points in the visual boundary. 

In Section~\ref{sec:boundary-amenability}, we establish our first main theorem, concerning boundary amenability of groups acting on piecewise hyperbolic $\cat$ simplicial complexes with finitely many isometry types of simplices. 

 Section~\ref{sec:background-me} reviews background on measured groupoids and their relationship to measure equivalence,  as well as how to derive measure equivalence rigidity results from a rigidity statement about self-couplings. 

In Section~\ref{sec:axioms}, we give a general axiomatic framework, involving an action of a group $G$ on a compact space and on a graph which plays the role of a curve graph for $G$, for proving measure equivalence rigidity statements. Our framework recovers Kida's work on mapping class groups, and will turn out to also be adapted to the Artin group setting. 

In Section~\ref{sec:theta}, we introduce the fixed set graph $\Theta_\Gamma$, and thoroughly study its vertex and edge stabilizers in order to check the axiomatic framework from the previous section in the context of $2$-dimensional Artin groups of hyperbolic type. 

Section~\ref{sec:combinatorial-rigidity} contains combinatorial rigidity statements regarding the automorphism groups of the Cayley complex, the modified Deligne complex and the fixed set graph, essentially coming from the work of Crisp \cite{crisp}.

Section~\ref{sec:icc} establishes the \emph{strong ICC} property for the automorphism group of the Deligne complex: this refines the fact that, under suitable assumptions on $\Gamma$, all nontrivial conjugacy classes of $G_{\Gamma}$ are infinite, and turns out to be useful for our arguments.

In Section~\ref{sec:theorems}, we complete the proof of all our main measure equivalence rigidity theorems, and their applications to orbit equivalence rigidity and lattice embeddings. The theorems related to von Neumann algebras are established in Section~\ref{sec:von-neumann}, and Theorem~\ref{theointro:cocycle} is proved in Section~\ref{sec:cocycle}. 

\paragraph*{Acknowledgments.}  This work was launched during the conference \emph{Non-positive curvature} held in May 2019 in Warsaw, as part of the Semester \emph{Geometric and Analytic Group Theory},  partially supported by the grant 346300 for IMPAN from the Simons Foundation and the matching 2015-2019 Polish MNiSW fund. We are grateful to the organizers of these events.
The first named author also acknowledges support from the Agence Nationale de la Recherche under Grant ANR-16-CE40-0006 DAGGER.

\section{Background on Artin groups of hyperbolic type}\label{sec:background}

\emph{This section contains background material on $2$-dimensional Artin groups of hyperbolic type, including basic definitions, the definition by Charney and Davis of the modified Deligne complex, and the definition by Martin and Przytycki of its coned-off version, a $\cat$ space on which the group acts nonelementarily by isometries.}

\subsection{Artin groups}

Let $\Gamma$ be a finite simple labeled graph, where every edge is labeled by an integer at least equal to $2$. We denote by $V\Gamma$ the vertex set of $\Gamma$. The \emph{Artin group with defining labeled graph $\Gamma$}, denoted $G_{\Gamma}$, is the group defined by the following presentation: it is generated by $V\Gamma$, with one relation of the form $$\underbrace{aba\cdots}_{m}=\underbrace{bab\cdots}_{m}$$ for each $a$ and $b$ spanning an edge labeled by $m$. The elements in $V\Gamma$ are called \emph{standard generators} of $G_\Gamma$. For $s,t\in V\Gamma$, let $m_{st}$ be equal to either $+\infty$ if $s$ and $t$ are not connected by an edge in $\Gamma$, or to the label of the edge if $s$ and $t$ are connected by an edge.

\begin{rk}
We warn the reader that two distinct labeled graphs $\Gamma$ and $\Gamma'$ may define isomorphic Artin groups $G_\Gamma$ and $G_{\Gamma'}$, see \cite[Theorem~1]{crisp}. In this paper, when we write ``Let $G=G_\Gamma$ be an Artin group with defining graph $\Gamma$'', we assume that $\Gamma$ is fixed once and for all as part of the data.   
\end{rk}

\begin{lemma}[{van der Lek \cite[Theorem~II.4.13]{Van1983homotopy}}]
	\label{lem:injective}
Let $\Gamma$ be a finite simple labeled graph, and let $\Gamma'\subseteq\Gamma$ be an induced subgraph -- i.e.\ an edge of $\Gamma$ is in $\Gamma'$ if its endpoints are in $\Gamma'$, and $\Gamma'$ is equipped with the induced labeling of $\Gamma$. Then the natural map $G_{\Gamma'}\to G_{\Gamma}$ is injective. Moreover, for induced subgraphs $\Gamma_1$ and $\Gamma_2$ of $\Gamma$, we have $G_{\Gamma_1}\cap G_{\Gamma_2}=G_{\Gamma_1\cap\Gamma_2}$ (when viewing them as subgroups of $G_\Gamma$).
\end{lemma}

By \cite[Theorem B]{CharneyDavis}, an Artin group has cohomological dimension at most $2$ if for each triangle in $\Gamma$ with sides labeled by $m,n,r$, we have $\frac{1}{m}+\frac{1}{n}+\frac{1}{r}\le 1$. Being of finite cohomological dimension, these Artin groups are torsion-free. When we say that an Artin group is \emph{$2$-dimensional}, we mean that its cohomological dimension is equal to $2$.

\subsection{The modified Deligne complex}\label{sec:background-deligne}

The following definition is due to Charney and Davis \cite{CharneyDavis}.

\begin{de}[Modified Deligne complex]\label{de:deligne} 
Let $G=G_\Gamma$ be a $2$-dimensional Artin group with defining labeled graph $\Gamma$. The \emph{modified Deligne complex} $\del$ is the geometric realization of the poset (ordered by inclusion) of all subsets of $G_\Gamma$ of the form $gG_{\Gamma'}$, where $\Gamma'$ is either the empty subgraph (in which case $G_{\Gamma'}$ is the trivial subgroup), a vertex of $\Gamma$, or an edge of $\Gamma$.
\end{de}	

The \emph{rank} of a vertex $gG_{\Gamma'}$ is the number of vertices in $\Gamma'$.	
It is clear that $\mathbb D_\Gamma$ is a $2$--dimensional simplicial complex, and $G_\Gamma$ acts on $\mathbb D_\Gamma$ without inversions, i.e.\ if an element of $G_\Gamma$ preserves a simplex of $\mathbb D_\Gamma$, then it fixes the simplex pointwise. We endow $\mathbb D_\Gamma$ with a piecewise Euclidean metric, called the \emph{Moussong metric}, such that each triangle $$\Delta(g_1,g_2G_s,g_3G_{st})$$ is a Euclidean triangle with angle $\pi/2$ at $g_2G_s$ and angle $\frac{\pi}{2n}$ at $g_3G_{st}$ with $n$ being the label of the edge $st$ of $\Gamma$. By \cite[Proposition 4.4.5]{CharneyDavis}, the complex $\mathbb D_\Gamma$ is $\CAT$ with this metric. As was observed in \cite[Lemma~6]{crisp}, the $G_\Gamma$-action on $\mathbb D_\Gamma$ is semisimple, i.e.\ every element $g\in G_\Gamma$ achieves its translation length $|g|_{\del}:=\inf_{x\in\del}d(x,g\cdot x)$ at some point of $\del$.	

The link of a rank 0 vertex $g$ in $\del$ can be identified with the barycentric subdivision of $\Gamma$ by sending $gG_s$ to the vertex $s$ of $\Gamma$ and $gG_{st}$ to the midpoint of the edge $\overline{st}$ of $\Gamma$. The above $\mathrm{CAT}(0)$ metric on $\del$ induces an angular metric on the link of $g$ (see \cite[p.~103, Chapter I.7.15]{BH}). If the edge $\overline{st}$ of $\Gamma$ is labeled by $n$, then the length of $\overline{st}$ under the angular metric is $\pi-\frac{\pi}{n}$, which is always at least $\pi/2$.

Note that if $\Gamma'\subset\Gamma$ is a full subgraph, then we can identify $\mathbb D_{\Gamma'}$ as a subcomplex of $\mathbb D_\Gamma$. The embedding $\mathbb D_{\Gamma'}\to\mathbb D_\Gamma$ is isometric with respect to the Moussong metric on $\mathbb D_{\Gamma'}$ and $\mathbb{D}_\Gamma$, see \cite[Lemma~5.1]{Charney}. By uniqueness of geodesics in $\CAT$ space, we deduce that $\mathbb D_{\Gamma'}$ is a convex subcomplex of $\mathbb D_\Gamma$. 

One can also metrize $\del$ as a CAT(-1) space, namely $\Delta(g_1,g_2G_s,g_3G_{st})$ as above is metrized as a hyperbolic triangle with the same angle at $g_2 G_s$ and $g_3G_{st}$, but slightly smaller angle at $g_1$. The angle at $g_1$ is chosen to ensure that $\del$ is locally CAT(-1) around $g_1$. This turns $\del$ into a $\cat$ space with finitely many isometry types of vertices, on which $G_\Gamma$ acts by isometries. This is a simple case of a more delicate construction discussed in Martin and Przytycki \cite[Section 3.1]{MP}. Any such CAT(-1) metric will be called an \emph{adapted CAT(-1) metric} (they are not unique).

\begin{lemma}
	\label{lemma:aut deligne}
Let $\alpha$ be a simplicial automorphism of $\del$. Then
\begin{enumerate}
	\item for every $k\in\{0,1,2\}$, every rank $k$ vertex is sent by $\alpha$ to a rank $k$ vertex;
	\item $\alpha$ can be represented as an isometry with respect to either the Moussong metric or any adapted CAT(-1) metric.
\end{enumerate}
\end{lemma}
\begin{proof}
Note that the link of a rank $0$ vertex is a finite graph, the link of a rank 1 vertex is a join of two discrete sets with one being finite and another infinite and the link of a rank 2 vertex is either a join of two infinite discrete sets (if this vertex represents $gG_{e}$ such that $e$ is an edge of $\Gamma$ labeled by $2$) or is an infinite graph which does not split as a join (if this vertex represents $gG_{e}$ such that $e$ has label at least $3$, as follows from the proof of \cite[Proposition~40]{crisp}). Thus the first assertion follows.

It also follows from \cite[Lemma~39]{crisp} that the link of vertex $gG_{e}$ and $g'G_{e'}$ are isomorphic if and only if $e$ and $e'$ are labeled by the same number, thus the second assertion follows from the way we metrize $\del$.
\end{proof}

A 2-dimensional Artin group is of \emph{hyperbolic type} if its modified Deligne complex is Gromov-hyperbolic. By \cite[Lemma~5]{crisp}, an Artin group $G_\Gamma$ is a 2-dimensional Artin group of hyperbolic type if and only if each induced 4-cycle (i.e.\ 4-cycle without diagonals) in $\Gamma$ has an edge with label at least $3$ and for each triangle in $\Gamma$ labeled by $m,n,r$, we have $\frac{1}{m}+\frac{1}{n}+\frac{1}{r}<1$.

\section{Geometry of spaces of negative curvature}\label{sec:geometry-cat}

\emph{In this section, we establish a few facts about the geometry of a $\cat$ space $X$ that will be useful in later sections. The two main results are a construction of a canonical map from the horoboundary of $X$ to $X\cup\partial_\infty X$, and a construction of a canonical map that associates to every triple of pairwise distinct points of $\partial_\infty X$, a `barycenter' in $X$. (The latter will only be used in the proof of the theorems on measure equivalence, and not for boundary amenability.) Most technicalities in the present section come from the necessity to check that these maps can be made Borel.} 

\subsection{Centers of bounded subsets in $\CAT$ spaces}

Let $X$ be a complete $\CAT$ metric space. Let $\calb(X)$ be the space of all nonempty bounded closed subsets of $X$, equipped with the Hausdorff topology. Every set $B\in\calb(X)$ has a \emph{circumcenter}, denoted by $\cent(B)$, defined in the following way. The \emph{radius} of $B$, denoted by $\rad(B)$, is the smallest $r\ge 0$ such that $B$ is contained in a closed ball of radius $r$. The \emph{circumcenter} of $B$, denoted by $\cent(B)$, is then defined as the unique point $x\in X$ such that $B$ is entirely contained in the closed ball $\overline{B}(x,\rad(B))$, see \cite[Proposition~II.2.7]{BH}.

\begin{lemma}\label{lemma:approximate-center}
Let $X$ be a complete $\CAT$ space. For every $R>0$ and every $\epsilon>0$, there exists $\delta>0$ such that for every $B\in\calb(X)$ of radius at most $R$ and every $c\in X$, if $B$ is contained in the closed ball $\overline{B}(c,\rad(B)+\delta)$, then $d(c,\cent(B))\le \epsilon$.
\end{lemma}

\begin{proof}
The proof is an adaptation of \cite[Proposition~II.2.7]{BH}. We let $\mathbb{E}^2$ be the Euclidean plane, and fix $O\in\mathbb{E}^2$. For every $r\in\mathbb{R}_+^*$ and every $\delta>0$, we let $A(r,\delta)$ be the closed $\delta$-neighborhood of the circle of radius $r$ centered at $O$ in $\mathbb{E}^2$. We can (and shall) choose $\delta>0$ such that for every $r\le R$, every segment contained in $A(r,\delta)$ has length at most $\epsilon/2$.  

Let $B\in\calb(X)$ be a closed bounded set of radius at most $R$, and let $c\in X$ be such that $B\subseteq\overline{B}(c,\rad(B)+\delta)$. For every $x\in B$, we let $\overline{\Delta}_x=\Delta(O,\overline{c}_x,\overline{\cent(B)}_x)$ be a comparison triangle to the geodesic triangle $\Delta(x,c,\cent(B))$, and we let $m_x$ be the midpoint of the side of $\overline{\Delta}_x$ opposite to $O$. Notice that for every $x\in B$, we have $d(c,x)\le\rad(B)+\delta$ and $d(\cent(B),x)\le\rad(B)$. This implies that $m_x\in \overline{B}_{\mathbb{E}^2}(O,\rad(B)+\delta)$. 

We claim that there exists $x\in B$ such that $m_x$ belongs to $A(\rad(B),\delta)$. Indeed, otherwise, for every $x\in B$, we have $d_{\mathbb{E}^2}(O,m_x)\le \rad(B)-\delta$. Let $m$ be the midpoint of the geodesic segment from $c$ to $\cent(B)$ in $X$. By $\CAT$ comparison, for every $x\in B$, we have $d(x,m)\le d_{\mathbb{E}^2}(O,m_x)$ whence $B\subseteq B(m,\rad(B)-\delta)$. This contradicts the definition of $\rad(B)$, thus showing our claim.

We now fix $x\in B$ as given by the above paragraph. Then the side of $\overline{\Delta}_x$ opposite to $O$ is at least half-contained in $A(\rad(B),\delta)$, so our choice of $\delta$ ensures that it has length at most $\epsilon$. This implies that $d(c,\cent(B))\le \epsilon$.
\end{proof}

\begin{lemma}\label{lemma:center-countable-bounded-set}
Let $X$ be a complete $\CAT$ space. Let $B$ be a countable nonempty closed bounded subset of $X$, and let $(B_n)_{n\in\mathbb{N}}$ be an increasing sequence of nonempty finite subsets of $B$ such that $B=\overline{\cup_n B_n}$.

Then $(\rad(B_n))_{n\in\mathbb{N}}$ converges to $\rad(B)$ and $(\cent(B_n))_{n\in\mathbb{N}}$ converges to $\cent(B)$. 
\end{lemma}

\begin{proof}
The sequence $(\rad(B_n))_{n\in\mathbb{N}}$ is a non-decreasing sequence which is bounded from above by $\rad(B)$, so it converges to some $r_\infty>0$, with $\rad(B)\ge r_\infty$. 

We claim that the sequence $(\cent(B_n))_{n\in\mathbb{N}}$ is a Cauchy sequence. Indeed, let $\epsilon>0$. Let $\delta>0$ be given by Lemma~\ref{lemma:approximate-center} (for $R=\rad(B)$ and $\epsilon$), and let $n_0\in\mathbb{N}$ be such that $|r_\infty-\rad(B_{n_0})|<\delta$. Let $n,m\ge n_0$, with $m\ge n$. The set $B_n$ is contained in $B_m$, and therefore it is contained in $\overline{B}(\cent(B_m),\rad(B_n)+\delta)$. Lemma~\ref{lemma:approximate-center} thus implies that $d(\cent(B_n),\cent(B_m))\le\epsilon$. This shows that the sequence $(\cent(B_n))_{n\in\mathbb{N}}$ is Cauchy, and as $X$ is complete it converges to a point $c_\infty\in X$.

We will now prove that for every $\delta'>0$, the set $B$ is contained in $\overline{B}(c_\infty,r_\infty+\delta')$; the conclusion will then follow from Lemma~\ref{lemma:approximate-center}. Indeed, let $n_0\in\mathbb{N}$ be such that for all $n\ge n_0$, we have $d(c_\infty,\cent(B_n))\le\delta'$. Then for every $n\ge n_0$, we have $B_n\subseteq \overline{B}(c_\infty,\rad(B_n)+\delta')\subseteq \overline{B}(c_\infty,r_\infty+\delta')$. As $B$ is the closure of the union of all the $B_n$, the conclusion follows.
\end{proof}

\subsection{Horofunction compactification of a $\cat$ space}

Let $X$ be a separable metric space, and let $x_0\in X$ be a basepoint. We denote by $\mathcal{C}(X,\mathbb{R})$ the space of all real-valued continuous functions on $X$, equipped with the topology of convergence on compact subsets of $X$. The map
\begin{displaymath}
\begin{array}{cccc}
X & \to & \calc(X,\mathbb{R})\\
x & \mapsto & (h_x: z\mapsto d(z,x)-d(x_0,x))
\end{array}
\end{displaymath}
\noindent is a topological embedding, see e.g.\ \cite[Lemma~3.2]{MT}. The closure $\horox$ of the image of this embedding is called the \emph{horofunction compactification} of $X$, and it is compact and metrizable \cite[Proposition~3.1]{MT}.  Here we insist that in order for $\horox$ to be compact, it was important to equip $\calc(X,\mathbb{R})$ with the topology of uniform convergence on compact subsets of $X$, and not just uniform convergence on bounded subsets of $X$, which is another natural choice that has often been considered in the literature, e.g.\ \cite{BGS,Bal,BH}. The horofunction compactification does not depend on the choice of the basepoint $x_0$. If $G$ is a group of isometries of $X$, then the $G$-action extends to a continuous action by homeomorphisms on $\horox$ by letting $g\cdot h(x):=h(g^{-1}x)-h(g^{-1}x_0)$, see \cite[Lemma~3.4]{MT}.

 In the present paper, whenever $X$ is a $\mathrm{CAT}(0)$ space, we let $\partial_\infty X$ be the visual boundary of $X$, equipped with the cone topology, see \cite[p.~263, II.8.5]{BH}. In the special case when $X$ is Gromov hyperbolic, $\partial_\infty X$ can also be identified with the Gromov boundary of $X$. The goal of the present section is to establish the following statement, which combines earlier works on the horocompactification of a negatively curved space, most notably \cite{BGS} (here we follow the account given in \cite{BH}) and more recent work of Maher and Tiozzo \cite[Section~3.2]{MT}. 
 
\begin{prop}\label{prop:horo-partition}
Let $X$ be a complete separable $\cat$ space, and let $G$ be a countable group acting discretely on $X$ by isometries. The sets $\horox_{\infty}$ and $\horox_{\bdd}$ made of horofunctions that are unbounded (resp.\ bounded) from below form a Borel $G$-invariant partition of $X$, and there exist a $G$-equivariant homeomorphism $\theta_\infty:\horox_{\infty}\to\partial_\infty X$ and a $G$-equivariant Borel map $\theta_{\bdd}:\horox_{\bdd}\to X$, with the following properties: 
\begin{enumerate}
\item $\theta_\infty$ continuously extends the identity map on $X$,
\item there exists $M>0$ such that for every $h\in\horox_{\bdd}$ and every $x\in X$, if $|h(x)-\inf(h)|\le 1$, then $d(x,\theta_{\bdd}(h))\le M$;
\item for every sequence $(h_n)_{n\in\mathbb{N}}$ of horofunctions in $\horox_{\bdd}$, if $(\theta_{\bdd}(h_n))_{n\in\mathbb{N}}$ converges to a point $\xi\in\partial_\infty X$ for the visual topology, then $(h_n)_{n\in\mathbb{N}}$ converges to $\theta_{\infty}^{-1}(\xi)$. 
\end{enumerate}
\end{prop} 

\begin{proof}
Let $x_0\in X$ be a basepoint. 
Let $X_\mathbb{Q}$ be a dense countable $G$-invariant subset of $X$. Since $\cald$ is countable, the map 
\begin{displaymath}
\begin{array}{cccc}
\beta: & \horox & \to & \mathbb{R}\cup\{-\infty\}\\
& h & \mapsto & \inf\limits_{d\in X_\mathbb{Q}} h(d)
\end{array}
\end{displaymath}
\noindent is Borel. As $h\in\horox_\infty$ if and only if $\beta(h)=-\infty$, we deduce that $\horox_\bdd$ and $\horox_\infty$ are Borel subsets of $\horox$.

In the next four paragraphs, we will explain how to build the map $\theta_\infty$ continuously extending the identity map on $X$. As $X$ is $\CAT$, 
every point $\xi\in\partial_\infty X$ is represented by a unique geodesic ray $r_\xi$ originating at $x_0$. For every $x\in X$, the limit $$\lim_{t\to +\infty} h_{r_\xi(t)}(x)=\lim_{t\to +\infty} (d(r_\xi(t),x)-t)$$ exists, and the convergence is uniform on compact subsets of $X$ (in fact on bounded subsets of $X$, see \cite[Lemma~II.8.18.(2)]{BH}). In the limit, we get a horofunction $h_\xi$, which belongs to $\horox_\infty$ since it is unbounded from below on the ray $r_\xi$. This yields a map $\psi:\partial_\infty X\to\horox_\infty$ (sending $\xi$ to $h_\xi$), which is continuous (see \cite[Theorem~II.8.13]{BH}). In addition $\psi$ satisfies the following property, that we record for future use (see \cite[Theorem~II.8.13]{BH}):
\begin{enumerate}
\item[$(\ast)$] For every sequence $(y_n)_{n\in\mathbb{N}}\in X^{\mathbb{N}}$ converging (in the visual topology) to a point $\xi\in\partial_\infty X$, the sequence $(h_{y_n})_{n\in\mathbb{N}}$ converges to $\psi(\xi)=h_\xi$ uniformly on compact sets.
\end{enumerate}

We will now provide an inverse to the map $\psi$, coming from the work of Maher and Tiozzo: by \cite[Lemma~3.10]{MT}, there exists a unique $G$-equivariant map $\varphi:\horox_\infty\to\partial_\infty X$, which has the property that for every $h\in\horox_\infty$ and every $(y_n)_{n\in\mathbb{N}}\in X^{\mathbb{N}}$, if $h(y_n)$ converges to $-\infty$, then $(y_n)_{n\in\mathbb{N}}$ converges to $\varphi(h)$ in the visual topology. By \cite[Proposition~3.14]{MT}, the map $\varphi$ is continuous. We will now prove that $\varphi$ and $\psi$ are inverse to each other: this will establish the first part of the proposition, by letting $\theta_\infty=\varphi$.   

First, for every $\xi\in\partial_\infty X$, the sequence $(r_\xi(n))_{n\in\mathbb{N}}$ converges to $\xi$ in the visual topology, and $h_\xi(r_\xi(n))$ converges to $-\infty$, so the definition of $\varphi$ ensures that $\xi=\varphi(h_\xi)$. This shows that $\varphi\circ\psi=\mathrm{id}$.

To see that $\psi\circ\varphi=\mathrm{id}$, let $h\in\horox_\infty$, and write $h$ as the limit of maps  $h_{y_n}$ with $y_n\in X$. By \cite[Proposition~3.14]{MT}, the sequence $(y_n)_{n\in\mathbb{N}}$ converges to a point $\xi=\varphi(h)$ in $\partial_\infty X$ in the visual topology, and it follows from Property~$(\ast)$ above that $h=h_\xi$. Thus $\psi\circ\varphi(h)=\psi(\xi)=h$. This concludes our proof of the first part of the proposition.

In the next four paragraphs, we will now connstruct a $G$-equivariant Borel map $\theta_{\bdd}:\horox_{\bdd}\to X$ and check that it satisfies the second conclusion from the proposition.  We denote by $\cald$ the $G$-orbit of $x_0$, a countable $G$-invariant subset of $X$.

Let $\calp_f(\cald)$ be the countable set of all finite  (possibly empty) subsets of $\cald$, equipped with the discrete topology. Let $\{d_i\}_{i\in\mathbb{N}}$ be an enumeration of $\cald$. For every $n\in\mathbb{N}$, let $\alpha_n:\horox_\bdd\to\calp_f(\cald)$ be the map which associates to every $h\in\horox_\bdd$ the finite (possibly empty) set made of all $d\in\{d_1,\dots,d_n\}$ satisfying $h(d)\le\inf_{d'\in\cald_0} h(d')+1$. 

We claim that for every $n\in\mathbb{N}$, the map $\alpha_n$ is Borel. Indeed, to check this, it is enough to notice that for every $d$ in the countable set $\cald$, the subset $K_{n,d}\subseteq\horox_{\bdd}$ made of all horofunctions $h$ such that $d\in\alpha_n(h)$ is Borel. If $d\notin\{d_1,\dots,d_n\}$, then this subset is empty so the conclusion is obvious. If $d\in\{d_1,\dots,d_n\}$, then $h\in K_{n,d}$ if and only if $h(d)\le\inf_{d'\in\cald_0}h(d')+1$. As the maps $h\mapsto h(d)$ and $h\mapsto\inf_{d'\in\cald_0}h(d')$ are Borel, the conclusion follows.

For every $h\in\horox_{\bdd}$, let $\alpha_{\infty}(h)$ be the set of all $d\in\cald$ such that $h(d)\le\inf_{\cald} h+1$. Equivalently $\alpha_\infty(h)$ is equal to the increasing union of all $\alpha_n(h)$ as $n$ ranges over $\mathbb{N}$. For every $h\in\horox_{\bdd}$, the set $\alpha_\infty(h)$ is a nonempty closed subset of $X$, and by \cite[Lemma~3.13]{MT} it has bounded diameter.

Now, for every $n\in\mathbb{N}$, let $\theta_n:\horox_\bdd\to X$ be the Borel map defined by letting $\theta_n(h):=\cent(\alpha_n(h))$ if $\alpha_n(h)\neq\emptyset$, and $\theta_n(h)=x_0$ otherwise. It follows from Lemma~\ref{lemma:center-countable-bounded-set} that for every $h\in\horox_\bdd$, the points $\theta_n(h)$ converge to a point $\theta(h)$, with $\theta(h)=\cent(\alpha_\infty(h))$. In particular $\theta$ can be viewed as a map $\horox_\bdd\to X$ which is Borel, being a pointwise limit of Borel maps. In addition, it is $G$-equivariant in view of the definition of the $G$-action on the space of horofunctions and the fact that $\cald$ is $G$-invariant. We are thus done by letting $\theta_{\bdd}=\theta$ (using the boundedness of $\alpha_\infty(h)$, it satisfies the second conclusion of the proposition).

We finally prove the last part of the proposition, namely: if $(h_n)_{n\in\mathbb{N}}$ is a sequence of horofunctions in $\horox_{\bdd}$, and if the sequence $(\theta_{\bdd}(h_n))_{n\in\mathbb{N}}$ converges to a point $\xi\in\partial_\infty X$ in the visual topology, then the horofunctions $h_n$ converge to $h_\xi$ (uniformly on compact sets). It is in fact enough to prove that every accumulation point of the sequence $(h_n)_{n\in\mathbb{N}}$ belongs to $\horox_\infty$: indeed, from this, if $h$ is an accumulation point of $(h_n)_{n\in\mathbb{N}}$, then \cite[Proposition~3.14]{MT} implies that $\theta_{\bdd}(h_n)$ converges to $\theta_\infty(h)$, showing that $h=h_\xi$.

To simplify notations, for every $n\in\mathbb{N}$, we let $y_n=\theta_{\bdd}(h_n)$. By hyperbolicity (and using the Morse lemma), there exists a constant $C>0$, only depending on the hyperbolicity constant of $X$, such that for every $n\in\mathbb{N}$, there exists $\theta(n)\in\mathbb{N}$ such that for every $m\ge \theta(n)$, the point $r_\xi(n)$ is at distance at most $C$ from any geodesic segment joining  $x_0=r_\xi(0)$ to $y_m$. Using \cite[Proposition~3.6]{MT} and the fact that $h_m$ is close to its infimum at $y_m$, one then derives that for all sufficiently large $m\in\mathbb{N}$, the value of $h_m(r_\xi(n))$ is equal to $-d(x_0,r_\xi(n))$ up to a bounded error. In the limit, we see that every accumulation point of the sequence $(h_n)_{n\in\mathbb{N}}$ is a horofunction which is unbounded from below, which concludes our proof.   
\end{proof}

Given a loxodromic isometry $g$ of $X$, we let $g^{-\infty}$ and $g^{+\infty}$ the repelling and attracting fixed points of $g$ in $\partial_\infty X$, also viewed as a subset of the horocompactification of $X$ through Proposition~\ref{prop:horo-partition}. We record the following consequence of Proposition~\ref{prop:horo-partition} for future use.  

\begin{lemma}\label{lemma:ns-dynamics-on-horo}
Let $X$ be a complete separable $\cat$ space, let $g$ be a loxodromic isometry of $X$, and let $h\in\horox\setminus\{g^{-\infty}\}$. Then $g^n\cdot h$ converges to $g^{+\infty}$ as $n$ goes to $+\infty$.
\end{lemma}

\begin{proof}
Let $\horox=\horox_\infty\dunion\horox_\bdd$ be the partition provided by Proposition~\ref{prop:horo-partition}. The result is clear if $h\in\horox_\infty$ (using the $G$-equivariant homeomorphism with $\partial_\infty X$), so we assume that $h\in\horox_{\bdd}$. Let $K=\{x\in X|h(x)\le \inf h+1\}$ be the \emph{almost minimizing set} of $h$: by \cite[Lemma~3.13]{MT}, this is a bounded nonempty subset of $X$. For every $n\in\mathbb{N}$, the almost minimizing set of the horofuntion $$g^n\cdot h:x\mapsto h(g^{-n}x)-h(g^{-n}x_0)$$ is  $g^nK$, and one derives that $\theta_{\bdd}(g^n\cdot h)$ converges to $h$ (where $\theta_{\bdd}$ is as in Proposition~\ref{prop:horo-partition}). It follows from  the last conclusion of Proposition~\ref{prop:horo-partition} that $g^n\cdot h$ converges to $g^{+\infty}$. 
\end{proof}

\subsection{A barycenter map}

We denote by $(\partial_\infty X)^{(3)}$ the subspace of $(\partial_\infty X)^3$ made of pairwise distinct triples.

\begin{prop}\label{prop:barycenter}
Let $X$ be a separable $\cat$ space, and let $G$ be a countable group acting on $X$ by isometries. 

Then there exists a $G$-equivariant Borel map $(\partial_\infty X)^{(3)}\to X$.
\end{prop}

\begin{proof}
Since $X$ is $\cat$, we can find
\begin{enumerate}
\item a constant $K_1>0$ such that for every triple $(\xi_1,\xi_2,\xi_3)$ of pairwise distinct boundary points, denoting by $B(\xi_1,\xi_2,\xi_3)$ the set of all points $x\in X$ such that for every $i\in\{1,2,3\}$, the point $x$ is at distance at most $K_1$ from a point on the geodesic from $\xi_{i}$ to $\xi_{i+1}$ (with indices taken modulo $3$), the set $B(\xi_1,\xi_2,\xi_3)$ is nonempty and of bounded diameter;
\item a constant $K_2>0$ such that for every triple $(x,y,z)\in X^3$, if $d(x,y)+d(y,z)\le d(x,z)+1$, then $y$ lies at distance at most $K_2$ from the geodesic from $x$ to $z$.
\end{enumerate}

The space $X\cup\partial_\infty X$ is metrizable \cite[Section~5]{Vai}; we let $d_{X\cup\partial_\infty X}$ be a distance on $X\cup\partial_\infty X$ which is compatible with the topology. Let $X_\mathbb{Q}$ be a dense countable $G$-invariant subset of $X$, let $\calp_f(X_\mathbb{Q})$ be the countable collection of all finite subsets of $X_\mathbb{Q}$ (including the empty set), endowed with the discrete topology, and let $\{x_i\}_{i\in\mathbb{N}}$ be an enumeration of $X_\mathbb{Q}$. 

For every $n\in\mathbb{N}$, let $\alpha_n:(\partial_\infty X)^{(3)}\to \calp_f(X_\mathbb{Q})$ be the map that sends $(\xi_1,\xi_2,\xi_3)$ to the set of all $x\in\{x_1,\dots,x_n\}$ such that for every $k\in\mathbb{N}$, there exists a tuple $(y_1,y_2,y_3,x_{23},x_{31})\in X_{\mathbb{Q}}^5$ such that, letting $x_{12}=x$, the following three facts hold:
\begin{enumerate}
\item for every $i\in\{1,2,3\}$, we have $d_{X\cup\partial_\infty X}(y_i,\xi_i)\le\frac{1}{k}$, and
\item for every $i\in\{1,2,3\}$ (taken modulo $3$), we have $d(y_i,x_{i,i+1})+d(x_{i,i+1},y_{i+1})\le d(y_i,y_{i+1})+1$, and
\item for every $i,j\in\{1,2,3\}$ (taken modulo $3$), we have $d(x_{i,i+1},x_{j,j+1})\le 2K_1+1$. 
\end{enumerate}

As $\alpha_n(\xi_1,\xi_2,\xi_2)$ is defined by using continuous functions, for a given $x\in X_{\mathbb Q}$, the subset $K_{x,n}$ made of elements of $(\partial_\infty X)^{(3)}$ whose $\alpha_n$-image contains $x$ is Borel. As $\calp_f(X_\mathbb{Q})$ is countable, each $\alpha_n$ is Borel. 

\begin{figure}
\centering
\includegraphics[scale=0.55]{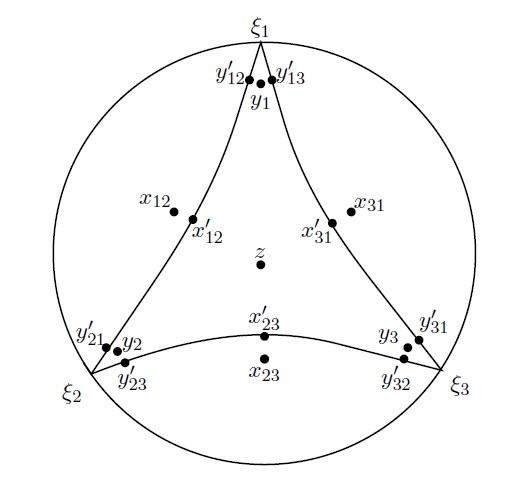}
\caption{Proof of the existence of the barycenter map.}
\label{fig:barycenter}
\end{figure}

We claim that for every $(\xi_1,\xi_2,\xi_3)\in (\partial_\infty X)^{(3)}$, there exists $n\in\mathbb{N}$ such that $\alpha_n(\xi_1,\xi_2,\xi_3)$ is nonempty -- the argument is illustrated in Figure~\ref{fig:barycenter}. Indeed, by definition of the constant $K_1$, we can find a point $z\in X$ such that for every $i\in\{1,2,3\}$, there exists a point $x'_{i,i+1}$ on the geodesic joining $\xi_i$ to $\xi_{i+1}$ with $d(z,x'_{i,i+1})\le K_1$. Since $X_\mathbb{Q}$ is dense, for every $i\in\{1,2,3\}$, we can find a point $x_{i,i+1}\in X_\mathbb{Q}$ with $d(x_{i,i+1},x'_{i,i+1})\le \frac{1}{6}$. Let now $n\in\mathbb{N}$ be such that $x_{1,2}\in\{x_1,\dots,x_n\}$; we will check that $\alpha_n(\xi_1,\xi_2,\xi_3)\neq\emptyset$. Let $k\in\mathbb{N}$. For every $i\in\{1,2,3\}$, we can find a point $y_i\in X_\mathbb{Q}$ with $d_{X\cup\partial_\infty X}(y_i,\xi_i)\le\frac{1}{k}$ and lying at distance at most $\frac{1}{6}$ (in $X$) both from a point $y'_{i,i-1}$ on the geodesic from $\xi_i$ to $\xi_{i-1}$ and from a point $y'_{i,i+1}$ on the geodesic from $\xi_i$ to $\xi_{i+1}$, so that $x'_{i,i+1}$ lies on the geodesic segment joining $y'_{i+1,i}$ to $y'_{i,i+1}$. It follows that 
\begin{displaymath}
\begin{array}{rl}
d(y_i,x_{i,i+1})+d(x_{i,i+1},y_{i+1}) & \le d(y_i,y'_{i,i+1})+d(y'_{i,i+1},x'_{i,i+1})+2d(x'_{i,i+1},x_{i,i+1})\\
& ~~+d(x'_{i,i+1},y'_{i+1,i})+d(y'_{i+1,i},y_{i+1})\\
& \le d(y'_{i,i+1},y'_{i+1,i})+\frac{4}{6}\\
& \le d(y_i,y_{i+1})+1.
\end{array}
\end{displaymath}
Finally, for every $i,j\in\{1,2,3\}$, we have $d(x_{i,i+1},x_{j,j+1})\le 2K_1+\frac{1}{3}$. This concludes the proof of our claim. 

For every $(\xi_1,\xi_2,\xi_3)\in(\partial_\infty X)^{(3)}$, we let $\alpha_\infty(\xi_1,\xi_2,\xi_3)$ be the closure of the increasing union of all $\alpha_n(\xi_1,\xi_2,\xi_3)$ as $n$ ranges over $\mathbb{N}$. This is a nonempty closed subset of $X$ of bounded diameter.

For every $n\in\mathbb{N}$, we then let $\theta_n:(\partial_\infty X)^{(3)}\to X$ be the map that sends $(\xi_1,\xi_2,\xi_3)$ to the circumcenter of the finite set $\alpha_n(\xi_1,\xi_2,\xi_3)$ if this set is nonempty, and to a fixed basepoint $x_0\in X$ otherwise. This map is Borel as it is a composition of a Borel map and a continuous map (the map $c:\calp_f(X_\mathbb{Q})\to X$ is continuous because $\calp_f(X_\mathbb{Q})$ is equipped with the discrete topology). Lemma~\ref{lemma:center-countable-bounded-set} ensures that for every triple $(\xi_1,\xi_2,\xi_3)\in (\partial_\infty X)^{(3)}$, the points $\theta_n(\xi_1,\xi_2,\xi_3)$ converge to a point $\theta(\xi_1,\xi_2,\xi_3)$ (equal to the circumcenter of $\alpha_\infty(\xi_1,\xi_2,\xi_3)$). The map $\theta:(\partial_\infty X)^{(3)}\to X$ is then Borel, being a pointwise limit of Borel maps. 

Finally, using the $G$-invariance of $X_\mathbb{Q}$ and the fact that the $G$-action by isometries on $X$ extends continuously to a $G$-action on $X\cup\partial_\infty X$, we see that all maps $\alpha_n$ are $G$-equivariant. Therefore all maps $\theta_n$ are $G$-equivariant, and in the limit $\theta$ is $G$-equivariant.
\end{proof}

\section{Boundary amenability}\label{sec:boundary-amenability}

\emph{In this section, we prove that every countable group $G$ acting discretely by isometries on a connected piecewise hyperbolic $\cat$ simplicial complex $X$ with countably many simplices that belong to finitely many isometry types, and boundary amenable vertex stabilizers, is boundary amenable. As a consequence, we derive that all $2$-dimensional Artin groups of hyperbolic type are boundary amenable, using their action on the modified Deligne complex. Also, we show that, under the assumption that all edge stabilizers in $X$ are amenable, the action of a group $G$ as above on the visual boundary $\partial_\infty X$ is Borel amenable: this will be useful in later sections to tackle the measure equivalence rigidity questions.}

\subsection{Piecewise hyperbolic simplicial complexes}

We recall the definition of a piecewise hyperbolic simplicial complex from \cite[Definition~I.7.2]{BH}. Let $\mathbb{H}^2$ be the hyperbolic plane. Let $\Lambda$ be a set; for every $\lambda\in\Lambda$, let $n_\lambda\in\mathbb{N}$, and let $S_\lambda$ be a closed geodesic simplex in $\mathbb{H}^{n_\lambda}$. Let $$X:=\bigcup_{\lambda\in\Lambda}(S_\lambda\times\{\lambda\}),$$ let $\sim$ be an equivalence relation on $X$, and let $K:=X/{\sim}$. Let $p:X\to K$ be the natural projection, and for every $\lambda\in\Lambda$, let $p_\lambda:S_\lambda\to K$ be defined by letting $p_\lambda(x):=p(x,\lambda)$ for every $x\in S_\lambda$. Following \cite[Definition~I.7.2]{BH}, we say that $K$ is a \emph{piecewise hyperbolic simplicial complex} if 
\begin{enumerate}
\item for every $\lambda\in\Lambda$, the map $p_\lambda$ is injective, and
\item for every $\lambda,\lambda'\in\Lambda$ such that $p_\lambda(S_\lambda)\cap p_{\lambda'}(S_{\lambda'})\neq\emptyset$, there exists an isometry $h_{\lambda,\lambda'}$ from a face $T_\lambda$ of $S_\lambda$ onto a face $T_{\lambda'}$ of $S_{\lambda'}$ such that for every $x\in S_\lambda$ and every $x'\in S_{\lambda'}$, one has $p(x,\lambda)=p(x',\lambda')$ if and only if $x\in T_\lambda$, $x'\in T_{\lambda'}$ and $x'=h_{\lambda,\lambda'}(x)$.
\end{enumerate}
Every piecewise hyperbolic simplicial complex naturally comes equipped with a pseudometric, where the distance between two points $x$ and $y$ is defined as the infimum of the total length of a concatenation of geodesic segments from $x$ to $y$, each of these segments being contained in a simplex, see \cite[Definition~I.7.4]{BH}. When $K$ is connected and has finitely many isometry types of simplices, this pseudometric is in fact a metric that turns $K$ into a complete geodesic metric space \cite[Theorem~I.7.19]{BH}. 

\subsection{Statements} 

Let $G$ be a countable group acting on a topological space $X$ by Borel automorphisms. The $G$-action on $X$ is \emph{Borel amenable} if there exists a sequence of Borel maps $\nu_n:X\to\Prob(G)$ which is \emph{asymptotically equivariant}, i.e.\ such that for all $x\in X$ and all $g\in G$, one has $||\nu_n(gx)-g\cdot\nu_n(x)||_1\to 0$ as $n$ goes to $+\infty$. Here the topology on $\Prob(G)$ is the subspace topology from $\ell^1(G)$, i.e.\ the topology of pointwise convergence. A countable group $G$ is \emph{boundary amenable} (or \emph{exact}) if it admits a Borel amenable action on some nonempty compact Hausdorff topological space (equipped with its Borel $\sigma$-algebra). The goal of the present section is to prove Theorems~\ref{theo:exactness} and~\ref{theo:amenability-boundary} below.

\begin{theo}\label{theo:exactness}
Let $X$ be a connected piecewise hyperbolic $\cat$ simplicial complex with countably many simplices that belong to finitely many isometry types. Let $G$ be a countable group acting discretely on $X$ by simplicial isometries. Assume that all vertex stabilizers are boundary amenable.

Then $G$ is boundary amenable.
\end{theo}

This implies in particular that $G$ satisfies the Novikov conjecture \cite{Yu,HR,Hig}. 

We believe that Theorem~\ref{theo:exactness} and its application to Artin groups (Theorem~\ref{theo:exactness-Artin} below) are of independent interest, but we will not use them in the remainder of the present paper. However, we will need the following theorem, which will follow from the same argument -- but none of the two theorems seems to follow from the other, and in particular we point out that the space $\partial_\infty X$ need not be compact in Theorem~\ref{theo:amenability-boundary} below.

\begin{theo}\label{theo:amenability-boundary}
Let $X$ be a connected piecewise hyperbolic $\cat$ simplicial complex with countably many simplices that belong to finitely many isometry types. Let $G$ be a countable group acting discretely on $X$ by simplicial isometries. Assume that all edge stabilizers are amenable.

Then the $G$-action on $\partial_\infty X$ is Borel amenable.
\end{theo}

Given a simplicial complex $X$, we denote by $\cals^{\ge 1}(X)$ the collection of all simplices of $X$ of dimension at least $1$ (i.e.\ excluding vertices). Using an argument due to Ozawa \cite[Proposition~11]{Oza}, we will see that Theorems~\ref{theo:exactness} and~\ref{theo:amenability-boundary} are consequences of the following proposition.

\begin{prop}\label{prop:proba-on-edges}
Let $X$ be a connected piecewise hyperbolic $\cat$ simplicial complex with countably many simplices that belong to finitely many isometry types. Let $G$ be a countable group acting discretely on $X$ by simplicial isometries. 

Then there exists an asymptotically equivariant sequence of Borel maps $\partial_\infty X\to\Prob(\cals^{\ge 1}(X))$.
\end{prop}

Before moving on to the proof of Proposition~\ref{prop:proba-on-edges}, we first explain how to derive Theorems~\ref{theo:exactness} and~\ref{theo:amenability-boundary}.

\begin{proof}[Proof of Theorem~\ref{theo:exactness}]
Recall that $\horox$ denotes the horofunction compactification of $X$. As $X$ is a complete metric space \cite[Theorem~I.7.19]{BH}, Proposition~\ref{prop:horo-partition} yields a Borel $G$-equivariant map $\horox\to X\cup\partial_\infty X$. Denote by $\cals(X)$ the countable collection of all simplices of $X$. Then there is a Borel $G$-equivariant map $X\to\Prob(\cals(X))$, sending a point $x$ to the Dirac measure at the unique simplex of minimal dimension that contains $x$. Also, Proposition~\ref{prop:proba-on-edges} yields an asymptotically equivariant sequence of Borel maps $\partial_\infty X\to\Prob(\cals(X))$. Combining the above maps, we get an asymptotically equivariant sequence of Borel maps $\nu_n:\horox\to\Prob(\cals(X))$. 

By assumption, the $G$-stabilizer of every vertex of $X$ is boundary amenable. As boundary amenability is stable under subgroups (as follows from \cite[Proposition~11]{Oza}) and finite-index extensions (see e.g.\ \cite{KW}), we deduce that the $G$-stabilizer of every element of $\cals(X)$ is boundary amenable. As $\horox$ is compact Hausdorff and $\cals(X)$ is countable, it thus follows from \cite[Proposition~C.1]{Kid-memoir} (attributed to Ozawa) that $G$ is boundary amenable.  
\end{proof} 

\begin{proof}[Proof of Theorem~\ref{theo:amenability-boundary}]
Our assumption that edge stabilizers are amenable implies that the stabilizer of every element of $\cals^{\ge 1}(X)$ is amenable. As $\cals^{\ge 1}(X)$ is countable, Theorem~\ref{theo:amenability-boundary} thus follows from Proposition~\ref{prop:proba-on-edges} together with  \cite[Proposition~2.12]{BGH} (which is a consequence of \cite[Proposition~11]{Oza}).
\end{proof}

\subsection{Construction of the measures} 

The rest of the section is devoted to the proof of Proposition~\ref{prop:proba-on-edges}. From now on, we let $X$ be a connected piecewise hyperbolic $\cat$ simplicial complex with countably many simplices that belong to finitely many isometry types, and $G$ be a countable group acting discretely by simplicial isometries on $X$.

\paragraph*{Cornered simplices.} By \cite[Corollary~I.7.29]{BH}, for every $x_0\in X$, every geodesic ray $r$ from $x_0$ to a point $\xi\in\partial_\infty X$ and every $t\in\mathbb{R}_+$, the image of $r_{|[0,t]}$ only meets finitely many simplices of $X$ in their interior. Denote by $\calp_f(\cals^{\ge 1}(X))$ the set of all finite subsets of $\cals^{\ge 1}(X)$. For every $x_0\in X$ and every $n\in\mathbb{N}$, we let $\Sigma'_{n,x_0}:\partial_\infty X\to\calp_f(\cals^{\ge 1}(X))$ be the map that sends $\xi\in\partial_\infty X$ to the set of all simplices of dimension at least $1$ that meet the image of $r_{|[0,n]}$ in their interior, where $r$ is the geodesic ray from $x_0$ to $\xi$. We then let $\Sigma_{n,x_0}:\partial_\infty X\to\calp_f(\cals^{\ge 1}(X))$ be the map defined by letting $\Sigma_{n,x_0}(\xi)$ be the set of all simplices in $\Sigma'_{n,x_0}(\xi)$ together with all their faces in $\cals^{\ge 1}(X)$.

Let $r$ be a geodesic ray in $X$, let $\sigma$ be a simplex of $X$, and let $\epsilon>0$. We say that $\sigma$ is \emph{$\epsilon$-cornered} by $r$ if the image of $r$ intersects $\sigma$ in a nonempty set contained in the $\epsilon$-ball around some vertex of $\sigma$.

\begin{lemma}\label{lemma:cornered}
There exist $\epsilon_0>0$ and $p_0\in [0,1)$ such that for every $x_0\in X$, every $\xi\in\partial_\infty X$ and every $n\in\mathbb{N}$, the proportion of simplices in $\Sigma'_{n,x_0}(\xi)$ which are $\epsilon_0$-cornered by the geodesic ray from $x_0$ to $\xi$ is at most $p_0$. 
\end{lemma}

\begin{proof} 
Let $x_0\in X$ and $\xi\in\partial_\infty X$, and let $r$ be the geodesic ray from $x_0$ to $\xi$. For every $n\in\mathbb{N}$, we denote by $\alpha_n$ the cardinality of $\Sigma'_{n,x_0}(\xi)$. By \cite[Corollary~I.7.30]{BH}, there exists $C>0$ -- which only depends on $X$ and not on $x_0$ and $\xi$ -- such that for every $n\in\mathbb{N}$, one has $\alpha_n\le Cn$. Let $\epsilon_0>0$ be chosen so that $2\epsilon_0 C<1$. For every $n\in\mathbb{N}$, let $\kappa_n$ be the proportion of simplices in $\Sigma'_{n,x_0}(\xi)$ that are $\epsilon_0$-cornered by $r$ -- so that there are $\kappa_n\alpha_n$ such simplices. We aim to show that $\kappa_n$ is uniformly bounded away from $1$ -- with a bound independent on $n$ and on $x_0$ and $\xi$.

Let $M>0$ be large enough so that $CM>1$ and every simplex in $X$ has diameter smaller than $M$. By decomposing the length $n$ of $r_{|[0,n]}$ into the total length spent in $\epsilon_0$-cornered simplices and the total length spent in other simplices, we get $$n\le 2\epsilon_0\kappa_n\alpha_n+ M(1-\kappa_n)\alpha_n,$$ i.e.\ $n\le \alpha_n[M-\kappa_n(M-2\epsilon_0)]$. On the other hand, as $\alpha_n\le Cn$, we get $$1\le CM-C\kappa_n(M-2\epsilon_0).$$ This leads to $$\kappa_n\le\frac{CM-1}{CM-2C\epsilon_0},$$ which is strictly smaller than $1$ by our choice of $\epsilon_0$ -- and the bound only depends on $X$, as desired. 
\end{proof}

We will now fix $\epsilon_0>0$ as in Lemma~\ref{lemma:cornered} once and for all, which we choose sufficiently small so that for every simplex $\sigma$ in $\cals(X)$, the balls of radius $\epsilon_0$ around the vertices of $\sigma$ are pairwise disjoint in $\sigma$.

\begin{figure}
\centering
\includegraphics[scale=0.6]{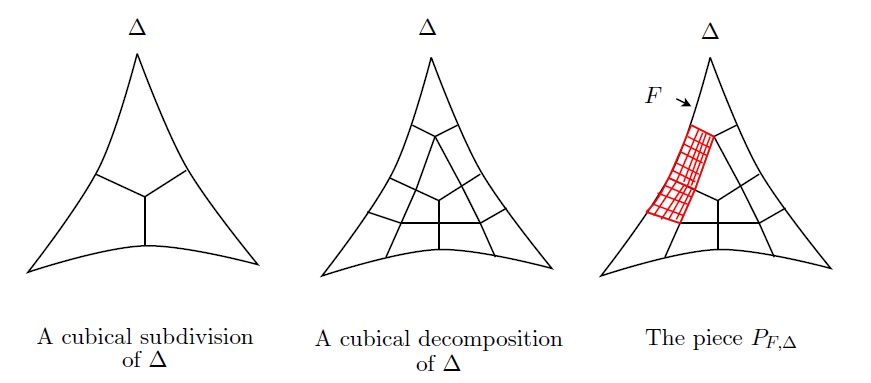}
\caption{Cubical subdivisions, cubical decompositions and pieces.}
\label{fig:cubical}
\end{figure}

\paragraph*{Cubical decompositions and pieces.}
Let $\sigma$ be a simplex of $X$. The following notions are illustrated in Figure~\ref{fig:cubical} in the case where $\sigma$ has dimension $2$. 

A \emph{cubical subdivision} of $\sigma$ is a decomposition of $\sigma$ as a cube complex whose vertices are the barycenters of $\sigma$ and of its faces, for which there is an edge $e$ between vertices corresponding to faces $\si_1$ and $\si_2$ precisely when $\si_1$ is codimension $1$ in $\si_2$; this edge $e$ is then homeomorphic to a segment contained in $\si_2$ which meets $\si_1$ only at its extremity. Notice that a cubical subdivision of $\sigma$ naturally induces a cubical subdivision of each of its faces.

In an analogous way, one can define a cubical subdivision of a cube instead of a simplex. A \emph{cubical decomposition} of $\sigma$ is then a decomposition of $\sigma$ as a cube complex obtained from a cubical subdivision of $\sigma$ by cubically subdividing each $n$-cube into $2^n$ subcubes. As above, every cubical decomposition of $\sigma$ naturally induces a cubical decomposition of each of its faces.

Given a simplex $\sigma'$ which is either equal to $\sigma$ or to a face of $\sigma$, and a cubical decomposition of $\sigma$, the \emph{piece} $P_{\sigma',\sigma}$ is the star of the barycenter of $\sigma'$ in the cubical decomposition.

We now equip each simplex of $X$ with a cubical decomposition, in such a way that
\begin{enumerate}
	\item if $\si\subseteq\si'$, then the decomposition of $\sigma$ is induced from the decomposition of $\sigma'$;
	\item if $\si_1$ and $\si_2$ are isometric simplices in $X$, then any isometry between them respects their cubical decompositions;
	\item the piece in $\si$ associated with a vertex $v$ of $\si$ is contained in the $\eps_0$-neighborhood of $v$.
\end{enumerate}

Given a simplex $\sigma\in\cals^{\ge 1}(X)$, we denote by $\st(\sigma)$ the \emph{star} of $\sigma$, defined as the union of $\sigma$ and of all simplices in $X$ that contain $\sigma$ as a face. We define the \emph{interior} of $\st(\si)$, denoted by $\mathring{\st}(\si)$, to be the collection of points $x$ in $\st(\si)$ such that the smallest simplex $\si_x$ containing $x$ satisfies $\si\subseteq\si_x$. Equivalently, a point $x\in X$ belongs to $\mathring{\st}(\si)$ if and only if $x$ belongs to the interior of a simplex which is either equal to $\si$ or contains $\si$ as a face.

\paragraph*{Weight functions.}
Let $w_0>0$ be chosen once and for all. For every simplex $\sigma\in\cals^{\ge 1}(X)$, we choose a function $f_{\sigma}:X\to \mathbb R_{\ge 0}$, in such a way that
\begin{enumerate}
\item $f_{\sigma}(x)>0$ for every $x$ lying in the interior of $\st(\sigma)$, and $f_{\sigma}(x)=0$ for every $x$ lying outside of $\st(\si)$;
	\item for every simplex $\si'$ that contains $\si$ as a face, and every $x\in P_{\si,\si'}$, we have $f_{\sigma}(x)\ge w_0$;
	\item there exists $M>0$ such that every function $f_{\sigma}$ is bounded from above by $M$, and there exists $L_0>0$ such that every function $f_{\sigma}$ is $L_0$-Lipschitz;
	\item for every $g\in G$, every $\sigma\in\cals^{\ge 1}(X)$ and every $x\in X$, one has $f_{g\sigma}(gx)=f_{\sigma}(x)$.
\end{enumerate}

Notice that the reason why we can choose a uniform $L_0$ in (3) is because $X$ has finitely many isometry types of simplices.

\begin{de}[Weight function]
Let $x_0\in X$, let $\sigma\in\cals^{\ge 1}(X)$. The \emph{weight function} of $\sigma$ based at $x_0$ is the function $w_{\sigma,x_0}:\partial_\infty X\to \mathbb{R}_{\ge 0}$ defined by letting $w_{\sigma,x_0}(\xi)$ be the maximal value of $f_{\sigma}(x)$ when $x$ varies over all points lying on the image of the geodesic ray from $x_0$ to $\xi$. 
\end{de}

Notice that the functions $w_{\sigma,x_0}$ all take finite (in fact bounded) values because the functions $f_\sigma$ are bounded by assumption. In addition, for every $g\in G$, every $x_0\in X$, every $\sigma\in\cals^{\ge 1}(X)$ and every $\xi\in\partial_\infty X$, one has $w_{g\sigma,gx_0}(g\xi)=w_{\sigma,x_0}(\xi)$.

\begin{lemma}\label{lemma:weight-significant-contribution}
There exists $p>0$ such that for every $x_0\in X$, every $\xi\in\partial_\infty X$ and every $n\in\mathbb{N}$, the proportion of simplices $\sigma$ in $\Sigma_{n,x_0}(\xi)$ with $w_{\sigma,x_0}(\xi)\ge w_0$ is at least $p$.
\end{lemma}

\begin{proof}
Recall that $\epsilon_0>0$ has been chosen as in Lemma~\ref{lemma:cornered}, and such that for every simplex $\sigma$ in $\cals(X)$, the balls of radius $\epsilon_0$ around the vertices of $\sigma$ are pairwise disjoint. In view of Lemma~\ref{lemma:cornered}, it is enough to show that for every $x_0\in X$ and every $\xi\in\partial_\infty X$, denoting by $r$ the geodesic ray from $x_0$ to $\xi$, for every simplex $\sigma'$ which meets $r$ in its interior and is not $\epsilon_0$-cornered by $r$, there exists a simplex $\sigma$ which is either equal to $\sigma'$ or to a face of $\sigma'$ such that $w_{\sigma,x_0}(\xi)\ge w_0$. Let $x$ be a point on the image of $r$ which is outside the $\eps_0$-neighborhood of every vertex of $\si'$. Then there exists a simplex $\si$, either equal to $\si'$ or to a face of $\si'$, such that the piece $P_{\si,\si'}$ contains $x$. Thus $f_{\si}(x)\ge w_0$ by the definition of $f_\sigma$. Hence $w_{\sigma,x_0}(\xi)\ge w_0$.
\end{proof}

\paragraph*{Definition of the measures.}

For every $x_0\in X$, every $n\in\mathbb{N}$, and every $\xi\in\partial_\infty X$, we let $$\Upsilon_{n,x_0}(\xi):=\sum_{\si\in\cals^{\ge 1}(X)}\mathbf{1}_{\Sigma_{n,x_0}(\xi)}(\sigma)w_{\si,x_0}(\xi).$$ We then define maps $\nu_{n,x_0}:\partial_\infty X\to\Prob(\cals^{\ge 1}(X))$ by letting $$\nu_{n,x_0}(\xi):=\frac{1}{\Upsilon_{n,x_0}(\xi)}\left(\sum_{\sigma\in\cals^{\ge 1}(X)}\mathbf{1}_{\Sigma_{n,x_0}(\xi)}(\sigma) w_{\sigma,x_0}(\xi)\delta_\sigma\right),$$ where $\delta_\sigma$ denotes the Dirac mass at $\sigma$. Notice that for every $g\in G$, every $x_0\in X$, every $n\in\mathbb{N}$ and every $\xi\in\partial_\infty X$, one has $\Upsilon_{n,gx_0}(g\xi)=\Upsilon_{n,x_0}(\xi)$ and $\nu_{n,gx_0}(g\xi)=g\nu_{n,x_0}(\xi)$. 

\subsection{The key estimate}

In this section, we establish the key estimate regarding the measures $\nu_{n,x}$ defined above for our proof of Proposition~\ref{prop:proba-on-edges}: this is Lemma~\ref{lemma:asymptotically-equivariant} below. We start with a general lemma for comparing probability measures. 

\begin{lemma}\label{lemma:calculus}
Let $\cald$ be a countable set. Let $\epsilon\in (0,1)$. For every $n\in\mathbb{N}$, let $A_n\subseteq\cald$ be a finite subset, let $\alpha_n,\alpha'_n:A_n\to\mathbb{R}_+$ be two functions, let $$\Upsilon_n:=\sum_{\sigma\in A_n}\alpha_n(\sigma),~~~~\Upsilon'_n:=\sum_{\sigma\in A_n}\alpha'_n(\sigma),$$ and let $\nu_n$ and $\nu'_n$ be the probability measures on $\cald$ defined by $$\nu_n:=\frac{1}{\Upsilon_n}\sum_{\sigma\in A_n}\alpha_{n}(\sigma)\delta_{\sigma},~~~~\nu'_n:=\frac{1}{\Upsilon'_n}\sum_{\sigma\in A_n}\alpha'_n(\sigma)\delta_\sigma.$$ Assume that there exist $K>0$ and $M>0$ such that 
\begin{enumerate}
\item for every $n\in\mathbb{N}$ and every $\sigma\in A_n$, one has $\alpha_{n}(\sigma)\le M$ and $\alpha'_n(\sigma)\le M$;
\item for every $n\in\mathbb{N}$, one has $\Upsilon_n\ge K |A_n|$ and $\Upsilon'_n\ge K |A_n|$, and $|A_n|\to +\infty$; 
\item the exists $m\in\mathbb{N}$ such that for every $n\in\mathbb{N}$ and all but at most $m$ elements $\sigma\in A_n$, one has $|\alpha_n(\sigma)-\alpha'_n(\sigma)|\le \min\{\frac{K}{3},\frac{K^2}{6M}\}\epsilon.$
\end{enumerate}
\noindent Then for all sufficiently large $n\in\mathbb{N}$, one has $||\nu_n-\nu'_n||_1\le\epsilon$.
\end{lemma}

\begin{proof}
Let $C:=\min\{\frac{K}{3},\frac{K^2}{6M}\}$. We first claim that for all sufficiently large $n\in\mathbb{N}$, one has $$\left|\frac{\Upsilon'_n}{\Upsilon_n}-1\right|\le \frac{2C}{K}\epsilon.$$ Indeed, for every $n\in\mathbb{N}$, let $A_n^*\subseteq A_n$ be the subset made of all $\sigma\in A_n$ such that $|\alpha_{n}(\sigma)-\alpha'_n(\sigma)|\le C\epsilon$. By assumption, for every $n\in\mathbb{N}$, the set $A_n\setminus A_n^*$ has cardinality at most $m$. For every $n\in\mathbb{N}$, we have $$|\Upsilon'_n-\Upsilon_n|\le\left(\sum_{\sigma\in A_n^*}|\alpha_n(\sigma)-\alpha'_n(\sigma)|\right)+2mM\le C\epsilon |A_n^*|+2mM,$$ so $$\left|\frac{\Upsilon'_n}{\Upsilon_n}-1\right|\le\frac{C}{K}\epsilon+\frac{2mM}{\Upsilon_n}.$$ Since by assumption $\Upsilon_n$ tends to $+\infty$ as $n$ goes to $+\infty$, the claim follows. 

For every $n\in\mathbb{N}$, we then have $$||\nu_n-\nu'_n||_1\le\frac{1}{\Upsilon_n}\left(\sum_{\sigma\in A_n}\left|\alpha_n(\sigma)-\frac{\Upsilon_n}{\Upsilon'_n}\alpha'_n(\sigma)\right|\right).$$ Since for all sufficiently large $n$ the ratio $\frac{\Upsilon_n}{\Upsilon_n'}$ is at most $3$, we deduce that for all sufficiently large $n\in\mathbb{N}$, one has
$$||\nu'_n-\nu_n||_1\le \frac{1}{\Upsilon_n}\left(\left(\sum_{\sigma\in A_n^*}\left|\alpha_n(\sigma)-\frac{\Upsilon_n}{\Upsilon'_n}\alpha'_n(\sigma)\right|\right)+4mM\right).$$ For every $n\in\mathbb{N}$ and every $\sigma\in A_n^*$, we have $$\left|\alpha_n(\sigma)-\frac{\Upsilon_n}{\Upsilon'_n}\alpha'_{n}(\sigma)\right|\le |\alpha_n(\sigma)-\alpha'_n(\sigma)|+\left|\frac{\Upsilon_n}{\Upsilon'_n}-1\right|\alpha'_n(\sigma)\le \left(C+\frac{2CM}{K}\right)\epsilon.$$ Therefore $$||\nu_n-\nu'_n||_1\le\frac{1}{\Upsilon_n}\left(C+\frac{2CM}{K}\right)\epsilon |A_n^*|+\frac{4mM}{\Upsilon_n}\le \left(\frac{C}{K}+\frac{2CM}{K^2}\right)\epsilon+\frac{4mM}{\Upsilon_n}\le \frac{2}{3}\epsilon+\frac{4mM}{\Upsilon_n}.$$ The lemma follows. 
\end{proof}

The following lemma is the key estimate in our proof of Proposition~\ref{prop:proba-on-edges}.

\begin{lemma}\label{lemma:asymptotically-equivariant}
For every $x_1,x_2\in X$ and every $\xi\in\partial_\infty X$, one has $$||\nu_{n,x_1}(\xi)-\nu_{n,x_2}(\xi)||_1\to 0$$ as $n$ goes to $+\infty$. 
\end{lemma}

\begin{proof}
We fix $x_1,x_2$ and $\xi$ once and for all. Let $\eps>0$. Using the criterion from Lemma~\ref{lemma:calculus}, we aim to prove that for every sufficiently large $n\in\mathbb{N}$, we have $$||\nu_{n,x_1}(\xi)-\nu_{n,x_2}(\xi)||_1\le \eps.$$
Using \cite[Corollary~I.7.30]{BH}, we can (and shall) choose $C\ge 1$ such that for every $n\in\mathbb{N}$, one has $$\frac{1}{C}|\Sigma_{n,x_1}(\xi)|\le |\Sigma_{n,x_2}(\xi)|\le C |\Sigma_{n,x_1}(\xi)|,$$ and these quantities all tend to $+\infty$ as $n$ goes to $+\infty$. Let $L_0>0$ be such that for every simplex $\sigma$, the function $f_{\sigma}$ is $L_0$-Lipschitz. Let $M$ be a global upper bound to all functions $f_{\sigma}$. Let $p>0$ be the number provided by Lemma~\ref{lemma:weight-significant-contribution}, and let $K:=\frac{pw_0}{C+1}$. Fix $\epsilon'\in (0,1)$ with $\epsilon'<\frac{1}{L_0}\min\{\frac{K}{3},\frac{K^2}{6M}\}\epsilon$. 

For every $i\in\{1,2\}$, let $r_i$ be the geodesic ray from $x_i$ to $\xi$. Since $X$ is $\cat$, there exist $t_0\in\mathbb{R}_+$ and $\tau\in\mathbb{R}$ such that for every $t\ge t_0$, one has $d(r_1(t),r_2(t+\tau))\le\epsilon'$. Up to reversing the roles of $x_1$ and $x_2$, we can (and shall) assume that $\tau\ge 0$. Let $D$ be the maximal diameter of a simplex in $X$ (this exists and it is finite as $X$ has only finitely many isometry types of simplices). For every $n\ge t_0+\tau+3D+1$, we let $$\Sigma_{n,x_1}^*(\xi):=\{\si\in\Sigma_{n,x_1}(\xi)|\mathring{\st}(\si)\cap r_1(\mathbb{R}_+)\subseteq r_1([t_0,n-\tau-(3D+1)])\}$$ and $$\Sigma_{n,x_2}^*(\xi):=\{\si\in\Sigma_{n,x_2}(\xi)|\mathring{\st}(\si)\cap r_2(\mathbb{R}_+)\subseteq r_2([t_0+\tau,n-(3D+1)])\}.$$ 

\begin{claim}\label{claim:bh}
There exists $m\ge 0$ such that for every $i\in\{1,2\}$ and every $n\ge t_0+\tau+3D+1$, one has $|\Sigma_{n,x_i}(\xi)\setminus\Sigma_{n,x_i}^*(\xi)|\le m.$
\end{claim}

\begin{proof}[Proof of Claim~\ref{claim:bh}]
We will prove the claim for $i=1$, the proof for $i=2$ being similar. As every simplex of $X$ has diameter at most $D$, the geodesic ray $r_1$ cannot meet the star of a simplex $\sigma$ both before time $n-3D$ and after time $n$. Consequently, if  $\sigma\in\Sigma_{n,x_1}(\xi)\setminus\Sigma_{n,x_1}^*(\xi)$, then $\mathring{\st}(\sigma)$ meets either $r_1([0,t_0])$ or $r_1([n-\tau-(3D+1),n])$. In addition, using \cite[Corollary~I.7.30]{BH}, the number of simplices $\sigma$ such that $\mathring{\st}(\sigma)$ meets $r_1([0,t_0])$ is finite, and the number of simplices $\sigma$ such that $\mathring{\st}(\sigma)$ meets $r_1([n-\tau-(3D+1),n])$ is bounded, with a bound independent of $n$. The claim follows.
\end{proof}

\begin{claim}\label{claim:upsilon}
For every $n\in\mathbb{N}$ and every $i\in\{1,2\}$, we have $\Upsilon_{n,x_i}(\xi)\ge K |\Sigma_{n,x_1}(\xi)\cup\Sigma_{n,x_2}(\xi)|$.
\end{claim}

\begin{proof}[Proof of Claim~\ref{claim:upsilon}]
By Lemma~\ref{lemma:weight-significant-contribution}, for every $n\in\mathbb{N}$ and every $i\in\{1,2\}$, we have $\Upsilon_{n,x_i}(\xi)\ge pw_0|\Sigma_{n,x_i}(\xi)|$. The conclusion follows because for every $i\in\{1,2\}$, we have $|\Sigma_{n,x_1}(\xi)\cup\Sigma_{n,x_2}(\xi)|\le (C+1)|\Sigma_{n,x_i}(\xi)|.$
\end{proof}

\begin{claim}\label{claim:lipschitz}
For every $n\ge t_0+\tau+3D+1$ and every $\sigma\in\Sigma_{n,x_1}^*(\xi)\cap\Sigma_{n,x_2}^*(\xi)$, one has $|w_{\sigma,x_1}(\xi)-w_{\sigma,x_2}(\xi)|\le \min\{\frac{K}{3},\frac{K^2}{6M}\}\epsilon$.
\end{claim}

\begin{proof}[Proof of Claim~\ref{claim:lipschitz}]
Let $x_\sigma$ be a point on the image of $r_1$ at which the function $f_\sigma$ attains its maximum. Then $x_{\sigma}\in\mathring{\st}(\si)$. As $\sigma$ belongs to $\Sigma_{n,x_1}^*(\xi)$, we have $x_\sigma\in r_1([t_0,+\infty))$. Then there exists a point $y_\sigma$ on the image of $r_2$ at distance at most $\epsilon'$ from $x_\sigma$. Thus $$w_{\sigma,x_1}(\xi)=f_{\sigma}(x_{\sigma})\le f_{\sigma}(y_{\sigma})+L_0\epsilon'\le w_{\sigma,x_2}(\xi)+L_0\epsilon'.$$ Reversing the roles of $x_1$ and $x_2$ shows that $|w_{\sigma,x_1}(\xi)-w_{\sigma,x_2}(\xi)|\le L_0\epsilon'$, which concludes the proof of the claim.
\end{proof}

\begin{claim}\label{claim:zero}
For every $n\ge t_0+\tau+3D+1$ and every $\sigma\in\Sigma_{n,x_1}^*(\xi)\setminus\Sigma_{n,x_2}(\xi)$, one has $w_{\sigma,x_2}(\xi)=0$ and $w_{\sigma,x_1}(\xi)\le \min\{\frac{K}{3},\frac{K^2}{6M}\}\epsilon$. Likewise, for every $n\ge t_0+\tau+3D+1$ and every $\sigma\in\Sigma_{n,x_2}^*(\xi)\setminus\Sigma_{n,x_1}(\xi)$, one has $w_{\sigma,x_1}(\xi)=0$ and $w_{\sigma,x_2}(\xi)\le \min\{\frac{K}{3},\frac{K^2}{6M}\}\epsilon$.  
\end{claim}

\begin{proof}[Proof of Claim~\ref{claim:zero}]
We will only prove the first assertion, the proof of the second being similar. Let $\sigma\in\Sigma_{n,x_1}^*(\xi)\setminus\Sigma_{n,x_2}(\xi)$. Let $\Sigma_{\infty,x_2}(\xi)$ denote the set of all $\sigma\in\cals^{\ge 1}(X)$ that appear as faces of simplices traversed in their interior by the ray $r_2$ (equivalently, this is the union over all $n\in\mathbb{N}$ of the sets $\Sigma_{n,x_2}(\xi)$). 

We first claim that $\sigma\notin\Sigma_{\infty,x_2}(\xi)$. Indeed, assume towards a contradiction that $\sigma\in\Sigma_{\infty,x_2}(\xi)$. This implies that there exists $t\in\mathbb{R}_+$ such that $\sigma$ is a face of a simplex containing $r_2(t)$ in its interior. In particular $d(r_2(t),\sigma)\le D$, where we recall that $D$ is the maximal diameter of a simplex in $X$. On the other hand, as $\sigma\in \Sigma_{n,x_1}^*(\xi)$, there exists $t'\in [t_0,n-\tau-(3D+1)]$ such that $d(r_1(t'),\sigma)\le D$, which implies that $d(r_2(t'+\tau),\sigma)\le D+\epsilon'\le D+1$. Since $\sigma$ has diameter at most $D$, we deduce that $d(r_2(t),r_2(t'+\tau))\le 3D+1$. As $r_2$ is a geodesic ray, we thus get that $t\le t'+\tau+3D+1\le n$. This shows that $\sigma\in\Sigma_{n,x_2}(\xi)$, a contradiction. 

We can now conclude the proof of Claim~\ref{claim:zero}. As $\sigma\notin\Sigma_{\infty,x_2}(\xi)$, we have $f_{\sigma}(r_2(t))=0$ for every $t\in\mathbb{R}_+$, which implies that $w_{\sigma,x_2}(\xi)=0$. In addition, let $x$ be a point lying on the image of $r_1$ at which the function $f_{\sigma}$ achieves its maximum. As $\sigma\in\Sigma_{n,x_1}^*(\xi)$, we have $x\in r_1([t_0,+\infty))$. Therefore, there exists a point $y$ on the image of $r_2$ at distance at most $\epsilon'$ from $x$. As $f_{\sigma}(y)=0$ and $f_{\sigma}$ is $L_0$-Lipschitz, we deduce that $f_{\sigma}(x)\le L_0\epsilon'$. Thus $w_{\sigma,x_1}(\xi)\le L_0\epsilon'\le\min\{\frac{K}{3},\frac{K^2}{6M}\}\epsilon$.
\end{proof}

We are now in position to apply the criterion provided by Lemma~\ref{lemma:calculus} to the measures $\nu_{n,x_1}(\xi)$ and $\nu_{n,x_2}(\xi)$. For every $n\in\mathbb{N}$, let $A_n:=\Sigma_{n,x_1}(\xi)\cup\Sigma_{n,x_2}(\xi)$. For every $\sigma\in A_n$, let 
$$\alpha_n(\sigma):=\left\{
\begin{array}{cl}
w_{\sigma,x_1}(\xi) & \text{if~} \si\in\Sigma_{n,x_1}(\xi)\\
0 & \text{otherwise}
\end{array}
\right.
~~~~~\alpha'_n(\sigma):=\left\{
\begin{array}{cl}
w_{\sigma,x_2}(\xi) & \text{if~} \si\in\Sigma_{n,x_2}(\xi)\\
0 & \text{otherwise}
\end{array}
\right.$$
 These are all bounded by $M$, showing that the first assumption from Lemma~\ref{lemma:calculus} holds. The second assumption has been checked in Claim~\ref{claim:upsilon}. To check the third assumption from Lemma~\ref{lemma:calculus}, notice that for every $n\in\mathbb{N}$, one has $$A_n\subseteq A_n^1\cup A_n^2\cup A_n^3\cup A_n^4\cup A_n^5,$$ with $$A_n^1:=\Sigma_{n,x_1}(\xi)\setminus\Sigma_{n,x_1}^*(\xi),~~~~ A_n^2:=\Sigma_{n,x_2}(\xi)\setminus\Sigma_{n,x_2}^*(\xi),$$ $$A_n^3:=\Sigma_{n,x_1}^*(\xi)\setminus\Sigma_{n,x_2}(\xi),~~~~ A_n^4:=\Sigma_{n,x_2}^*(\xi)\setminus\Sigma_{n,x_1}(\xi),$$ $$A_n^5:=\Sigma_{n,x_1}^*(\xi)\cap\Sigma_{n,x_2}^*(\xi).$$ The first two sets have uniformly bounded cardinality by Claim~\ref{claim:bh}, and for every $\sigma\in A_n^3\cup A_n^4\cup A_n^5$, one has $|w_{\sigma,x_1}(\xi)-w_{\sigma,x_2}(\xi)|\le\min\{\frac{K}{3},\frac{K^2}{6M}\}\epsilon$ (as follows from Claim~\ref{claim:zero} for $A_n^3$ and $A_n^4$, and from Claim~\ref{claim:lipschitz} for $A_n^5$). We can therefore apply Lemma~\ref{lemma:calculus}, which concludes our proof.
\end{proof}

\subsection{Measurability considerations and end of the proof} 

Before completing the proof of Proposition~\ref{prop:proba-on-edges}, we need to check the  measurability of all our constructions. We start with a general lemma in $\cat$ geometry.

\begin{lemma}\label{lemma:cat-rays}
Let $X$ be a $\cat$ space. Let $x_0\in X$, let $K\subseteq X$ be a compact set, and let $B_{x_0,K}\subseteq \partial_\infty X$ be the subset made of all boundary points $\xi$ such that the geodesic ray from $x_0$ to $\xi$ intersects $K$. 

Then $B_{x_0,K}$ is closed in $\partial_\infty X$.
\end{lemma}

\begin{proof}
As the cone topology on $\partial_\infty X$ is metrizable \cite[Section~5]{Vai}, it suffices to prove that for every sequence $(\xi_n)_{n\in\mathbb{N}}$ of boundary points in $B_{x_0,K}$ converging to a point $\xi\in\partial_\infty X$, we have $\xi\in B_{x_0,K}$. For every $n\in\mathbb{N}$, let $x_n\in K$ be a point lying on the geodesic ray from $x_0$ to $\xi_n$. Let $y_n$ be the point on the geodesic ray $\overline{x_0\xi}$ from $x_0$ to $\xi$ which is closest to $x_n$ (this is unique by $\cat$ geometry). Then the distance $d(x_n,y_n)$ converges to $0$ as $n$ goes to $+\infty$. As $K$ is compact, up to passing to a subsequence we can assume that $(x_n)_{n\in\mathbb{N}}$ converges to a point $x_\infty\in K$. Then $(y_n)_{n\in\mathbb{N}}$ also converges to $x_\infty$. Therefore $x_\infty\in K\cap \overline{x_0\xi}$, which finishes the proof.
\end{proof}

\begin{lemma}\label{lemma:list-of-simplices-borel}
For every $x_0\in X$ and every $n\in\mathbb{N}$, the map $\Sigma_{n,x_0}:\partial_\infty X\to\calp_f(\cals^{\ge 1}(X))$ is Borel.
\end{lemma}

\begin{proof}
It is enough to prove that the map $\Sigma'_{n,x_0}$ is Borel; for that, it is enough to prove that for every $\sigma\in\cals^{\ge 1}(X)$, the subset $A_{\sigma}\subseteq\partial_\infty X$ made of all boundary points $\xi$ such that the geodesic ray from $x_0$ to $\xi$ passes through the interior of $\sigma$ is Borel. The interior of $\sigma$ is an increasing union of countably many compact sets $K_{n,\sigma}$ with $n\in\mathbb{N}$. Therefore, with the notation from Lemma~\ref{lemma:cat-rays}, the set $A_\sigma$ is equal to the countable union over $n\in\mathbb{N}$ of all sets $B_{x_0,K_{n,\sigma}}$. As each of these sets is closed (Lemma~\ref{lemma:cat-rays}), the lemma follows.  
\end{proof}

\begin{lemma}\label{lemma:weight-continuous}
For every $\si\in\cals^{\ge 1}(X)$ and every $x_0\in X$, the map $w_{\si,x_0}:\partial_\infty X\to\mathbb{R}_{\ge 0}$ is continuous.
\end{lemma}

\begin{proof}
Notice that for every $\sigma\in\cals^{\ge 1}(X)$, there exists $M\in\mathbb{R}_+$ such that for every geodesic ray $r$ in $X$ originating at $x_0$, and every $t\ge M$, one has $f_{\si}(r(t))=0$. The continuity of the map $w_{\si,x_0}$ therefore follows from the fact that $f_{\si}$ is Lipschitz, together with the fact that for every $\epsilon>0$ and every sequence $(\xi_n)_{n\in\mathbb{N}}\in (\partial_\infty X)^{\mathbb{N}}$ converging to $\xi$, denoting by $r_n$ the geodesic ray from $x_0$ to $\xi_n$ and by $r$ the geodesic ray from $x_0$ to $\xi$, then for all sufficiently large $n\in\mathbb{N}$, the segments $r_n([0,M])$ and $r([0,M])$ are at Hausdorff distance at most $\epsilon$ from each other.
\end{proof}

As a consequence of Lemmas~\ref{lemma:list-of-simplices-borel} and~\ref{lemma:weight-continuous}, we obtain the following fact.

\begin{cor}\label{cor:probability-borel}
For every $x_0\in X$ and every $n\in\mathbb{N}$, the map $\nu_{n,x_0}:\partial_\infty X\to\Prob(\cals(X))$ is Borel.
\qed
\end{cor}

We are now in position to complete our proof of Proposition~\ref{prop:proba-on-edges}.

\begin{proof}[Proof of Proposition~\ref{prop:proba-on-edges}]
The maps $\nu_{n,x_0}:\partial_\infty X\to\Prob(\cals^{\ge 1}(X))$ are Borel (Corollary~\ref{cor:probability-borel}). They are asymptotically $G$-equivariant because for every $g\in G$ and every $\xi\in\partial_\infty X$, one has $$||\nu_{n,x_0}(g\xi)-g\cdot\nu_{n,x_0}(\xi)||_1=||\nu_{n,x_0}(g\xi)-\nu_{n,gx_0}(g\xi)||_1,$$ which converges to $0$ as $n$ goes to $+\infty$ in view of Lemma~\ref{lemma:asymptotically-equivariant}.
\end{proof}

\subsection{Specification to Artin groups}

 We now reach the first main theorem of this paper for Artin groups.

\begin{theo}\label{theo:exactness-Artin}
Every $2$-dimensional Artin group of hyperbolic type is boundary amenable.
\end{theo}

\begin{proof}
Let $G=G_\Gamma$ be a $2$-dimensional Artin group of hyperbolic type. Recall from Section~\ref{sec:background} or \cite[Section~3.1]{MP} that the modified Deligne complex $\del$ can be equipped with the structure of a $\cat$ simplicial complex with finitely many isometry types of simplices, so that $G$ acts discretely and by  simplicial isometries on $\del$. Vertex stabilizers are boundary amenable: indeed, the stabilizer of a rank $0$ vertex is trivial, the stabilizer of a rank $1$ vertex is infinite cyclic, and the stabilizer of a rank $2$ vertex is isomorphic to an Artin group whose defining graph os an edge, so it is boundary amenable as the quotient by its center is virtually free, see \cite[Section~2]{crisp}. The conclusion therefore follows from Theorem~\ref{theo:exactness}.
\end{proof}

We also attain the following statement, which will be useful in the sequel of the paper.

\begin{theo}\label{theo:amenability-boundary-artin}
Let $G=G_\Gamma$ be a $2$-dimensional Artin group of hyperbolic type with defining graph $\Gamma$. Then the $G$-action on the Gromov boundary of the modified Deligne complex $\del$ is Borel amenable. 
\end{theo}

\begin{proof}
This follows from Theorem~\ref{theo:amenability-boundary}, using the fact that edge stabilizers for the $G$-action on $\del$ are amenable (since every edge contains a vertex of rank at most $1$, whose stabilizer is trivial or infinite cyclic). 
\end{proof}

\section{Measured groupoids and measure equivalence rigidity}\label{sec:background-me}

\emph{We now review standard techniques for proving measure equivalence rigidity theorems from consdering self measure equivalence couplings, and provide some background on Borel and measured groupoids and their use in measured group theory.}

\subsection{Universal factors and measure equivalence rigidity}
 We first recall the notion of a measure equivalence coupling between two countable groups, due to Gromov \cite{Gro}, from the introduction. Recall that a \emph{standard measure space} $\Sigma$ is a measurable space associated to a Polish space (i.e.\ a complete separable metric space), equipped with a $\sigma$-finite positive Borel measure. Given an action of a group $G$ on $\Sigma$, we say that a Borel subset $Y\subseteq\Sigma$ is a \emph{fundamental domain} for the $G$-action if $\Sigma$ is the disjoint union of all translates $gY$ with $g\in G$ (with this definition, the existence of a fundamental domain implies the freeness of the action).

\begin{de}[Measure equivalence coupling]
Let $G_1$ and $G_2$ be two countable groups. A \emph{measure equivalence coupling} between $G_1$ and $G_2$ is a standard measure space $\Sigma$ equipped with a measure-preserving action of $G_1\times G_2$ by Borel automorphisms such that for every $i\in\{1,2\}$, the $G_i$-action on $\Sigma$ is free and has a finite measure fundamental domain.
\end{de}

 When $G_1=G_2$, we simply say that $\Sigma$ is a \emph{self measure equivalence coupling} of $G_1=G_2$. There are slight variations over the above definition in the literature, where one only imposes the action to be essentially free and the fundamental domains to be defined ``up to measure $0$'', as in Kida's work \cite[Definition~1.1]{Kid-survey}. However, up to replacing the coupling by a conull invariant subset, one can easily transform a coupling in the sense of Kida to a coupling in the above sense.

\begin{de}[Universal $(G_1,G_2)$-factor]
Let $G_1$ and $G_2$ be two countable discrete groups, and let $\Omega$ be a standard Borel space equipped with an action of  $G_1\times G_2$ by Borel automorphisms. 

The Borel space $\Omega$ is a \emph{universal $(G_1,G_2)$-factor} if for every measure equivalence coupling $\Sigma$ between $G_1$ and $G_2$, there exists a Borel map $\Psi:\Sigma\to\Omega$ which is almost $(G_1\times G_2)$-equivariant, i.e.\ such that for all $(g_1,g_2)\in G_1\times G_2$ and a.e.\ $y\in\Sigma$, one has $\Psi((g_1,g_2)y)=(g_1,g_2)\Psi(y)$.
\end{de}

When $G_1=G_2$, Theorem~\ref{theo:bfs} below justifies  (taking $\Omega=\mathsf{G}$) why finding universal factors is useful for proving measure equivalence rigidity results. Following \cite[Definition~2.2]{BFS}, we say that a Polish group $\mathsf{G}$ is \emph{strongly ICC} if there is no conjugation-invariant probability measure on  $\mathsf{G}\setminus\{e\}$. When $\mathsf{G}$ is a countable discrete group, this is equivalent to the more classical \emph{ICC} (for \emph{infinite conjugacy classes}) condition, requiring that every nontrivial conjugacy class of  $\mathsf{G}$ is infinite. 

\begin{theo}[Bader--Furman--Sauer \cite{BFS}]\label{theo:bfs}
 Let $\mathsf{G}$ be a locally compact second countable group, and let $\Lambda$ be a lattice in $\mathsf{G}$. Assume that $\mathsf{G}$ is strongly ICC, and that for every self measure-equivalence coupling $\Sigma$ of $\Lambda$, there exists a $(\Lambda\times \Lambda)$-equivariant Borel map $\Sigma\to \mathsf{G}$, where the action of $\Lambda\times \Lambda$ on $\mathsf{G}$ is via $(\lambda_1,\lambda_2)\cdot g=\lambda_1 g \lambda_2^{-1}$.

Then for every countable group $\Lambda'$ which is measure equivalent to $\Lambda$, there exists a homomorphism $\rho:\Lambda'\to \mathsf{G}$ with finite kernel whose image is a lattice in $\mathsf{G}$. 
\end{theo}

\begin{proof}
We include a short proof, which is essentially a definition chase, to explain why this theorem is a consequence of \cite[Theorem~2.6]{BFS}. The notion of measure equivalence extends to an equivalence relation on the set of all unimodular locally compact second countable groups: for this we refer to \cite[0.5.E]{Gro} or \cite[Definition~1.1]{BFS}, and we only mention that a lattice in $\mathsf{G}$ is always measure equivalent to $\mathsf{G}$.  

As $\mathsf{G}$ is strongly ICC, it follows from \cite[Lemma~A.6]{BFS} that $\mathsf{G}$ is also strongly ICC relative to $\Lambda$, i.e.\ the only probability measure on $\mathsf{G}$ which is invariant under conjugation by $\Lambda$ is the Dirac mass at $e$. Our assumption on self-couplings of $\Lambda$ says that $\Lambda$ is taut relative to the inclusion $\Lambda\hookrightarrow \mathsf{G}$ in the sense of \cite[Definition~1.10]{BFS} (see indeed \cite[Lemma~A.8]{BFS}, using the strong ICC property, for the uniqueness of the tautening map $\Phi$ from \cite[Definition~1.10]{BFS}). It thus follows from \cite[Proposition~2.9]{BFS} that $\mathsf{G}$ is taut (in the sense of \cite[Definition~1.3]{BFS}, i.e.\ relative to itself).

Let now $\Lambda'$ be a countable group which is measure equivalent to $\Lambda$, whence to $\mathsf{G}$. We can therefore apply \cite[Theorem~2.6]{BFS} (with  $G=\mathcal{G}=\mathsf{G}$ and $H=\Lambda$) to find a continuous homomorphism $\rho:\Lambda'\to \mathsf{G}$ with compact (whence finite) kernel and closed image, and a nonzero Borel measure $\bar{m}$ on $\mathsf{G}$ which is invariant under the left $\mathsf{G}$-action and the right $\rho(\Lambda')$-action such that $\bar{m}$ descends to a finite measure on $\mathsf{G}/\rho(\Lambda')$. As $\rho(\Lambda')$ is discrete, $\bar{m}$ is locally finite (i.e.\ every point has an open neighborhood with finite $\bar{m}$-measure), hence $\bar{m}$ is also regular (locally finite Borel measures on Polish spaces are regular, see e.g.\ \cite[Theorem~26.3 and Corollary~26.4]{Bau}). Hence $\bar{m}$ is a Haar measure on $\mathsf{G}$, and $\rho(\Lambda')$ is a lattice in $\mathsf{G}$.
\end{proof}

We recall from the introduction that two countable groups $G_1,G_2$ are \emph{almost isomorphic} if there exist finite-index subgroups $G_i^0\subseteq G_i$ and finite normal subgroups $F_i\unlhd G_i^0$ such that $G_1^0/F_1$ is isomorphic to $G_2^0/F_2$.

\begin{cor}\label{cor:bfs}
Let $G$ be a countable group, and let  $\Lambda\subseteq G$ be a finite-index subgroup. Assume that $\Lambda$ is ICC, and that for every self measure-equivalence coupling $\Sigma$ of $\Lambda$, there exists a $(\Lambda\times\Lambda)$-equivariant Borel map $\Sigma\to G$, where the action of $\Lambda\times\Lambda$ on $G$ is via $(\lambda_1,\lambda_2)\cdot g=\lambda_1 g \lambda_2^{-1}$.

Then every countable group $\Lambda'$ which is measure equivalent to $\Lambda$, is in fact almost isomorphic to $\Lambda$. 
\end{cor}

\begin{proof}
Let $K\unlhd G$ be the FC-center of $G$, i.e.\ the characteristic subgroup made of all elements whose conjugacy class is finite. As $\Lambda$ is ICC, we have $K\cap\Lambda=\{1\}$, so $K$ is finite. The quotient $\overline{G}=G/K$ is therefore ICC, because every infinite conjugacy class in $G$ projects to an infinite conjugacy class in $\overline{G}$. Let $\overline{\Lambda}$ be the image of $\Lambda$, a finite-index subgroup of $\overline{G}$ isomorphic to $\Lambda$. The ICC group $\overline{G}$ is also a universal $(\Lambda,\Lambda)$-factor, to which we can apply Theorem~\ref{theo:bfs} to get a homomorphism $\Lambda'\to\overline{G}$ with finite kernel $F$ and finite-index image. Let  $\Lambda''\subseteq\Lambda'$ be a finite-index subgroup whose image is contained in $\overline{\Lambda}$. Then $\Lambda''/(F\cap\Lambda'')$ is isomorphic to a finite-index subgroup of $\Lambda$.
\end{proof}

\subsection{Background on measured groupoids} 

We recall that a \emph{standard Borel space} is a Borel space associated to a Polish space. 

\begin{de}[Borel groupoid]
Let $Y$ be a standard Borel space. A \emph{Borel groupoid over $Y$} is a standard Borel space $\calg$, together with the following data:
\begin{enumerate} 
\item Borel maps $s,r:\calg\to Y$ (called the \emph{source map} and the \emph{range map}, respectively); 
\item a Borel map $*$ from $\{(g_1,g_2)\in\calg\times\calg|s(g_1)=r(g_2)\}$ to $\calg$ such that for every $g_1,g_2,g_3\in\calg$, the following two facts hold: 
\begin{enumerate}
\item if $s(g_1)=r(g_2)$, then $s(g_1\ast g_2)=s(g_2)$ and $r(g_1*g_2)=r(g_1)$, and 
\item if $s(g_1)=r(g_2)$ and $s(g_2)=r(g_3)$, then $(g_1*g_2)*g_3=g_1*(g_2*g_3)$;
\end{enumerate} 
\item a Borel map $e:Y\to\calg$ such that for every $y\in Y$ and every $g\in\calg$, the following three facts hold:
\begin{enumerate}
\item one has $s(e(y))=r(e(y))=y$, 
\item if $s(g)=y$, then $g\ast (e(y))=g$,
\item if $r(g)=y$, then $(e(y))\ast g=g$; 
\end{enumerate}
\item a Borel map $\iota:\calg\to\calg$ such that for every $g\in\calg$, the following four facts hold:
\begin{enumerate} 
\item $s(\iota(g))=r(g)$, 
\item $r(\iota(g))=s(g)$,
\item $g\ast\iota(g)=e(r(g))$,
\item $\iota(g)\ast g=e(s(g))$.
\end{enumerate}
\end{enumerate}
\end{de}

The space $Y$ is called the \emph{base space} of the groupoid. We will only consider Borel groupoids which are \emph{discrete}, i.e.\ such that for every $y\in Y$, there are at most countably many $g\in \calg$ with $s(g)=y$. 
The groupoid $\calg$ is \emph{trivial} if $\calg=\{e(y)|y\in Y\}$.

\begin{de}[Borel subgroupoid]
Let $\calg$ be a Borel groupoid over a base space $Y$. A \emph{Borel subgroupoid} of $\calg$ is a Borel subset $\calh\subseteq\calg$ which satisfies the following properties: 
\begin{enumerate}
\item for every $(h_1,h_2)\in \calh^2$ such that $s(h_1)=r(h_2)$, one has $h_1\ast h_2\in\calh$,
\item for every $h\in\calh$, one has $\iota(h)\in\calh$, and
\item for every $y\in Y$, one has $e(y)\in\calh$.
\end{enumerate}
\end{de}

Every Borel subgroupoid of $\calg$ naturally comes with the structure of a Borel groupoid over the same base space $Y$, induced by that of $\calg$. We will always endow subgroupoids with this induced structure.

\begin{de}[Restriction of a groupoid]
Let $\calg$ be a Borel groupoid over a base space $Y$, and let $A\subseteq Y$ be a Borel subset. Then the groupoid structure on $\calg$ naturally turns $$\calg_{|A}:=\{g\in\calg|s(g),r(g)\in A\}$$ into a Borel groupoid over the base space $A$, which is called the \emph{restriction} of $\calg$ to $A$. 
\end{de} 
 
\paragraph*{Measured groupoids.} Let $\calg$ be a discrete Borel groupoid over a base space $Y$. If $Y$ carries a $\sigma$-finite positive Borel measure $\mu$ (in which case we say that $(Y,\mu)$ is a \emph{standard measure space}), we say that $\mu$ is \emph{$\calg$-quasi-invariant} if for every Borel subset $B\subseteq \calg$ such that $s_{|B}$ and $r_{|B}$ are Borel isomorphisms to Borel subsets of $Y$, one has $\mu(s(B))=0$ if and only if $\mu(r(B))=0$. In this case we say that $\calg$ is a \emph{measured groupoid} over $(Y,\mu)$. Every Borel subgroupoid of a measured groupoid naturally comes with the structure of a measured groupoid, and likewise the restriction of a measured groupoid to a Borel subset of positive measure naturally comes with the structure of a measured groupoid. Notice that for every measure $\mu$ as above, there always exists a finite measure $\mu'$ with the same null sets as $\mu$, so without loss of generality we will often assume in the sequel that $\mu$ has finite measure.

A measured groupoid $\calg$ over $(Y,\mu)$ is \emph{of infinite type} if for every Borel subset $Y'\subseteq Y$ of positive measure, and a.e.\ $y\in Y'$, there are infinitely many $g\in\calg_{|Y'}$ such that $s(g)=y$.

\paragraph*{Stable containment, stable equivalence.} Let $\calg$ be a measured groupoid over a base space $Y$, and let $\calh$ and $\calh'$ be two measured subgroupoids of $\calg$. We say that $\calh$ is \emph{stably contained} in $\calh'$ if there exist a conull Borel subset $Y^*\subseteq Y$ and a partition $Y^*=\dunion_{i\in I} Y_i$ into at most countably many Borel subsets, such that for every $i\in I$, we have $\calh_{|Y_i}\subseteq\calh'_{|Y_i}$. We say that $\calh$ and $\calh'$ are \emph{stably equivalent} if each is stably contained in the other (equivalently, there exists a partition $Y^*=\dunion_{i\in I}Y_i$ as above such that for every $i\in I$, one has $\calh_{|Y_i}=\calh'_{|Y_i}$).

\paragraph*{Cocycles.} 

Let $\calg$ be a measured groupoid over a base space $Y$, and let $G$ be a countable discrete group. A Borel map $\rho:\calg\to G$ is a \emph{cocycle} if there exists a conull Borel subset $Y^*\subseteq Y$ such that for all $g_1,g_2\in \calg_{|Y^*}$, if $s(g_1)=r(g_2)$, then $\rho(g_1*g_2)=\rho(g_1)\rho(g_2)$. We say that the cocycle $\rho$ is \emph{strict} if the above identity holds for all $g_1,g_2\in\calg$ with $s(g_1)=r(g_2)$, without passing to a conull subset $Y^*$. 

The \emph{kernel} of a strict cocycle $\rho$ is the subgroupoid of $\calg$ made of all elements $g\in\calg$ such that $\rho(g)=1_G$ (the fact that it is a subgroupoid follows from the fact that the cocycle is strict, otherwise one needs to define the kernel as a subgroupoid of $\calg_{|Y^*}$). A strict cocycle $\rho$ is \emph{action-type} if its kernel is trivial, and for every infinite subgroup $H\subseteq G$, the subgroupoid $\rho^{-1}(H)$ is of infinite type (once again, the fact that $\rho^{-1}(H)$ is a subgroupoid follows from the fact that the cocycle $\rho$ is strict).  

Let now $\Omega$ be a Polish space equipped with an action of $G$ by Borel automorphisms. A Borel map $\phi:Y\to\Omega$ is \emph{$(\calg,\rho)$-equivariant} if there exists a conull Borel subset $Y^*\subseteq Y$ such that for all $g\in\calg_{|Y^*}$, one has $\phi(r(g))=\rho(g)\phi(s(g))$. 
We say that an element $\omega\in\Omega$ is \emph{$(\calg,\rho)$-invariant} if the constant map with value $\omega$ is $(\calg,\rho)$-equivariant; in other words, there exists a conull Borel subset $Y^*\subseteq Y$ such that $\rho(\calg_{|Y^*})\subseteq\Stab_G(\omega)$. 

Slightly abusing notation, in the context of restrictions of groupoids, we will often denote a cocycle and its restriction in the same way. For example, we will write that a Borel map $U\to\Omega$ is $(\calg_{|U},\rho)$-equivariant when we should really write $(\calg_{|U},\rho_{|\calg_{|U}})$-equivariant. 

Given a Borel map $\phi:Y\to\Omega$ and a strict cocycle $\rho:\calg\to G$ to a countable group $G$, the \emph{$(\calg,\rho)$-stabilizer} of $\phi$ is the subgroupoid of $\calg$ made of all $g\in\calg$ such that $\phi(r(g))=\rho(g)\phi(s(g))$ -- notice that the fact that this is a subgroupoid of $\calg$ uses the assumption that the cocycle $\rho$ is strict.

\paragraph*{Example: Restrictions of group actions.}
The most important example to keep in mind for the present paper is the following. Let $(Y,\nu)$ be a standard measure space, and let $G$ be a countable group acting by Borel automorphisms in a quasi-measure-preserving way on $(Y,\nu)$ (i.e.\ for every $g\in G$ and every Borel subset $A\subseteq Y$, one has $\nu(A)=0$ if and only if $\nu(gA)=0$). Then the direct sum $G\times Y$, equipped with the product of the counting measure on $G$ and the measure $\nu$ on $Y$, has a natural structure of a measured groupoid over the base space $Y$, with $s(g,y)=y$ and $r(g,y)=gy$, where the composition law is given by $(g_1,g_2 y)\ast (g_2,y)=(g_1g_2,y)$, the identity map is given by $e(y)=(1_G,y)$, and the inverse map is given by $\iota(g,y)=(g^{-1},gy)$. This groupoid is denoted $G\ltimes Y$. The map $\rho:G\ltimes Y\to G$ defined by letting $\rho(g,y)=g$ is a strict cocycle. When $\nu$ is a finite measure and the $G$-action on $Y$ is measure-preserving, the cocycle $\rho$ is action-type by \cite[Proposition~2.26]{Kid-survey}. 

More generally, if $A\subseteq Y$ is a Borel subset of positive measure, then we can consider the groupoid $(G\ltimes Y)_{|A}$, which also comes with a natural cocycle towards $G$.

\paragraph*{Amenable groupoids.} 

Given a separable real Banach space $B$, we denote by $B^\ast$ the dual Banach space of $B$ (with the weak-$*$ topology), by $B_1^*$ the closed unit ball of $B^*$, by $\Conv(B_1^\ast)$ the set of all nonempty closed convex subsets of $B_1^*$, and by $\Isom(B)$ the group of all linear isometries of $B$, endowed with the strong operator topology. The adjoint operator of every isometry $T\in\Isom(B)$ preserves $B_1^*$; we denote by $T^\ast$ the restriction to $B_1^\ast$ of the adjoint operator of $T$. 

Let now $\calg$ be a measured groupoid over a base space $Y$, and let $\rho:\calg\to\Isom(B)$ be a cocycle. The \emph{adjoint cocycle} of $\rho$ is the cocycle $\rho^*:\calg\to\Isom(B_1^*)$ defined by letting $\rho^{\ast}(g):=(\rho(g)^{-1})^\ast$ for every $g\in\calg$. Given a separable real Banach space $B$, a \emph{measurable $B$-convex field over $Y$} is a map $K:Y\to\Conv(B_1^*)$ such that $$\{(y,k)\in Y\times B_1^\ast| k\in K(y)\}$$ is Borel. A \emph{section} of $K$ is a Borel map $s:Y\to B_1^\ast$ such that $s(y)\in K(y)$ for a.e.\ $y\in Y$. 

The following notion is a generalization to the groupoid setting of Zimmer's notion of an amenable group action \cite{Zim}.

\begin{de}[Amenable measured groupoid]
A measured groupoid $\calg$ over a base space $Y$ is \emph{amenable} if for every separable real Banach space $B$ and every cocycle $\rho:\calg\to\Isom(B)$, every $(\calg,\rho^\ast)$-equivariant measurable $B$-convex field over $Y$ has a  $(\calg,\rho^\ast)$-equivariant section. 
\end{de}

\paragraph*{Normalized subgroupoids.} 
Let $\calg$ be a measured groupoid over a base space $Y$, and let $\calh$ be a measured subgroupoid of $\calg$. Given a Borel subset $B\subseteq\calg$, we say that $\calh$ is \emph{$B$-invariant} if there exists a conull Borel subset $Y^*\subseteq Y$ such that for all $g_1,g_2\in B\cap\calg_{|Y^*}$ and all $h\in\calg_{|Y^*}$ with $s(h)=s(g_1)$ and $r(h)=s(g_2)$, one has $h\in\calh$ if and only if $g_2*h*\iota(g_1)\in\calh$. 

\begin{de}[Normalized subgroupoid]
Let $\calg$ be a measured groupoid over a base space $Y$, and let $\calh$ and $\calh'$ be measured subgroupoids of $\calg$.

The subgroupoid $\calh$ is \emph{normalized} by $\calh'$ if $\calh'$ can be written as a union of countably many Borel subsets $B_n$ of $\calg$ in such a way that $\calh$ is $B_n$-invariant for every $n\in\mathbb{N}$.  

The subgroupoid $\calh$ is \emph{stably normalized} by $\calh'$ if there exists a countable Borel partition $Y=\dunion_{n\in\mathbb{N}}Y_n$ such that $\calh_{|Y_n}$ is normalized by $\calh'_{|Y_n}$ for every $n\in\mathbb{N}$.
\end{de}

\begin{rk}\label{rk-normal}
In view of a theorem of Lusin and Novikov (see \cite[Theorem~18.10]{Kec}), the subsets $B_n$ in the definition of $\calh$ being normalized by $\calh'$ can always be chosen so that $s_{|B_n}$ and $r_{|B_n}$ are Borel isomorphisms onto Borel subsets $s(B_n)$ and $r(B_n)$.
\end{rk}

\subsection{From measured groupoids to universal factors}

The following proposition explains how measured groupoids naturally arise in the study of measure equivalence between countable groups, and gives a criterion in terms of measured groupoids to check that a space is a universal $(G_1,G_2)$-factor. The argument relies on a construction from \cite[Section~3]{Fur-oe}.

\begin{prop}\label{prop:me-groupoid}
Let $G_1$ and $G_2$ be two countable discrete groups which are measure equivalent. 
\begin{enumerate}
\item There exist a measured groupoid $\calg$ over a finite measure space $Y$, and two strict action-type cocycles $\calg\to G_1$ and $\calg\to G_2$.
\item Let $\Omega$ be a standard Borel space, and assume that for every measured groupoid $\calg$ over a finite measure space $Y$, and every pair of strict action-type cocycles $\rho_1:\calg\to G_1$ and $\rho_2:\calg\to G_2$, there exist a Borel map $\phi:Y\to\Omega$ and a conull Borel subset $Y^*\subseteq Y$ such that for all $g\in\calg_{|Y^*}$, one has $\phi(r(g))=(\rho_1(g),\rho_2(g))\phi(s(g))$. Then $\Omega$ is a universal $(G_1,G_2)$-factor.
\end{enumerate}
\end{prop}

\begin{proof} 
We first prove Assertion~1. Let $\Sigma$ be a measure equivalence coupling between $G_1$ and $G_2$.
Every fundamental domain $X$ for the $G_1$-action on $\Sigma$ comes equipped with a natural measure-preserving $G_2$-action, by letting $g_2\cdot x$ be the unique point of $X$ in the same $G_1$-orbit as $g_2x\in\Sigma$. Similarly, every  fundamental domain $Y$ for the $G_2$-action on $\Sigma$ comes equipped with a measure-preserving $G_1$-action defined in a similar way. By \cite[Lemma~2.27]{Kid-survey}, we can choose the Borel fundamental domains $X$ and $Y$ in such a way that $G_1\cdot (X\cap Y)= Y$ and $G_2\cdot (X\cap Y) = X$ up to null sets. We choose $X$ and $Y$ as such, and let $A:=X\cap Y$.

The groupoids $(G_1\ltimes Y)_{|A}$ and $(G_2\ltimes X)_{|A}$ are isomorphic via the map sending $(g_1,y)$ to $(g_2,y)$, where $g_2$ is the unique element of $G_2$ such that $(g_1,g_2)y\in Y$. Letting $\calg:=(G_1\ltimes Y)_{|A}$, we thus have two natural strict action-type cocycles $\rho_i:\calg\to G_i$ for $i\in\{1,2\}$ -- see \cite[Proposition~2.26]{Kid} for the fact that these cocycles are action-type. This concludes the proof of Assertion~1.

To prove Assertion~2, let $\calg$ be the groupoid constructed just above, and $\rho_1:\calg\to G_1$ and $\rho_2:\calg\to G_2$ be the constructed cocycles. By assumption, there exist a Borel map $\phi:A\to\Omega$ and a conull Borel subset $A^*\subseteq A$ such that for all $g\in\calg_{|A^*}$, one has $\phi(r(g))=(\rho_1(g),\rho_2(g))\phi(s(g))$. 

We claim that the assignment $(g_1,g_2)x\mapsto (g_1,g_2)\phi(x)$ (with $x\in A^*$) yields a well-defined Borel map $\psi:\Sigma\to\Omega$ (in fact this map is only defined almost everywhere but can be extended arbitrarily to the whole $\Sigma$), and this map is almost $(G_1\times G_2)$-equivariant. 
To prove the claim, it is enough to check that for all  $x,x'\in A^*$, if $(g_1,g_2)x=(h_1,h_2)x'$, then $(g_1,g_2)\phi(x)=(h_1,h_2)\phi(x')$. Indeed, in this case, we have $(h_1^{-1}g_1)(h_2^{-1}g_2)x=x'$, and it follows that there exists $g\in\calg$ with $\rho_1(g)=h_1^{-1}g_1$, $\rho_2(g)=h_2^{-1}g_2$ and with $s(g)=x$ and $r(g)=x'$. The given equivariance of $\phi$ ensures that $\phi(x')=(h_1^{-1}g_1,h_2^{-1}g_2)\phi(x)$. In other words $(g_1,g_2)\phi(x)=(h_1,h_2)\phi(x')$, as desired. 
\end{proof}

\section{Curve-graph-like graphs and universal factors}\label{sec:axioms}

\emph{In this section, we abstract Kida's arguments from his proof of measure equivalence rigidity of mapping class groups, and give a general framework to prove a rigidity statement regarding self-couplings of a countable group $G$ -- which leads to strong measure equivalence rigidity results through Theorem~\ref{theo:bfs}. This framework requires having a $G$-action on a graph -- analogous to the curve graph -- 
with some algebraic control on vertex and edge stabilizers. It also involves a dynamical setting with a $G$-action on a compact space -- analogous to Thurston's compactification of the Teichmüller space -- and a $G$-action on a space that plays the role of the boundary of the curve graph. This setting will be checked to hold in the realm of $2$-dimensional Artin groups of hyperbolic type in later sections.}

\subsection{Curve-graph-like $G$-graphs and statement of the main result}

Given a group $G$ and a $G$-set $S$, we let $\Stab_G(S)$ be the setwise stabilizer of $S$ in $G$ and $\pstab_G(S)$ be its elementwise stabilizer. We let $\calp_f(S)$ be the set of all nonempty finite subsets of $S$, and $\calp_{\flex}(S)$ be the set of all subsets of $S$ whose elementwise $G$-stabilizer is infinite. Given a Polish space $\Delta$, we denote by $\Delta^{(3)}$ the set of all pairwise distinct triples of $\Delta$, equipped with the topology induced from that of $\Delta^3$. A graph $\Theta$ is \emph{simple} if it contains no edge-loop and no multiple edges between distinct vertices. This is equivalent to every automorphism of $\Theta$ being determined by its action on the vertex set of $\Theta$.

\begin{de}[Curve-graph-like $G$-graphs]\label{de:curve-graph-like}
Let $G$ be a group. A $G$-graph $\Theta$ is \emph{curve-graph-like} if it is simple, has countable vertex set, and there exists a $G$-invariant subset $\calp_\ast(V(\Theta))\subseteq \calp_f(V(\Theta))$ such that the following properties hold.
\begin{enumerate}
\item\textbf{(Stabilizers versus pointwise stabilizers)} For every $S\in\calp_\ast(V(\Theta))$, we have $\stab_{G}(S)=\pstab_{G}(S)$.
\item\textbf{(Dynamical setting)} There exist a compact nonempty metrizable space $K$ equipped with a $G$-action by homeomorphisms, a Borel $G$-invariant partition $K=K_{\bdd}\cup K_\infty$, and a Polish $G$-space $\Delta$ such that
\begin{enumerate}
\item there is a Borel $G$-equivariant map $K_{\bdd}\to \calp_\ast(V(\Theta))$,
\item there is a Borel $G$-equivariant map $K_\infty\to \Delta$,
\item there is a Borel $G$-equivariant map $\Delta^{(3)}\to \calp_\ast(V(\Theta))$,
\item the $G$-action on $\Delta$ is Borel amenable.
\end{enumerate}
\item\textbf{(Canonical reduction systems and chain condition)} The following hold:
\begin{enumerate}
\item there is a $G$-equivariant map $\calp_{\flex}(V(\Theta))\to \calp_\ast(V(\Theta))$,
\item there is a bound on the length of an increasing (for inclusion) chain of subsets $X_i$ of $V(\Theta)$ such that $\pstab_G(X_1)\supsetneq\pstab_G(X_2)\supsetneq\cdots\supsetneq\pstab_G(X_k)$. 
\end{enumerate}
\item\textbf{(Structure of vertex and edge stabilizers of $\Theta$)} The following hold:
\begin{enumerate}
\item for every $v\in V(\Theta)$, the group $\Stab_G(v)$ contains both a nonabelian free subgroup and an infinite amenable normal subgroup,
\item given any two distinct vertices $v,w\in V(\Theta)$, there exists an infinite order element in $G$ that fixes $v$ but not $w$,
\item given any two distinct vertices $v,w\in V(\Theta)$, the following are equivalent:
\begin{enumerate}
\item the vertices $v$ and $w$ are joined by an edge in $\Theta$;
\item for every vertex $v'$ of $\Theta$, if every element of $\Stab_G(v)\cap\Stab_G(w)$ has a power that fixes $v'$, then $v'\in\{v,w\}$.
\end{enumerate}
\end{enumerate} 
\end{enumerate}
\end{de}

\begin{rk}\label{rk:curve-graph-finite-index}
Notice that if a $G$-graph $\Theta$ is curve-graph-like, then for every finite-index subgroup $G'\subseteq G$, the graph $\Theta$ viewed as as $G'$-graph is again a curve-graph-like $G'$-graph (for Assumption~2.(d), the fact that the $G'$-action on $\Delta$ is Borel amenable follows for instance from \cite[Proposition~2.12]{BGH}). 
\end{rk}

In order to motivate this definition, it is worth describing the model case, which comes from mapping class groups of finite-type hyperbolic surfaces acting on their corresponding curve graph.

\begin{ex}[Mapping class groups and curve graphs]\label{ex:mcg}
Let $\Sigma$ be a connected orientable surface of finite type, i.e.\ $\Sigma$ is obtained from a closed surface of finite genus by possibly removing finitely many points. Assume that $\Sigma$ is neither a sphere with at most four punctures, nor a torus with at most one puncture. Let $\Mod(\Sigma)$ be the mapping class group of $\Sigma$, and let $G$ be its finite-index subgroup made of all mapping classes acting trivially on the homology of the surface with $\mathbb{Z}/3\mathbb{Z}$ coefficients. Let $\calc(\Sigma)$ be the curve graph of $\Sigma$, and let $\calp_\ast(V(\calc(\Sigma)))=\calp_f(V(\calc(\Sigma)))$ be the set of all nonempty finite subsets of $V(\calc(\Sigma))$. We will now check that $\calc(\Sigma)$ is a curve-graph-like $G$-graph. We warn the reader that with our definition $\calc(\Sigma)$ is not a curve-graph-like $\Mod(\Sigma)$-graph because Assumption~1 is only true after passing to the finite-index subgroup $G$; but this is harmless for our purpose as commensurable groups are measure equivalent.

Assumption~1 asserts that whenever an element $g\in G$ permutes a finite family of isotopy classes of essential simple closed curves on $\Sigma$, then $g$ fixes each of them: this is true in $G$ (but not in $\Mod(\Sigma)$) by work of Ivanov \cite[Corollary~3.6]{Iva2}. 

Interestingly, Assumption~2 (Dynamical setting) is exactly the condition introduced by Guirardel, Lécureux and the first named author in \cite{GHL} (where it is called \emph{geometric rigidity}) to get cocycle superrigidity results. In the mapping class group setting, it is checked as follows. The space $K$ is the compact space $\PML$ of projective measured laminations on $\Sigma$, and we have $\PML=\PML_\bdd\cup\PML_\infty$, where $\PML_\infty$ consists of the subspace made of \emph{arational} laminations. The space $\Delta$ is the Gromov boundary $\partial_\infty\calc(\Sigma)$. By Klarreich's theorem describing the boundary of the curve graph \cite{Kla}, there is a natural $\Mod(\Sigma)$-equivariant (whence $G$-equivariant) map from $\PML_\infty$ to $\Delta=\partial_\infty\calc(\Sigma)$. Also, to every nonarational lamination $\lambda$, one can canonically associate a nonempty finite set of pairwise disjoint isotopy classes of essential simple closed curves on $\Sigma$, namely, the set of all curves that are either isolated in $\lambda$ or arise as a boundary of a subsurface filled by $\lambda$. Amenability of the mapping class group action (whence of the $G$-action) on $\Delta=\partial_\infty\calc(\Sigma)$ was proved by Hamenstädt in \cite{Ham} and independently -- with a different proof -- by Kida in \cite{Kid-memoir} (see also \cite{PrSa} for a recent short proof). Finally, the map assigning a nonempty finite set of isotopy classes of essential simple closed curves to every triple of pairwise distinct points in $\partial_\infty\calc(\Sigma)$ was constructed by Kida in \cite[Section~4.1.2]{Kid-memoir}.

Assumption~3.(a) is the so-called \emph{canonical reduction system}, and is checked as follows: if $X$ is an infinite collection of isotopy classes of essential simple closed curves on $\Sigma$ whose elementwise stabilizer is infinite, then the curves in $X$ do not fill the whole surface, and therefore one can canonically associate to $X$ a nonempty finite set of isotopy classes of essential simple closed curves, made of all the boundaries of the filled subsurfaces (including annular subsurfaces). This defines a Borel $G$-equivariant map $\calp_{\flex}(V(\calc(\Sigma)))\to\calp_\ast(V(\calc(\Theta)))$, as required. To check Assumption~3.(b), observe that the elementwise stabilizer of $X$ is the subgroup of $G$ made of mapping classes that act trivially on each of the filled subsurfaces. As there is a bound on the size of an increasing (for inclusion) family of isotopy classes of subsurfaces of $\Sigma$, we get the desired chain condition.

Finally, Assumption~4 provides algebraic properties of stabilizers of vertices and edges in $\calc(\Sigma)$. For Assumption~4.(a), notice that the stabilizer of a curve contains a nonabelian free group (the complement of a curve is a subsurface which has a connected component that supports two independent pseudo-Anosov diffeomorphisms because we have ruled out low-complexity surfaces). It also contains the intersection with $G$ of the infinite cyclic subgroup generated by the corresponding Dehn twist as a normal subgroup. Assumption~4.(b) follows from the fact that for every essential simple closed curve $c$, every connected component of $\Sigma\setminus c$ which carries an essential simple closed curve non-isotopic to $c$, also supports a pseudo-Anosov diffeomorphism. For Assumption~4.(c), if two curves are disjoint, then their common stabilizer does not fix any third curve. But if they intersect, then their common stabilizer fixes every curve in the subsurface they fill.
\end{ex}

\begin{rk}
The interest of the setting provided by Definition~\ref{de:curve-graph-like} is that it holds both for mapping class groups, as explained right above, and for $2$-dimensional Artin groups of hyperbolic type, working with Crisp's fixed set graph (see Section~\ref{sec:theta} below) -- thereby revealing striking analogies between these two classes of groups. We hope that it could perhaps also hold for other classes of groups; one natural class to wonder about being higher-dimensional Artin groups of hyperbolic type, though so far we are lacking negative curvature complexes to arrange the dynamical setting. 

In general, we do not expect the assumptions from Definition~\ref{de:curve-graph-like} to always hold as such: one often needs to adjust the characterizations of vertex and edge stabilizers to the specific situation at hand. Still, Definition~\ref{de:curve-graph-like} might serve as a guideline: it is for instance noticeable that the dynamical setting and canonical reduction systems are also key in the recent proof of measure equivalence rigidity for $\Out(F_N)$ by Guirardel and the first named author \cite{GH}, or for the measure equivalence classification statement established in \cite{HH} for right-angled Artin groups. Generalizing this latter case, another situation of interest is that of Artin groups of type FC, where the graph of parabolic subgroups studied by Morris-Wright in \cite{MW} could serve as a curve graph analogue -- using the cubical structure of the modified Deligne complex, as in \cite{HH}, to arrange the dynamical setting.
\end{rk}

Given two simple graphs $\Theta_1$ and $\Theta_2$ with countable vertex sets, we denote by $\mathrm{Maps}(\Theta_1\to\Theta_2)$ the set of all maps from $V(\Theta_1)$ to $V(\Theta_2)$ that preserve adjacency. When $\Theta_1$ and $\Theta_2$ are isomorphic, we denote by $\Isom(\Theta_1\to\Theta_2)$ the set of all maps in $\mathrm{Maps}(\Theta_1\to\Theta_2)$ that are bijective from $V(\Theta_1)$ to $V(\Theta_2)$ and have an inverse in $\mathrm{Maps}(\Theta_2\to\Theta_1)$ (in other words, preserving both adjacency and non-adjacency). The set $\Isom(\Theta_1\to\Theta_2)$ has a natural Polish topology, where a basis of open sets is given by all sets of the form $U_{(v_1,\dots,v_n),(w_1,\dots,w_n)}$, with $v_1,\dots,v_n\in V(\Theta_1)$ and $w_1,\dots,w_n\in V(\Theta_2)$, made of all maps $f$ that send $v_i$ to $w_i$ for all $i\in\{1,\dots,n\}$. We always equip $\Isom(\Theta_1\to\Theta_2)$ with this topology and the corresponding Borel structure. If $\Theta_1$ is equipped with an action of a group $G_1$ and $\Theta_2$ is equipped with an action of a group $G_2$, then $\Isom(\Theta_1\to\Theta_2)$ is equipped with an action of $G_1\times G_2$ defined by letting $(g_1,g_2)\cdot\theta(v)=g_2\theta(g_1^{-1}v)$ for every $\theta\in\Isom(\Theta_1\to\Theta_2)$ and every $v\in\Theta_1$. Given a graph $\Theta$, we denote by $\Aut(\Theta)$ the set of all graph isomorphisms from $\Theta$ to itself, which is naturally a group under the composition law. The goal of the present section is to prove the following theorem. 

\begin{theo}\label{theo:blueprint}
Let $G_1$ and $G_2$ be two countable groups. Assume that for every $i\in\{1,2\}$, there exists a curve-graph-like $G_i$-graph, and let $\Theta_i$ be such a graph. 

If $G_1$ and $G_2$ are measure equivalent, then the graphs $\Theta_1$ and $\Theta_2$ are isomorphic, and $\Isom(\Theta_1\to\Theta_2)$ is a universal $(G_1,G_2)$-factor. 
\end{theo}

\begin{rk}
This implies in particular that if $G$ is a countable group, and if a curve-graph-like $G$-graph exists, then it is unique up to graph isomorphism.
\end{rk}

In particular, when $G_1=G_2$, Theorem~\ref{theo:blueprint} provides a universal factor for self-couplings of the group. In view of Theorem~\ref{theo:bfs}, the above theorem is thus a key statement towards proving measure equivalence rigidity results. In particular, Kida's proof of the measure equivalence superrigidity of mapping class groups essentially consists of combining Theorems~\ref{theo:bfs} and~\ref{theo:blueprint}, applied to the setting recalled in Example~\ref{ex:mcg}, together with Ivanov's theorem that the natural map from the extended mapping class group to the automorphism group of the curve graph is almost an isomorphism \cite{Iva} in this setting. In later sections, we will show how to also apply this setting to $2$-dimensional Artin groups of hyperbolic type; in particular, we will need to describe an analogue of the curve graph, already considered by Crisp in \cite{crisp}, in Section~\ref{sec:theta}, and check that it satisfies the assumptions from Definition~\ref{de:curve-graph-like}. Thus, Theorem~\ref{theo:blueprint} provides a common framework for mapping class groups and Artin groups.  

\subsection{Reducible subgroupoids}

From now on, we will adopt all notations from the definition of a curve-graph-like graph (Definition~\ref{de:curve-graph-like}). 

\begin{de}[Reducibility of a cocycle]
Let $G$ be a countable group, and let $\Theta$ be a curve-graph-like $G$-graph. Let $\calg$ be a discrete measured groupoid over a finite measure space $Y$, and let $\rho:\calg\to G$ be a strict cocycle. We say that $(\calg,\rho)$ is 
\begin{enumerate}
\item \emph{reducible} if some vertex of $\Theta$ is $(\calg,\rho)$-invariant;
\item \emph{stably reducible} if there exists a partition $Y=\dunion_{i\in I}Y_i$ into at most countably many Borel subsets such that $(\calg_{|Y_i},\rho)$ is reducible for every $i\in I$;
\item \emph{nowhere reducible} if for every Borel subset $U\subseteq Y$ of positive measure, no vertex of $\Theta$ is $(\calg_{|U},\rho)$-invariant.
\end{enumerate} 
\end{de}

The terminology comes from the analogy with the mapping class group setting, where Kida defined an irreducible subgroupoid coming with a cocycle to the mapping class group as a subgroupoid that does not admit any invariant isotopy class of essential simple closed curve \cite{Kid-memoir} -- generalizing the notion of an irreducible subgroup of the mapping class group, which is one that does not virtually fix the isotopy class of any essential simple closed curve on the surface. 

\begin{lemma}\label{lemma:finite-sets}
Let $G$ be a countable group, and let $\Theta$ be a curve-graph-like $G$-graph. Let $\calg$ be a discrete measured groupoid over a finite measure space $Y$, and let $\rho:\calg\to G$ be a strict cocycle. Assume that $(\calg,\rho)$ is nowhere reducible.

Then for every Borel subset $U\subseteq Y$ of positive measure, no element of $\calp_\ast(V(\Theta))$ is $(\calg_{|U},\rho)$-invariant.
\end{lemma}

\begin{proof}
Let $X\in\calp_\ast(V(\Theta))$, and assume towards a contradiction that $X$ is $(\calg_{|U},\rho)$-invariant for some Borel subset $U\subseteq Y$ of positive measure. This means that there exists a conull Borel subset $U^*\subseteq U$ such that for all $g\in\calg_{|U^*}$, one has $\rho(g)\in\Stab_G(X)$. Using Assumption~1 from the definition of a curve-graph-like $G$-graph, we have $\Stab_G(X)=\pstab_G(X)$. Therefore, for every $x\in X$, the constant map with value $x$ is also $(\calg_{|U},\rho)$-invariant, contradicting that $(\calg,\rho)$ is nowhere reducible. 
\end{proof}

\begin{lemma}\label{lemma:reducible-vs-irreducible}
Let $G$ be a countable group, let $\Theta$ be a curve-graph-like $G$-graph. Let $\calg$ be a discrete measured groupoid over a finite measure space $Y$, and let $\rho:\calg\to G$ be a strict cocycle. Let $\calh$ be a measured subgroupoid of $\calg$. 

Then there exists a Borel partition $Y=Y_0\dunion Y_1$ such that $(\calh_{|Y_0},\rho)$ is nowhere reducible, while $(\calh_{|Y_1},\rho)$ is stably reducible. 
\end{lemma}

\begin{proof}
Let $Y_1\subseteq Y$ be a Borel subset of maximal measure such that $(\calh_{|Y_1},\rho)$ is stably reducible: this exists because if $(U_n)_{n\in\mathbb{N}}$ is a measure-maximizing sequence of such subsets, letting $U=\bigcup_{n\in\mathbb{N}}U_n$, one sees that $(\calh_{|U},\rho)$ is again stably reducible. Letting $Y_0=Y\setminus Y_1$ proves the lemma: indeed, for every $V\subseteq Y_0$ of positive measure, $(\calg_{|V},\rho)$ is not reducible, as otherwise $Y'_1=Y_1\cup V$ would be a subset of larger measure than $Y_1$ such that $(\calh_{|Y'_1},\rho)$ is stably reducible.
\end{proof}

We now turn to analysing stably reducible subgroupoids; our next proposition associates a canonical vertex in $V(\Theta)$ to every such subgroupoid, after a countable Borel partition of the base space.

\begin{prop}\label{prop:crs}
Let $G$ be a countable group, and let $\Theta$ be a curve-graph-like $G$-graph. Let $\calg$ be a discrete measured groupoid over a finite measure space $Y$, and let $\rho:\calg\to G$ be a strict cocycle with trivial kernel. Let $\calh$ and $\calh'$ be two measured subgroupoids of $\calg$, and assume that $\calh$ is of infinite type and is normalized by $\calh'$. Assume that $(\calh,\rho)$ is stably reducible.

Then there exists a partition $Y=\dunion_{i\in I}Y_i$ into at most countably many Borel subsets such that for every $i\in I$, there exists $v_i\in V(\Theta)$ which is both $(\calh_{|Y_i},\rho)$-invariant and $(\calh'_{|Y_i},\rho)$-invariant.
\end{prop} 

\begin{proof}
Up to a first countable partition of the base space, we can assume that $(\calh,\rho)$ is reducible. Let $X\subseteq V(\Theta)$ be the nonempty subset made of all vertices that are $(\calh,\rho)$-invariant. Arguing as in the proof of Lemma~\ref{lemma:reducible-vs-irreducible}, we can find a first partition $Y=U_0\dunion (\dunion_{j\in J}U_j)$ into at most countably many Borel subsets such that
\begin{enumerate}
\item for every Borel subset $U\subseteq U_0$ of positive measure, every vertex of $\Theta$ which is $(\calh_{|U},\rho)$-invariant belongs to $X$, and
\item for every $j\in J$, the set $X_{j}\subseteq V(\Theta)$ made of all vertices which are $(\calh_{|U_{j}},\rho)$-invariant contains $X$ strictly.
\end{enumerate} 
 Observe that for every $j\in J$, we have $\pstab_G(X_j)\subsetneq\pstab_G(X)$: this is because there exists a conull Borel subset $Y^*\subseteq Y$ such that $\rho(\calh_{|Y^*})\subseteq\pstab_G(X)$, while $\rho(\calh_{|Y^*})\nsubseteq\pstab_G(X_i)$. 

We then repeat the same argument as above inductively on each of the countably many subsets $U_j$. The chain condition given by Assumption~$3.(b)$ from Definition~\ref{de:curve-graph-like} ensures that this process stops after finitely many iterations, and in the end we reach a partition $Y=\dunion_{i\in I} Y_i$ into at most countably many Borel subsets such that for every $i\in I$, there exists a set $X_i\subseteq V(\Theta)$ with the following two properties:
\begin{enumerate}
\item every vertex in $X_i$ is $(\calh_{|Y_i},\rho)$-invariant,
\item for every Borel subset $U\subseteq Y_i$ of positive measure, no vertex in $V(\Theta)\setminus X_i$ is $(\calh_{|U},\rho)$-invariant. 
\end{enumerate} 
For every $i\in I$, we now use the canonicity of the set $X_i$ and the fact that $\calh_{|Y_i}$ is normalized by $\calh'_{|Y_i}$ to derive that $X_i$ is (globally, not necessarily pointwise) $(\calh'_{|Y_i},\rho)$-invariant -- in other words, there exists a conull Borel subset $Y_i^*\subseteq Y_i$ such that $\rho(\calh'_{|Y_i^*})\subseteq\Stab_G(X_i)$.

Indeed, since $\calh_{|Y_i}$ is normalized by $\calh'_{|Y_i}$, we can use Remark~\ref{rk-normal} to write $\calh'_{|Y_i}$ as a countable union of Borel subsets $B_n\subset \calh'_{|Y_i}$, each defining an isomorphism between two Borel subsets $U_n,V_n\subseteq Y_i$, so that for every $h\in\calg_{|Y_i}$ and all $\phi,\psi\in B_n$, assuming that the composition $\phi\ast h \ast\iota(\psi)$ is well-defined, 
one has $h\in\calh$ if and only if
$\phi\ast h\ast \iota(\psi)\in \calh$.
By subdividing the sets $B_n$ if needed, we can assume in addition that $\rho(B_n)$ is a single element $g_n\in G$.

Let $n\in\mathbb{N}$ be such that $U_n$ and $V_n$ have positive measure. Every vertex $v\in X_i$ is $(\calh_{|U_n},\rho)$-invariant, so $g_nv$ is $(\calh_{|V_n},\rho)$-invariant. The second defining property of $X_i$ thus ensures that $g_nv\in X_i$. Therefore $g_n\in\Stab_G(X_i)$.
Let now $Y_i^*\subseteq Y_i$ be a conull Borel subset which intersects all subsets $U_n$ and $V_n$ of measure zero trivially.
We then have $\rho(\calh'_{|Y_i^*})\subseteq\Stab_G(X_i)$, as desired.

Next, we observe that for every $i\in I$, we have $X_i\in\calp_{\flex}(V(\Theta))$. Indeed, by definition of $X_i$, there exists a conull Borel subset $Y_i^*\subseteq Y_i$ such that $\rho(\calh_{|Y_i^*})\subseteq\pstab_G(X_i)$. As $\rho$ has trivial kernel and $\calh$ is of infinite type, this shows that $\pstab_G(X_i)$ is infinite, i.e.\ $X_i\in\calp_{\flex}(V(\Theta))$.

By Assumption~$3.(a)$ from Definition~\ref{de:curve-graph-like}, we can therefore canonically associate a nonempty finite subset $F_i\subseteq V(\Theta)$ to the set $X_i$. It follows that $F_i$ is both $(\calh_{|Y_i},\rho)$-invariant and $(\calh'_{|Y_i},\rho)$-invariant. Finally, using Assumption~$1$ from Definition~\ref{de:curve-graph-like}, every vertex in $F_i$ is in fact both $(\calh_{|Y_i},\rho)$-invariant and $(\calh'_{|Y_i},\rho)$-invariant, thereby concluding our proof.  
\end{proof}

\subsection{Irreducible amenable subgroupoids}

The goal of this section is to prove the following proposition, using an idea that dates back to work of Adams \cite{Ada} and was implemented in the mapping class group setting by Kida \cite{Kid-memoir}.

\begin{prop}\label{prop:irreducible-amenable}
Let $G$ be a countable group, and assume that there exists a curve-graph-like $G$-graph. Let $\calg$ be a measured groupoid over a finite measure space $Y$, and let $\rho:\calg\to G$ be a strict action-type cocycle. Let $\calh$ and $\calh'$ be two measured subgroupoids of $\calg$, with $\calh$ of infinite type and normalized by $\calh'$. 

If $(\calh,\rho)$ is nowhere reducible and $\calh$ is amenable, then $\calh'$ is amenable.
\end{prop}

In all statements below, whenever $X$ is a Polish space (which includes the case where $X$ is a countable discrete space), we denote by $\Prob(X)$ the space of all probability measures on $X$, equipped with the topology generated by the maps $\mu\mapsto\int_X fd\mu$, where $f$ varies over all real-valued bounded continuous functions on $X$. By \cite[Theorem~17.24]{Kec}, the Borel $\sigma$-algebra on $\Prob(X)$ is then the $\sigma$-algebra generated by the maps $\mu\mapsto\mu(A)$, where $A$ varies over the Borel subsets of $X$. Consequently, every Borel map $f:X\to Y$ between Polish spaces induces a Borel map $\Prob(X)\to\Prob(Y)$.

\begin{lemma}\label{lemma:no-invariant}
Let $G$ be a countable group, and let $\Theta$ be a curve-graph-like $G$-graph. Let $\calg$ be a measured groupoid over a finite measure space $Y$, and let $\rho:\calg\to G$ be a strict cocycle with trivial kernel. Let $\calh$ be a measured subgroupoid of $\calg$ of infinite type.

If $(\calh,\rho)$ is nowhere reducible, then for every Borel subset $U\subseteq Y$ of positive measure, there is no $(\calh_{|U},\rho)$-equivariant Borel map $U\to\Prob(\calp_\ast(V(\Theta)))$.
\end{lemma}

\begin{proof}
Assume towards a contradiction that there exist a Borel subset $U\subseteq Y$ of positive measure and an $(\calh_{|U},\rho)$-equivariant Borel map $U\to\Prob(\calp_\ast(V(\Theta)))$. As $\calp_\ast(V(\Theta))$ is contained in $\calp_f(V(\Theta))$, it is countable. Therefore, there exists a $G$-equivariant Borel map $\Prob(\calp_\ast(V(\Theta)))\to\calp_f(V(\Theta))$, sending a probability measure $\mu$ to the union of all subsets $X\in\calp_\ast(V(\Theta))$ such that $\mu(X)$ is maximal. We thus get an $(\calh_{|U},\rho)$-equivariant Borel map $U\to\calp_f(V(\Theta))$. Since $\calp_f(V(\Theta))$ is countable, we can find a Borel subset $V\subseteq U$ of positive measure where this map is constant, and we then have an $(\calh_{|V},\rho)$-invariant finite subset $X\in\calp_f(V(\Theta))$. This means that (up to replacing $V$ by a conull Borel subset) one has $\rho(\calh_{|V})\subseteq\Stab_G(X)$. As $\rho$ has trivial kernel and $\calh$ is of infinite type, we deduce that $X\in\calp_{\flex}(V(\Theta))$. Using the $G$-equivariant map $\calp_{\flex}(V(\Theta))\to \calp_\ast(V(\Theta))$ coming from Assumption~3.(a) of the definition of a curve-graph-like $G$-graph, we deduce an $(\calh_{|V},\rho)$-invariant element of $\calp_\ast(V(\Theta))$. As $(\calh,\rho)$ is nowhere reducible, this contradicts Lemma~\ref{lemma:finite-sets}.
\end{proof}

In the following lemma, the space $\Delta$ is a space coming from Assumption~$2$ in the definition of a curve-graph-like $G$-graph (Definition~\ref{de:curve-graph-like}). 

\begin{lemma}\label{lemma:no-triple}
Let $G$ be a countable group, and let $\Theta$ be a curve-graph-like $G$-graph. Let $\calg$ be a measured groupoid over a finite measure space $Y$, and let $\rho:\calg\to G$ be a strict cocycle with trivial kernel. Let $\calh$ be a measured subgroupoid of $\calg$ of infinite type, and assume that $(\calh,\rho)$ is nowhere reducible.

Then for every Borel subset $U\subseteq Y$ of positive measure, every $(\calh_{|U},\rho)$-equivariant Borel map $\mu:U\to\Prob(\Delta)$, and a.e.\ $y\in U$, the support of the probability measure $\mu(y)$ has cardinality at most $2$. 
\end{lemma}

\begin{proof}
Let $U\subseteq Y$ be a Borel subset of positive measure, and let $\mu:U\to\Prob(\Delta)$ be a Borel map. Assume towards a contradiction that the conclusion of the lemma fails for $\mu$. Using Fubini's theorem, we see that there exists a Borel subset $V\subseteq U$ of positive measure such that for a.e.\ $y\in V$, one has $\mu(y)\otimes\mu(y)\otimes\mu(y)(\Delta^{(3)})>0$ -- where we recall that $\Delta^{(3)}\subset \Delta\times\Delta\times \Delta$ consists of all triples with mutually distinct coordinates. By restricting $\mu(y)\otimes\mu(y)\otimes\mu(y)$ to $\Delta^3$ and renormalizing, we thus get an $(\calh_{|U},\rho)$-equivariant Borel map $U\to\Prob(\Delta^{(3)})$. Using the $G$-equivariant Borel map $\Delta^{(3)}\to \calp_\ast(V(\Theta))$ given in the definition of a curve-graph-like $G$-graph (Assumption~2.(c)), we deduce an $(\calh_{|V},\rho)$-equivariant Borel map $V\to\Prob(\calp_\ast(V(\Theta)))$, contradicting Lemma~\ref{lemma:no-invariant}.
\end{proof}

Given a set $\Delta$ and $k\in\mathbb{N}$, we denote by $\calp_{\le k}(\Delta)$ the set of all  nonempty subsets of $\Delta$ of cardinality at most $k$.

\begin{lemma}\label{lemma:amenability-finite-sets}
Let $G$ be a countable group, and let $\Delta$ be a Polish $G$-space. Assume that the $G$-action on $\Delta$ is Borel amenable. Then for every $k\in\mathbb{N}^\ast$, the $G$-action on $\calp_{\le k}(\Delta)$ is Borel amenable.
\end{lemma}

\begin{proof}
Let $\nu_n:\Delta\to\Prob(G)$ be an asymptotically equivariant sequence of Borel maps. Then the sequence $\nu_n^k:\calp_{\le k}(\Delta)\to\Prob(G)$ sending a set $X\in\calp_{\le k}(\Delta)$ to $\frac{1}{|X|}\sum_{x\in X}\nu_n(x)$ is asymptotically equivariant. 
\end{proof}

We are now in position to complete our proof of Proposition~\ref{prop:irreducible-amenable}.

\begin{proof}[Proof of Proposition~\ref{prop:irreducible-amenable}]
Let $\Theta$ be a curve-graph-like $G$-graph, and let $K=K_{\bdd}\cup K_\infty$ and $\Delta$ be as in Assumption~2 of Definition~\ref{de:curve-graph-like}. As $\calh$ is amenable and $K$ is compact and metrizable, it follows from  \cite[Proposition~4.14]{Kid-survey} that there exists an  $(\calh,\rho)$-equivariant Borel map $Y\to\Prob(K)$. We will denote by $\nu_y$ the image measure of a point $y\in Y$ under this map. 

We first claim that for a.e.\ $y\in Y$, one has $\nu_y(K_\infty)=1$. Indeed, assume otherwise. Then there exists a Borel subset $U\subseteq Y$ of positive measure such that for a.e.\ $y\in U$, the probability measure $\nu_y$ gives positive measure to $K_{\bdd}$. After renormalizing, we thus get an $(\calh_{|U},\rho)$-equivariant Borel map $U\to\Prob(K_{\bdd})$. Using the map $K_\bdd\to \calp_\ast(V(\Theta))$ from the definition of a curve-graph-like $G$-graph (Assumption~2.(a)), we deduce an $(\calh_{|U},\rho)$-equivariant Borel map $U\to\Prob(\calp_\ast(V(\Theta)))$, contradicting Lemma~\ref{lemma:no-invariant}.

Using the map $K_\infty\to\Delta$ provided by the definition of a curve-graph-like $G$-graph (Assumption~2.(b)), we thus get an $(\calh,\rho)$-equivariant Borel map $\overline{\nu}:Y\to\Prob(\Delta)$. Combining Lemma~\ref{lemma:no-triple} with \cite[Lemma~3.2]{Ada}, there exists an $(\calh,\rho)$-equivariant Borel map $\theta:Y\to\calp_{\le 2}(\Delta)$ which is maximal, i.e.\ such that for every other $(\calh,\rho)$-equivariant Borel map $\theta':Y\to\calp_{\le 2}(\Delta)$ and a.e.\ $y\in Y$, one has $\theta'(y)\subseteq\theta(y)$. Using \cite[Lemma~3.3]{Ada}, one sees that such a maximal map $\theta$ is essentially unique (notice that \cite[Lemmas~3.2 and~3.3]{Ada} are stated in the context of hyperbolic groups, but their proof remains valid in our setting). As $\calh$ is normalized by $\calh'$, the map $\theta$ is both $(\calh,\rho)$-equivariant and $(\calh',\rho)$-equivariant. As the $G$-action on $\Delta$ is Borel amenable, so is the action on $\calp_{\le 2}(\Delta)$ (Lemma~\ref{lemma:amenability-finite-sets}). Since in addition $\rho$ has trivial kernel, it follows that $\calh'$ is amenable, see e.g.\ \cite[Proposition~4.33]{Kid-memoir}. This concludes our proof. 
\end{proof}

\subsection{Characterizing stabilizers of vertices of $\Theta$}

Let $\calg$ be a measured groupoid over a finite measure space $Y$. We say that $\calg$ is \emph{everywhere nonamenable} if for every Borel subset $U\subseteq Y$ of positive measure, the groupoid $\calg_{|U}$ is nonamenable. 

\begin{prop}\label{prop:stabilizer-property}
Let $G$ be a countable group, and let $\Theta$ be a curve-graph-like $G$-graph. Let $\calg$ be a discrete measured groupoid over a finite measure space $Y$, and let $\rho:\calg\to G$ be a strict action-type cocycle. 

Then for every $v\in V(\Theta)$, the groupoid $\rho^{-1}(\Stab_G(v))$ is everywhere nonamenable and normalizes an amenable subgroupoid of $\calg$ of infinite type. 
\end{prop}

\begin{proof} 
Let $v\in V(\Theta)$, and let $\calh=\rho^{-1}(\Stab_G(v))$. 
By Assumption~3.(a) from the definition of a curve-graph-like $G$-graph, the group $\Stab_G(v)$ contains a nonabelian free subgroup. As $\rho$ is action-type, it thus follows from \cite[Lemma~3.20]{Kid} that $\calh$ is everywhere nonamenable. 

By Assumption~3.(a) from the definition of a curve-graph-like $G$-graph, the group $\Stab_G(v)$ contains an infinite normal amenable subgroup $A$. Since $\rho$ is action-type, it follows that $\rho^{-1}(A)$ is a subgroupoid of $\calg$ of infinite type, it is amenable and normalized by $\calh$. This concludes our proof.
\end{proof}

\begin{prop}\label{prop:stabilizer-characterization}
Let $G$ be a countable group, and let $\Theta$ be a curve-graph-like $G$-graph. Let $\calg$ be a discrete measured groupoid over a finite measure space $Y$, and let $\rho:\calg\to G$ be a strict action-type cocycle. Let $\calh$ be a measured subgroupoid of $\calg$. Assume that $\calh$ is everywhere nonamenable and stably normalizes an amenable subgroupoid of $\calg$ of infinite type.

Then $\calh$ is stably reducible. 
\end{prop}

\begin{proof}
Up to partitioning $Y$ into at most countably many Borel subsets, we can assume that there exists an amenable subgroupoid $\cala$ which is normalized (and not just stably normalized) by $\calh$. By Lemma~\ref{lemma:reducible-vs-irreducible}, there exists a Borel partition $Y=Y_0\dunion Y_1$ such that $(\cala_{|Y_0},\rho)$ is nowhere reducible, while $(\cala_{|Y_1},\rho)$ is stably reducible. 

We observe that $Y_0$ is a null set: otherwise, as $\cala_{|Y_0}$ is of infinite type and normalized by $\calh_{|Y_0}$, Proposition~\ref{prop:irreducible-amenable} would imply that $\calh_{|Y_0}$ is amenable, a contradiction. This shows that $(\cala,\rho)$ is stably reducible, and Proposition~\ref{prop:crs} then implies $(\calh,\rho)$ is also stably reducible.
\end{proof}

\begin{lemma}\label{lemma:uniqueness}
Let $G$ be a countable group, and let $\Theta$ be a curve-graph-like $G$-graph. Let $\calg$ be a discrete measured groupoid over a finite measure space $Y$, and let $\rho:\calg\to G$ be a strict action-type cocycle. Let $v,v'\in V(\Theta)$.

If there exists a Borel subset $U\subseteq Y$ of positive measure such that $\rho^{-1}(\Stab_G(v))_{|U}$ is contained in $\rho^{-1}(\Stab_{G}(v'))_{|U}$, then $v=v'$.
\end{lemma}

\begin{proof}
Write $\calh=\rho^{-1}(\Stab_G(v))$ and $\calh'=\rho^{-1}(\Stab_G(v'))$, and assume towards a contradiction that $v\neq v'$. Using Assumption~4.(b) from the definition of a curve-graph-like $G$-graph, we can find an infinite-order element $g\in G$ that fixes $v$ but does not fix $v'$ --  and in fact, using Assumption~1, no nontrivial power of $g$ fixes $v'$. Then $\rho^{-1}(\langle g\rangle)$ is contained in $\calh$. But as $\rho$ is action-type, for every Borel subset $U\subseteq Y$ of positive measure, the subgroupoid $\rho^{-1}(\langle g\rangle)_{|U}$ is not contained in $\calh'_{|U}$, so $\calh_{|U}$ is not contained in $\calh'_{|U}$. This contradiction concludes our proof. 
\end{proof}

By combining all the above statements, we reach the following proposition.

\begin{prop}\label{prop:characterization}
Let $G$ be a countable group, and let $\Theta$ be a curve-graph-like $G$-graph. Let $\calg$ be a discrete measured groupoid over a finite measure space $Y$, let $\rho:\calg\to G$ be a strict action-type cocycle, and let $\calh$ be a measured subgroupoid of $\calg$. The following assertions are equivalent.
\begin{enumerate}
\item[(i)] There exists a partition $Y^*=\dunion_{i\in I}Y_i$ of a conull Borel subset $Y^*\subseteq Y$ into at most countably many Borel subsets such that for every $i\in I$, there exists a vertex $v_i\in V(\Theta)$ such that $\calh_{|Y_i}=\rho^{-1}(\Stab_G(v_i))_{|Y_i}$.
\item[(ii)] The following two properties hold:
\begin{enumerate}
\item[$(P_1)$] The groupoid $\calh$ is everywhere nonamenable and stably normalizes an amenable normal subgroupoid of infinite type of $\calg$.
\item[$(P_2)$] The groupoid $\calh$ is stably maximal with respect to Property~$(P_1)$, i.e.\ if $\calh'$ is another subgroupoid of $\calg$ that satisfies $(P_1)$, and if $\calh$ is stably contained in $\calh'$, then $\calh$ is stably equivalent to $\calh'$.
\end{enumerate}
\end{enumerate}
\end{prop}

\begin{proof}
We first prove that $(i)\Rightarrow (ii)$. If $(i)$ holds, then Proposition~\ref{prop:stabilizer-property} shows that $\calh$ satisfies Property~$(P_1)$. To show Property~$(P_2)$, let $\calh'$ be another subgroupoid of $\calg$ that satisfies $(P_1)$ such that $\calh$ is stably contained in $\calh'$. Then Proposition~\ref{prop:stabilizer-characterization} implies that, after refining the partition coming from $(i)$, we can assume that for every $i\in I$, there is a vertex $v'_i\in V(\Theta)$ which is $(\calh'_{|Y_i},\rho)$-invariant. In particular, there exists a conull Borel subset $Y_i^*\subseteq Y_i$ such that $\rho^{-1}(\Stab_G(v_i))_{|Y_i^*}$ is contained in $\rho^{-1}(\Stab_G(v'_i))_{|Y_i^*}$, so Lemma~\ref{lemma:uniqueness} implies that $v_i=v'_i$ whenever $Y_i$ has positive measure. Therefore $\calh'$ is stably contained in $\calh$, showing that $\calh$ satifies Property~$(P_2)$.

We now prove that $(ii)\Rightarrow (i)$. Assume that $\calh$ satisfies Properties~$(P_1)$ and~$(P_2)$. Since $\calh$ satisfies Property~$(P_1)$, Proposition~\ref{prop:stabilizer-characterization} ensures that there exists a partition $Y=\dunion_{i\in I}Y_i$ into at most countably many Borel subsets such that for every $i\in I$, there exists a vertex $v_i\in V(\Theta)$ which is $(\calh_{|Y_i},\rho)$-invariant. Let $\calh'$ be a measured subgroupoid of $\calg$ such that for every $i\in I$, one has $\calh'_{|Y_i}=\rho^{-1}(\Stab_G(v_i))$. Then $\calh$ is stably contained in $\calh'$, and Proposition~\ref{prop:stabilizer-property} implies that $\calh'$ also satisfies Property~$(P_1)$. Therefore, Property~$(P_2)$ ensures that $\calh$ is stably equivalent to $\calh'$, so $(i)$ holds. 
\end{proof}

\subsection{Characterizing adjacency}

\begin{prop}\label{prop:adjacency}
Let $G$ be a countable group, and let $\Theta$ be a curve-graph-like $G$-graph. Let $\calg$ be a discrete measured groupoid over a finite measure space $Y$, and let $\rho:\calg\to G$ be a strict action-type cocycle.  Let $v,v'\in V(\Theta)$, and let $\calh=\rho^{-1}(\Stab_G(v))$ and $\calh'=\rho^{-1}(\Stab_G(v'))$. The following assertions are equivalent.
\begin{enumerate}
\item[(i)] The vertices $v$ and $v'$ are either equal or adjacent.
\item[(ii)] Let $\calh''$ be a measured subgroupoid of $\calg$ which satisfies Properties~$(P_1)$ and~$(P_2)$ from Proposition~\ref{prop:characterization}. If $\calh\cap\calh'$ is stably contained in $\calh''$, then there exists a Borel partition $Y=Y_1\dunion Y_2$ such that $\calh''_{|Y_1}$ is stably contained in $\calh_{|Y_1}$ and $\calh''_{|Y_2}$ is stably contained in $\calh'_{|Y_2}$.
\end{enumerate} 
\end{prop}

\begin{proof}
We first prove that $(i)\Rightarrow (ii)$. 
Let $\calh''$ be a measured subgroupoid of $\calg$ that satisfies Properties~$(P_1)$ and~$(P_2)$ from Proposition~\ref{prop:characterization}, and assume that $\calh\cap\calh'$ is stably contained in $\calh''$. Proposition~\ref{prop:characterization} (applied to $\calh''$) ensures that there exist a conull Borel subset $Y^*\subseteq Y$ and a partition $Y^*=\dunion_{i\in I}U_i$ into at most countably many Borel subsets such that for every $i\in I$, one has $\calh''_{|Y_i}=\rho^{-1}(\Stab_G(v''_i))_{|Y_i}$ for some $v''_i\in V(\Theta)$. Up to passing to a further conull Borel subset and refining this partition, we can assume that for every $i\in I$, one has $(\calh\cap\calh')_{|U_i}\subseteq\calh''_{|U_i}$. Notice that $(\calh\cap\calh')_{|U_i}$ is equal to $\rho^{-1}(\Stab_G(v)\cap\Stab_G(v'))_{|U_i}$, and it stabilizes $v''_i$. As $\rho$ is action-type, every element in $\Stab_G(v)\cap\Stab_G(v')$ has a power that belongs to $\rho((\calh\cap\calh')_{|U_i})$.  But as $v$ and $v'$ are equal or adjacent in $\Theta$, the only vertices of $\Theta$ which are stabilized by a power of every element of $\Stab_G(v)\cap\Stab_G(v')$ are $v$ and $v'$: indeed, this follows from Assumption~4.(c) from the definition of a curve-graph-like $G$-graph (Definition~\ref{de:curve-graph-like}) if $v\neq v'$, and from the combination of Assumptions~1 and~4.(b) if $v=v'$. Therefore $v''_i\in\{v,v'\}$. Assertion~$(ii)$ follows by letting $Y_1$ be the union of all $U_i$ for which $v''_i=v$ (together with $Y\setminus Y^*$) and letting $Y_2$ be the union of all other $U_i$.

We now prove that $\neg (i)\Rightarrow \neg (ii)$. So assume that $v$ and $v'$ are neither equal nor adjacent in $\Theta$.
 Assumption~4.(c) from the definition of a curve-graph-like $G$-graph ensures that we can then find a vertex $v''\in V(\Theta)$ distinct from both $v$ and $v'$ such that every element in $\Stab_G(v)\cap\Stab_G(v')$ has a power contained in $\Stab_G(v'')$; combined with Assumption~1 from this definition, this shows that $\Stab_G(v)\cap\Stab_G(v')\subseteq\Stab_G(v'')$. Let $\calh''=\rho^{-1}(\Stab_G(v''))$. Then $\calh\cap\calh'$ is (stably) contained in $\calh''$. Moreover, $\calh''$ satisfies Properties~$(P_1)$ and~$(P_2)$ by Proposition~\ref{prop:characterization}. But Lemma~\ref{lemma:uniqueness} ensures that for every Borel subset $U\subseteq Y$ of positive measure, the subgroupoid $\calh''_{|U}$ is neither contained in $\calh_{|U}$ nor in  $\calh'_{|U}$. Therefore $(ii)$ fails. 
\end{proof}

\subsection{End of the proof}

The goal of the present section is to complete our proof of Theorem~\ref{theo:blueprint}.

\begin{lemma}\label{lemma:rho-rho'}
Let $G_1$ and $G_2$ be two countable groups. For every $i\in\{1,2\}$, let $\Theta_i$ be a curve-graph-like $G_i$-graph. Let $\calg$ be a discrete measured groupoid over a finite measure space $Y$, and let $\rho_1:\calg\to G_1$ and $\rho_2:\calg\to G_2$ be two strict action-type cocycles.

 For every Borel map $\phi_1:Y\to V(\Theta_1)$, there exists an essentially unique Borel map $\phi_2:Y\to V(\Theta_2)$ such that the $(\calg,\rho_1)$-stabilizer of $\phi_1$ is stably equivalent to the $(\calg,\rho_2)$-stabilizer of $\phi_2$.  
\end{lemma}

In this situation, we say that $\phi_2$ is the $(\calg,\rho_1,\rho_2)$-image of $\phi_1$.

\begin{proof}
The essential uniqueness of $\phi_2$ follows from Lemma~\ref{lemma:uniqueness}. Indeed, if $\phi_2$ and $\phi'_2$ are two such maps which are not essentially equal, then we can find a Borel subset $U\subseteq Y$ of positive measure such that $\phi_2$ and $\phi'_2$ are constant on $U$, with distinct respective values $v$ and $v'$. As the $(\calg,\rho_2)$-stabilizers of $\phi_2$ and $\phi'_2$ are stably equivalent, we can then find a Borel subset $V\subseteq U$ of positive measure such that $\rho_2^{-1}(\Stab_{G_2}(v))_{|V}=\rho_2^{-1}(\Stab_{G_2}(v'))_{|V}$, contradicting Lemma~\ref{lemma:uniqueness}.

We now prove the existence of $\phi_2$. Let $\calh$ be the $(\calg,\rho_1)$-stabilizer of $\phi_1$, a subgroupoid of $\calg$. By the implication $(i)\Rightarrow (ii)$ of Proposition~\ref{prop:characterization} (applied to the cocycle $\rho_1$), the groupoid $\calh$ satisfies Properties~$(P_1)$ and~$(P_2)$. By the implication $(ii)\Rightarrow (i)$ of Proposition~\ref{prop:characterization} (applied to the cocycle $\rho_2$), there exists a Borel map  $\phi_2:Y\to V(\Theta_2)$ such that $\calh$ is stably equivalent to the $(\calg,\rho_2)$-stabilizer of $\phi_2$.
\end{proof}

As a consequence of our characterization of adjacency given in Proposition~\ref{prop:adjacency}, we also get the following.

\begin{lemma}\label{lemma:groupoid-isomorphism-adjacence}
Let $G_1$ and $G_2$ be two countable groups. For every $i\in\{1,2\}$, let $\Theta_i$ be a curve-graph-like $G_i$-graph. Let $\calg$ be a discrete measured groupoid over a finite measure space $Y$, and let $\rho_1:\calg\to G_1$ and $\rho_2:\calg\to G_2$ be two strict action-type cocycles.

Let $v_1,v'_1\in V(\Theta_1)$ and $v_2,v'_2\in V(\Theta_2)$. Assume that $\rho_1^{-1}(\Stab_{G_1}(v_1))=\rho_2^{-1}(\Stab_{G_2}(v_2))$ and that $\rho_1^{-1}(\Stab_{G_1}(v'_1))=\rho_2^{-1}(\Stab_{G_2}(v'_2))$. Then $v_1$ and $v'_1$ are equal or adjacent if and only if $v_2$ and $v'_2$ are equal or adjacent.
\end{lemma}

\begin{proof}
By symmetry, it is enough to assume that $v_1$ and $v'_1$ are either equal or adjacent, and prove that so are $v_2$ and $v'_2$. Let $\calh=\rho_1^{-1}(\Stab_{G_1}(v_1))=\rho_2^{-1}(\Stab_{G_2}(v_2))$ and $\calh'=\rho_1^{-1}(\Stab_{G_1}(v'_1))=\rho_2^{-1}(\Stab_{G_2}(v'_2))$. Let $\calh''$ be a subgroupoid of $\calg$ that satisfies Properties~$(P_1)$ and~$(P_2)$ from Proposition~\ref{prop:characterization}. Assume that $\calh\cap\calh'$ is stably contained in $\calh''$. Then by the implication $(i)\Rightarrow (ii)$ from Proposition~\ref{prop:adjacency} applied to $\rho_1$, there exists a Borel partition $Y=Y_1\dunion Y_2$ such that $\calh''_{|Y_1}$ is stably contained in $\calh_{|Y_1}$ and $\calh''_{|Y_2}$ is stably contained in $\calh'_{|Y_2}$. By the implication $(ii)\Rightarrow (i)$ from Proposition~\ref{prop:adjacency} applied to $\rho_2$, it follows that $v_2$ and $v'_2$ are either equal or adjacent, as required.
\end{proof}

We are now in position to complete our proof of Theorem~\ref{theo:blueprint}.

\begin{proof}[Proof of Theorem~\ref{theo:blueprint}] 
We first prove that if $G_1$ and $G_2$ are measure equivalent, then the graphs $\Theta_1$ and $\Theta_2$ are isomorphic. By the first assertion from Proposition~\ref{prop:me-groupoid}, we can find a measured groupoid $\calg$ over a finite measure space $Y$ of positive measure, and two strict action-type cocycles $\rho_1:\calg\to G_1$ and $\rho_2:\calg\to G_2$. 
 
Given a vertex $v\in V(\Theta_1)$, we let $\phi_{1,v}:Y\to V(\Theta_1)$ be the constant map with value $v$. We let $\phi_{2,v}:Y\to V(\Theta_2)$ be the $(\calg,\rho_1,\rho_2)$-image of $\phi_{1,v}$. We define a Borel map $\theta:Y\times V(\Theta_1)\to V(\Theta_2)$ by letting $\theta(y,v):=\phi_{2,v}(y)$ for every $y\in Y$ and every $v\in V(\Theta_1)$. The fact that $\Theta_1$ and $\Theta_2$ are isomorphic will follow from the following claim.

\begin{claim2} 
For a.e.\ $y\in Y$, the map $\theta(y,\cdot)$ is a graph isomorphism from $\Theta_1$ to $\Theta_2$. 
\end{claim2}

\begin{proof}[Proof of Claim]
First, for a.e.\ $y\in Y$, the map $v\mapsto\theta(y,v)$ is injective. Indeed, otherwise, there exists a Borel subset $U\subseteq Y$ of positive measure and two distinct vertices $v,v'\in V(\Theta_1)$ such that for a.e.\ $y\in U$, one has $\theta(y,v)=\theta(y,v')$. This means that the $(\calg_{|U},\rho_1)$-stabilizer $\calh$ of the constant map with value $v$ and the $(\calg_{|U},\rho_1)$-stabilizer $\calh'$ of the constant map with value $v'$ are both stably equivalent to the $(\calg_{|U},\rho_2)$-stabilizer of $\theta(\cdot,v)$ -- in particular $\calh$ and $\calh'$ are stably equivalent. This is a contradiction to Lemma~\ref{lemma:uniqueness}.

By Lemma~\ref{lemma:groupoid-isomorphism-adjacence}, for any two distinct vertices $v,v'\in V(\Theta_1)$, and a.e.\ $y\in Y$, the vertices $\theta(y,v)$ and $\theta(y,v')$ are adjacent if and only if $v$ and $v'$ are. As $V(\Theta_1)$ is countable, it follows that for a.e.\ $y\in Y$, the map $\theta(y,\cdot)$ preserves both adjacency and non-adjacency. 

Lastly, for a.e.\ $y\in Y$, the map $\theta(y,\cdot):V(\Theta_1)\to V(\Theta_2)$ is surjective. Indeed, let $w\in V(\Theta_2)$, and let $\psi_{2,w}:Y\to V(\Theta_2)$ be the constant map with value $w$. Let $\psi_{1,w}:Y\to V(\Theta_1)$ be the $(\calg,\rho_2,\rho_1)$-image of $\psi_{2,w}$. Consider a Borel partition $Y=\dunion_{n\in\mathbb{N}}Y_n$ such that for every $n\in\mathbb{N}$, the map $(\psi_{1,w})_{|Y_n}$ is constant, with value a vertex $w_n$. Then $(\psi_{2,w})_{|Y_n}$ is the $(\calg_{|Y_n},\rho_1,\rho_2)$-image of $(\phi_{1,w_n})_{|Y_n}$. This implies that the restriction to $Y_n$ of the $(\calg,\rho_1,\rho_2)$-image of $\phi_{1,w_n}$ is constant with value $w$. Therefore, for a.e.\ $y\in Y_n$, we have $\theta(y,w_n)=w$. In particular, for a.e.\ $y\in Y_n$, the vertex $w$ belongs to the image of $\theta(y,\cdot)$. As this is true for every $n\in\mathbb{N}$, and as $V(\Theta_2)$ is countable, our claim follows.
\end{proof}

We have just built a map $\tilde{\theta}:Y\to\Isom(\Theta_1\to\Theta_2)$ and this map is Borel by our choice of topology on $\Isom(\Theta_1\to\Theta_2)$. In view of Proposition~\ref{prop:me-groupoid}, in order to prove that $\Isom(\Theta_1\to\Theta_2)$ is a universal $(G_1,G_2)$-factor, we are left showing that $\tilde{\theta}$ is almost $(\rho_1,\rho_2)$-equivariant, i.e.\ there exists a conull Borel subset $Y^*\subseteq Y$ such that for all $g\in\calg_{|Y^*}$, one has $\tilde{\theta}(r(g))=(\rho_1(g),\rho_2(g))\tilde{\theta}(s(g))$. This amounts to proving that for all $g\in\calg_{|Y^*}$ and every $v\in V(\Theta_1)$, one has $$\theta(r(g),\rho_1(g)v)=\rho_2(g)\theta(s(g),v).$$ Given $g_1\in G_1$, $g_2\in  G_2$, $v_1\in V(\Theta_1)$ and $v_2\in V(\Theta_2)$, we let $$\calg_{g_1,g_2,v_1,v_2}:=\{g\in\calg|\rho_1(g)=g_1,\rho_2(g)=g_2,\phi_{2,v_1}(s(g))=v_2\}.$$ Notice that for a fixed $v_1$, the subsets $\calg_{g_1,g_2,v_1,v_2}$ form a Borel partition of $\calg$. Therefore, it is enough to fix $g_1\in G_1$, $g_2\in G_2$, $v_1\in V(\Theta_1)$ and $v_2\in V(\Theta_2)$ and prove that there exists a conull Borel subset $Y^*\subseteq Y$ such that for all $g\in \calg_{g_1,g_2,v_1,v_2}\cap\calg_{|Y^*}$, one has $\theta(r(g),\rho_1(g)v_1)=\rho_2(g)\theta(s(g),v_1)$, in other words $\theta(r(g),g_1 v_1)=g_2v_2$. Using the theorem of  Lusin and Novikov mentioned in Remark~\ref{rk-normal}, it is in fact enough to prove this fact for a Borel subset of $\calg_{g_1,g_2,v_1,v_2}$ that determines a partial isomorphism between two Borel subsets $A,B\subseteq Y$ -- indeed $\calg_{g_1,g_2,v_1,v_2}$ can be written as a countable union of such subsets.

It is enough to show that for a.e.\ $y\in B$, one has $\theta(y,g_1v_1)=g_2v_2$, in other words we need to show 
\begin{center}
($\ast$) the $(\calg_{|B},\rho_1)$-stabilizer of $\phi_{1,g_1v_1}$ is stably equivalent to the $(\calg_{|B},\rho_2)$-stabilizer of $\psi_{2,g_2v_2}$ (recall that $\psi_{2,w}$ denotes the constant map $Y\to V(\Theta_2)$ value $w$).
\end{center}
  Notice that the $(\calg_{|B},\rho_1)$-stabilizer of the constant map $\phi_{1,g_1v_1}$ is equal to the subgroupoid made of all $gh(g')^{-1}$ with $g,g'\in\calg_{g_1,g_2,v_1,v_2}$ and $h$ in the $(\calg_{|A},\rho_1)$-stabilizer of $\phi_{1,v_1}$. Likewise, the  $(\calg_{|B},\rho_2)$-stabilizer of  $\psi_{2,g_2v_2}$ is equal to the subgroupoid made of all $gh'(g')^{-1}$ with $g,g'\in\calg_{g_1,g_2,v_1,v_2}$ and $h'$ in the $(\calg_{|A},\rho_2)$-stabilizer of  $\psi_{2,v_2}$. By definition of $\calg_{g_1,g_2,v_1,v_2}$, the $(\calg_{|A},\rho_1)$-stabilizer of $\phi_{1,v_1}$ is stably equivalent to the $(\calg_{|A},\rho_2)$-stabilizer of $\psi_{2,v_2}$. Now ($\ast$) follows and the proof is complete.
\end{proof}

\section{Curve graphs for Artin groups}\label{sec:theta}

\emph{Building upon earlier work of Crisp \cite{crisp}, we now introduce an analogue of the curve graph for a $2$-dimensional Artin group $G_\Gamma$ of hyperbolic type. This graph can be understood in two different ways: geometrically, the \emph{fixed set graph} encodes the intersection pattern of standard trees in the modified Deligne complex $\del$ of $G_\Gamma$; algebraically, the \emph{graph of standard abelian subgroups} encodes the commutation pattern of various natural abelian subgroups of $G_\Gamma$. These graphs are intimately related: the graph of standard abelian subgroups turns out to be $G_\Gamma$-equivariantly isomorphic to the core subgraph of the fixed set graph. We will thoroughly study vertex and edge stabilizers of these graphs: this will enable us to check that the fixed set graph is indeed a curve-graph-like $G$-graph as defined in the previous section. Some arguments in this section are inspired from the first four sections of Crisp's work \cite{crisp}, though we intend the make the exposition here self-contained.}

\subsection{Standard trees and stabilizers in the Deligne complex}

As usual, we let $G=G_\Gamma$ be a $2$-dimensional Artin group with underlying graph $\Gamma$, and we denote by $\del$ its modified Deligne complex, equipped with the Moussong metric which is $\CAT$ (see Section~\ref{sec:background}). For every $g\in G$, we denote by $F_g$ the fixed point set of $g$ with respect to the action of $G$ on $\del$. When $g$ is conjugate to a standard generator, then $F_g$ is a 1-dimensional convex subcomplex of $\mathbb D_\Gamma$ whose vertices are alternately rank 1 and rank 2 vertices of $\mathbb D_\Gamma$ -- hence $F_g$ is a tree which comes with a natural bipartite structure. We make the following definition.

\begin{de}[Standard tree]
A \emph{standard tree} is a subset of $\del$ equal to $F_g$ for some element $g\in G$ that is conjugate to a standard generator. 
\end{de}

When $g$ is conjugate to a standard generator which does not correspond to an isolated vertex of $\Gamma$, then $F_g$ has more than one point. The following lemma due to Crisp will be used extensively in the present section. 

\begin{lemma}[{Crisp \cite[Lemma~8]{crisp}}]\label{lemma:fix-sets}
Let $G_{\Gamma}$ be a $2$-dimensional Artin group of hyperbolic type with defining graph $\Gamma$, and let $g\in G_{\Gamma}\setminus\{1\}$.
\begin{enumerate}
\item If $g=hs^mh^{-1}$ is conjugate to a nontrivial power of a standard generator $s$, then $F_g=F_{hsh^{-1}}$ is a standard tree.
\item If $g$ is conjugate into an edge subgroup but not conjugate to a power of any standard generator, then $F_g$ is reduced to a single rank $2$ vertex.
\item If $g$ is not conjugate into an edge subgroup, then $F_g=\emptyset$.
\end{enumerate}
\end{lemma}

The following lemma will also be useful.

\begin{lemma}\label{lemma:intersection-standard-subtrees}
Let $G_{\Gamma}$ be a $2$-dimensional Artin group with defining graph $\Gamma$, and let $T_1$ and $T_2$ be two standard trees in $\mathbb{D}_\Gamma$. If $T_1\cap T_2$ contains a rank one vertex, then $T_1=T_2$.
\end{lemma}

\begin{proof}
Let $x\in T_1\cap T_2$ be a rank one vertex. Then there exist an element $g\in G_{\Gamma}$ and a unique standard generator $s$ of $G_{\Gamma}$ such that $x=g\langle s\rangle$. Let $h:=gsg^{-1}$, so that $\stab(x)=\langle h\rangle$. For every $i\in\{1,2\}$, let $g_i$ be a conjugate of a standard generator such that $T_i=F_{g_i}$.  Then $g_i\in \stab(x)$, so $g_i=h^{m_i}$ for some $m_i\neq 0$. Using the first assertion from Lemma~\ref{lemma:fix-sets}, we get $T_1=F_{h^{m_1}}=F_{h}=F_{h^{m_2}}=T_2$.
\end{proof}

\subsection{The fixed set graph}

A rank 2 vertex $gG_{st}$ of the modified Deligne complex $\del$ is \emph{irreducible} if $m_{st}\ge 3$ and is \emph{reducible} if $m_{st}=2$. Following Crisp \cite[Section~4]{crisp}, we make the following definition.

\begin{de}[Fixed set graph]
	\label{de:theta}
Let $G_{\Gamma}$ be a 2-dimensional Artin group of hyperbolic type with defining graph $\Gamma$. 
\begin{enumerate}
\item The \emph{modified fixed set graph} $\Theta'_{\Gamma}$ is the graph having one vertex (of \emph{type I}) for every standard tree of $\del$, one vertex (of \emph{type II}) for every irreducible rank 2 vertex of $\del$, where two different vertices $v_1$ and $v_2$ are joined by an edge if either
\begin{enumerate}
	\item $v_1$ is of type I, $v_2$ is of type II and the standard tree associated with $v_1$ contains the rank 2 vertex associated with $v_2$, or
	\item both $v_1$ and $v_2$ are of type I, and  there exist commuting elements $g_1$ and $g_2$ such that for every $i\in\{1,2\}$, the standard tree corresponding to $v_i$ is $F_{g_i}$.
\end{enumerate}
\item The \emph{fixed set graph} $\Theta_{\Gamma}$ is 
the subgraph obtained from $\Theta'_\Gamma$ by removing all vertices of valence $0$ or $1$ and all edges containing a valence one vertex.
\end{enumerate}
\end{de}

\begin{rk}\label{rk:theta-theta'}
Vertices of type II of $\Theta'_\Gamma$ always belong to $\Theta_\Gamma$. On the other hand, a vertex of $\Theta'_\Gamma$ of type I belongs to $\Theta_\Gamma$ if and only if it corresponds to a standard tree $F_g$ where $g$ is conjugate to a standard generator associated to a vertex of $\Gamma$ of valence at least $2$. Indeed, the if direction is clear, and the only if direction can easily be deduced from \cite[Remark~10 and Lemma~11]{crisp} -- details are left to the reader.
\end{rk}

Notice that when $G_{\Gamma}$ is of large type (i.e.\ all edge labels in $\Gamma$ are at least $3$), the definition of the fixed set graph coincides with the one given by Crisp \cite[Section~4]{crisp}: in this case $\Theta'_{\Gamma}$ is a bipartite graph having one vertex (of type I) for every standard tree $T$, one vertex (of type II) for every vertex $x$ of $\mathbb{D}_\Gamma$ that belongs to more than one standard tree, with an edge joining $x$ to $T$ whenever $x\in T$. 
In general, we prefer the above definition, which will enable us to relate $\Theta'_{\Gamma}$ to the graph of standard abelian subgroups described later in this section (Definition~\ref{de:curve complex}) and prove that it is curve-graph-like.

We also mention that natural analogues of the curve graph exist beyond the class of Artin groups of hyperbolic type, see e.g.\ \cite[Definitions~10.3 and~10.7]{HO}; however as such, Definition~\ref{de:theta} is not so well-adapted past this setting. 

\begin{lemma}\label{lemma:equivariant-map}
Let $G=G_\Gamma$ be a $2$-dimensional Artin group of hyperbolic type with connected defining graph $\Gamma$. If $\Gamma$ is not an single edge labeled by $2$, then there is a $G$-equivariant Borel map $f:\del\to V(\Theta_\Gamma)$. 
\end{lemma}

\begin{proof}
A vertex $v$ of $\del$ is \emph{special} if either $v$ is a rank 2 vertex represented by the coset $gG_e$ such that $e$ has label at least $3$, or $v$ is a rank 1 vertex represented by $g\langle a\rangle$ for $a\in V\Gamma$ such that $a$ has valence at least $2$ in $\Gamma$. Each special vertex $v$ corresponds to an element $\bar v$ of $V(\Theta_\Gamma)$. The assumption on $\Gamma$ implies that each rank 2 vertex $v$ of $\del$ is either a special vertex itself, or else it corresponds to an edge of $\Gamma$ labeled by $2$ and is adjacent to a rank 1 special vertex $w$. Moreover, in the latter case, we have $\stab_G(v)\subseteq \stab_{G}(\bar w)$ where $\bar w $ is the vertex of $\Theta_\Gamma$ corresponding ot $w$. Each rank 1 vertex $v$ is either a special vertex itself, or is adjacent to a rank 2 special vertex $w_1$ with $\stab_G(v)\subset\stab_G(\bar w_1)$, or is adjacent to a rank 2 non-special vertex $w$ such that $\st(w)$ contain a special rank 1 vertex $w_2$ with $\stab_G(v)\subset\stab_G(\bar w_2)$. 
	
For each vertex $v\in\del$, let $\Lambda_v$ be the collection of vertices of $\del$ which are adjacent to $v$ and have rank bigger than $v$. Define $\st^+(v)$ be the full subcomplex of $\del$ spanned by $v$ and $\Lambda_v$. Let $O_v= \st^+(v)\setminus(\cup_{w\in \Lambda_v} \st^+(w))$. Then the collection of $O_v$ with $v$ ranging over vertices of $\del$ forms a Borel partition of $\del$. Note that for any point $p\in O_v$, we have $\stab(p)=\stab(v)$. 

Now let $V$ be a finite set of representatives of the $G$-orbits of vertices in $\del$. By the previous paragraph, for every $v\in V$, we can find a special vertex $w$ in $\del$ such that $\stab_G(v)\subseteq \stab_G(\bar w)$. This allows us to define a $G$-equivariant Borel map $\del\to V(\Theta_\Gamma)$ by letting $f(x)=\bar w$ whenever $v\in V$ and $x\in O_v$, and extending equivariantly to the whole $\del$.
\end{proof}  

\subsection{Standard abelian subgroups}

In the sequel, we will often blur the distinction between vertices of the defining graph $\Gamma$ and the associated standard generators of $G_\Gamma$. Given two vertices $s,t$ of $\Gamma$ that span an edge in $\Gamma$ with $2<m_{st}<\infty$, we let $Z_{st}$ be the center of $G_{st}$: this is an infinite cyclic subgroup generated by $(st)^{m_{st}}$ if $m_{st}$ is odd, and by $(st)^{m_{st}/2}$ if $m_{st}$ is even, as $G_{st}$ is itself an Artin group whose defining graph is an edge (Theorem~\ref{lem:injective}).

\begin{de}[Standard abelian subgroups]
	\label{de:sas}
Let $G_{\Gamma}$ be a 2-dimensional Artin group of hyperbolic type with defining graph $\Gamma$. A \emph{standard abelian subgroup} of $G_{\Gamma}$ is either the trivial subgroup or one of the following:
\begin{enumerate}
	\item[(1)] a conjugate of $G_{st}$ for $s,t\in V\Gamma$ with $m_{st}=2$;
	\item[(2)] a conjugate of $\langle s,Z_{st} \rangle$ where $s,t\in V\Gamma$ satisfy $2<m_{st}<\infty$;
	\item[(3)] a conjugate of $Z_{st}$ where $s,t\in V\Gamma$ satisfy $2<m_{st}<\infty$;
	\item[(4)] a conjugate of $G_s$ for $s\in V\Gamma$.
\end{enumerate}
\end{de}

\begin{lemma}
	\label{lem:type 4}
Let $s,t\in V\Gamma$ with $m_{st}<\infty$. 
\begin{enumerate}
\item Let $H$ be a standard abelian subgroup of type (3) in $G_{\Gamma}$. If $H\subseteq G_{st}$, then $m_{st}>2$ and $H=Z_{st}$. 
	\item Let $H$ be a standard abelian subgroup of type (4) in $G_{\Gamma}$. If $H\subseteq G_{st}$, then $H$ is conjugate inside $G_{st}$ to either $G_s$ or $G_t$.	
\end{enumerate}

\end{lemma}

\begin{proof}
We first prove Assertion~1. Note that if $H=g Z_{s't'}g^{-1}$, then $H$ fixes the vertex $y'$ of $\del$ corresponding to $gG_{s't'}$. Moreover, the fix point set of $H$ is $\{y'\}$ by the second point in Lemma~\ref{lemma:fix-sets}. On the other hand, the fix point set of $G_{st}$ is a singleton $\{y\}$. Since $H\subseteq G_{st}$, we have $y=y'$. Therefore $G_{st}=G_{s't'}$ and $g\in G_{s't'}$. Assertion~1 follows.

We now prove Assertion~2. Let $H$ be a standard abelian subgroup of type $(4)$.  
Then $H$ is equal to the $G_\Gamma$-stabilizer of some rank 1 vertex $x\in\mathbb D_\Gamma$. The fixed point set of $H$ is a standard tree $T$. 
Let $y\in\mathbb D_\Gamma$ be the vertex corresponding to $G_{st}$. Then $y\in T$ -- in particular $T$ contains at least two distinct vertices. Let $x'\in T$ be a rank 1 vertex adjacent to $y$ (this always exists by the bipartite structure on $T$). Then there exists $g\in G_{st}$ such that $x'=gG_s$ or $x'=gG_t$, so $\Stab(x')$ is conjugate in $G_{st}$ to either $G_s$ or $G_t$. Let now $T'$ be the standard tree which is the fixed point set of $\stab(x')$. By Lemma~\ref{lemma:intersection-standard-subtrees}, we have $T=T'$. So $\Stab(x)$ fixes $x'$ and $\Stab(x')$ fixes $x$, showing that $H=\stab(x)=\stab(x')$. So $H$ is conjugate inside $G_{st}$ to either $G_s$ or $G_t$.
\end{proof}

\begin{lemma}\label{lemma:commuting-sas}
Two distinct rank $1$ standard abelian subgroups $H_1$ and $H_2$ commute if and only if they generate a rank $2$ standard abelian subgroup. In this case, either 
\begin{enumerate}
\item there exist $s,t\in V\Gamma$ with $m_{st}=2$ and $g\in G_\Gamma$ such that $H_1=g\langle s\rangle g^{-1}$ and $H_2=g\langle t\rangle g^{-1}$ (and in this case $\langle H_1,H_2\rangle=gG_{st}g^{-1}$ is of type (1)), or 
\item there exist $s,t\in V\Gamma$ with $2<m_{st}<+\infty$ and $g\in G_\Gamma$ such that $\{H_1,H_2\}=\{gZ_{st}g^{-1},g\langle s\rangle g^{-1}\}$ (and in this case $\langle H_1,H_2\rangle=g\langle s,Z_{st}\rangle g^{-1}$ is of type (2)). 
\end{enumerate}
\end{lemma}

\begin{proof}
Let $H_1$ and $H_2$ be two distinct commuting rank $1$ standard abelian subgroups; we aim to show that $\langle H_1, H_2\rangle$ is a rank $2$ standard abelian subgroup.  Notice that both $H_1$ and $H_2$ are of type $(3)$ or $(4)$.

We first assume that one of the $H_i$ (say $H_1$) is of type $(3)$; up to conjugation $H_1=Z_{st}$ for some $s,t\in V\Gamma$ with $2<m_{st}<\infty$. By Lemma~\ref{lemma:fix-sets}, the fixed point set of $H_1$ in $\del$ is reduced to the vertex $x$ associated to $G_{st}$. Thus $H_2$ fixes $x$, so $H_2\subseteq G_{st}$. If $H_2$ is of type (3), the first assertion of Lemma~\ref{lem:type 4} implies that $H_2=H_1$, a contradiction. Thus $H_2$ of type (4). The second assertion of Lemma~\ref{lem:type 4} thus implies that $H_2$ is conjugate inside $G_{st}$ to either $G_s$ or $G_t$. Thus $\langle H_1,H_2\rangle$ is a standard abelian subgroup of type (2). 

We now assume that both $H_1$ and $H_2$ are of type (4). For every $i\in\{1,2\}$, the fix point set of $H_i$ is standard subtree $T_i$. Notice that $T_1\neq T_2$, because for every $i\in\{1,2\}$, the rank one vertex associated to $H_i$ belongs to $T_i$ but not to $T_j$ with $j\neq i$. As $H_1$ and $H_2$ commute, $H_1$ stabilizes $T_2$. Since $H_1$ acts elliptically on $\mathbb D_\Gamma$, in particular it acts elliptically on $T_2$ and therefore fixes a point in $T_2$. Hence $T_1\cap T_2\neq \emptyset$. By Lemma~\ref{lemma:intersection-standard-subtrees}, the intersection $T_1\cap T_2$ is a single vertex of rank 2. Up to conjugation, we assume this vertex corresponds to $G_{st}$. Both $H_1$ and $H_2$ fix this vertex, so $H_1$ and $H_2$ are both contained in $G_{st}$. By the first assertion of Lemma~\ref{lem:type 4}, for every $i\in\{1,2\}$, the group $H_i$ is conjugate inside $G_{st}$ to either $G_s$ or $G_t$. As $H_1$ and $H_2$ commute, it follows from Lemma~\ref{lemma:crisp-commute} below due to Crisp that $m_{st}=2$ and $\langle H_1,H_2\rangle=G_{st}$.
\end{proof}

\begin{lemma}[{Crisp \cite[Lemma~7(iii)]{crisp}}]\label{lemma:crisp-commute}
Let $s$ and $t$ be two standard generators with $2<m_{st}<\infty$. Let $u,v\in G_{st}$ be two elements, each of which is conjugate in $G_{st}$ to either $s$ or $t$. If there exist $k,l\in\mathbb{Z}\setminus\{0\}$ such that $u^k$ and $v^l$ commute, then $u=v$.
\end{lemma}

\subsection{Intersections and chains of standard abelian subgroups}

The goal of the present section is to bound the length of a decreasing chain of standard abelian subgroups (Corollary~\ref{cor:chain} below).

\begin{lemma}
	\label{lem:normalizer}
	Let $s$ and $t$ be two standard generators with $2<m_{st}<\infty$. 
	\begin{enumerate}
		\item If $h\in G_{st}$ satisfies $h^{-1}s^kh\in G_s$ for some integer $k\neq 0$, then $h\in \langle s,Z_{st}\rangle$. In particular, the normalizer $N_{G_{st}}(s)$ is equal to $\langle s,Z_{st}\rangle$.
		\item Let $u,v\in G_{st}$ be two elements, each of which is conjugate in $G_{st}$ to either $s$ or $t$. If $F_u=F_v$, then $\langle u,Z_{st}\rangle=\langle v,Z_{st}\rangle$. If $F_u\neq F_v$, then $\langle u,Z_{st}\rangle\cap \langle v,Z_{st}\rangle=\langle Z_{st}\rangle$.
	\end{enumerate}
\end{lemma}

\begin{proof}
	The first assertion is a special case of \cite[Theorem 5.2]{P}. The first part of (2) is clear as $F_u=F_v$ if and only if $\langle u\rangle=\langle v\rangle$. For the second part of (2), take $h\in \langle u,Z_{st}\rangle\cap \langle v,Z_{st}\rangle$. Then $h=u^{k_1}Z^{k_2}_{st}=v^{k_3}Z^{k_4}_{st}$. It suffices to show $k_1=0$ or $k_3=0$. Suppose $k_1\neq 0$ and $k_3\neq 0$. Then $u^{k_1}=v^{k_3}Z^{k_4-k_2}_{st}$, which implies that $u^{k_1}$ commutes with $v^{k_3}$. But Lemma~\ref{lemma:crisp-commute} then implies that $u=v$, which contradicts $F_u\neq F_v$.
\end{proof}

\begin{lemma}
	\label{lem:sas intersection}
	The intersection of two standard abelian subgroups is a standard abelian subgroup. Moreover, if a standard abelian subgroup is finite index in another standard abelian subgroup, then they are equal.
\end{lemma}

\begin{proof}
 Let $H_1$ and $H_2$ be two standard abelian subgroups. 	We first assume that either $H_1$ or $H_2$, say $H_1$, has rank $1$. The lemma then follows from the following claim: if $H_1$ has a finite index subgroup $H'_1$ contained in $H_2$, then $H_1\subseteq H_2$. We now prove the claim.
	
If $H_2$ also has rank 1, then $H_1$ and $H_2$ are commensurable, so have the same fixed point sets in $\del$ (Lemma~\ref{lemma:fix-sets}).
So either they are both of type (3), their common fixed point is unique and (up to conjugation) corresponds to some $G_{st}$, and Lemma~\ref{lem:type 4}(1) yields $H_1=H_2=Z_{st}$. Else, they are both of type (4) and fix a common rank 1 vertex of $\mathbb{D}_\Gamma$, and again $H_1=H_2$.
	
	If $H_2$ has rank 2, it fixes a unique rank 2 vertex $x\in \mathbb D_\Gamma$ (Lemma~\ref{lemma:fix-sets}). Thus $H'_1$ fixes $x$, and in fact $H_1$ fixes $x$ (as follows from Lemma~\ref{lemma:fix-sets}(1) if $H_1$ is of type (4) and from Lemma~\ref{lemma:fix-sets}(2) if $H_1$ is of type (3)). Up to conjugation $x$ corresponds to some $G_{st}$, and $H_1,H_2\subseteq G_{st}$. By Lemma~\ref{lem:type 4}, we are reduced to the special case where $G_{\Gamma}=G_{st}$. When $G_{\Gamma}=G_{st}$, our claim follows from Lemma~\ref{lemma:crisp-commute} if $H_1$ is of type (4) and $H_2$ is of type (2), and is trivial in all other cases.
	
	This concludes the proof of our claim, and finishes the proof of the lemma when either $H_1$ or $H_2$ has rank $1$. We now assume that both $H_1$ and $H_2$ are of rank 2 -- so the fix point set of $H_i$ in $\del$ is reduced to a rank two vertex $x_i$. 
	
	We first assume that $x_1\neq x_2$. We assume in addition that $H_1\cap H_2$ is nontrivial, as otherwise the conclusion is obvious. The common stabilizer of $x_1$ and $x_2$ is therefore infinite, so Lemma~\ref{lemma:fix-sets} implies that $x_1$ and $x_2$ are contained in a standard tree, and $H_1\cap H_2$ then fixes the geodesic segment $\overline{x_1x_2}$ pointwise. As $x_1\neq x_2$, the segment $\overline{x_1x_2}$ contains a rank one vertex $x$, and $H_1\cap H_2\subseteq \stab(x)$. Let $H=\stab(x)$, which is also a standard abelian subgroup. Then $H$ has a finite index subgroup contained in both $H_1$ and $H_2$. Hence $H\subseteq H_1$ and $H\subseteq H_2$ by the claim established at the beginning of this proof. Therefore $H_1\cap H_2=H$. 
	
	We finally assume that $x_1=x_2$. Then either both $H_1$ and $H_2$ are of type (1), in which case $H_1=H_2$, or they are both of type (2), of the form $\langle u,Z_{st}\rangle$ and $\langle v,Z_{st}\rangle$. Using the second assertion of Lemma~\ref{lem:type 4}, we see that $u$ and $v$ are both conjugate in $G_{st}$ to either $s$ or $t$. We can then conclude using the second assertion from Lemma~\ref{lem:normalizer}.
\end{proof}

\begin{cor}
	\label{cor:chain}
	Let $H_1\supseteq H_2\supseteq H_3\supseteq\cdots$ be a non-ascending chain of standard abelian subgroups. Then there exist at most two indices $i$ such that $H_i\supsetneq H_{i+1}$.
	\qed
\end{cor}

\begin{lemma}\label{lemma:rank-one-in-rank-two}
Let $H\subseteq G_{\Gamma}$ be a rank two standard abelian subgroup. Then $H$ contains exactly two rank one standard abelian subgroups.
\end{lemma}

\begin{proof}
If $H$ is of type (1), then up to conjugation $H=G_{st}$ with $m_{st}=2$, and Lemma~\ref{lem:type 4} implies that $\langle s\rangle$ and $\langle t\rangle$ are the only rank one abelian subgroups of $G_{st}$.

 If $H$ is of type (2), then up to conjugation $H=\langle s, Z_{st}\rangle$ with $2<m_{st}<\infty$. Let $H'\subseteq H$ be a standard abelian subgroup of rank 1. By Lemma~\ref{lem:type 4}, either $H'=\langle Z_{st}\rangle$, or $H'=\langle h\rangle$ where $h$ is conjugate (in $G_{st}$) to either $s$ or $t$. 
As $h$ commutes with $s$, Lemma~\ref{lemma:crisp-commute} yields $h=s$. Therefore the only rank one standard abelian subgroups contained in $H$ are $Z_{st}$ and $\langle s\rangle$.
\end{proof}

\subsection{The graph of standard abelian subgroups}

\begin{de}[Graph of standard abelian subgroups]
	\label{de:curve complex}
Let $G_{\Gamma}$ be a 2-dimensional Artin group of hyperbolic type with defining graph $\Gamma$. The \emph{graph of standard abelian subgroups} is the graph $\Theta_\Gamma^{\mathrm{sas}}$ whose vertices are the rank $1$ standard abelian subgroups of $G_{\Gamma}$, where two vertices are adjacent if the corresponding subgroups commute. 
\end{de}

The group $G_\Gamma$ acts on $\Theta_\Gamma^{\mathrm{sas}}$ by conjugation.

\begin{prop}
\label{prop:curve complex}	
There exists a $G_{\Gamma}$-equivariant isomorphism between the graph of standard abelian subgroups $\Theta_\Gamma^{\mathrm{sas}}$ and the modified fixed set graph $\Theta'_\Gamma$.
\end{prop}

\begin{proof}
We first define a map $\Phi:V\Theta_\Gamma^{\mathrm{sas}}\to V\Theta'_\Gamma$. Let $H\in V\Theta_\Gamma^{\mathrm{sas}}$ be a rank $1$ standard abelian subgroup. If $H$ is of type (3), then it  fixes a unique (irreducible) rank 2 vertex $x\in\mathbb D_\Gamma$ (Lemma~\ref{lemma:fix-sets}(2)), and we let $\Phi(H)=x$ (viewed as a vertex of $\Theta'_\Gamma$). If $H$ is of type (4), then the fixed point set of $H$ in $\del$ is a standard subtree $F$, and we let $\Phi(H)=F$. 

Lemma~\ref{lemma:commuting-sas} ensures that $\Phi$ sends adjacent vertices of $\Theta_\Gamma^{\mathrm{sas}}$ to adjacent vertices of $\Theta'_\Gamma$, so induces a graph map from $\Theta_\Gamma^{\mathrm{sas}}$ to $\Theta'_\Gamma$. This map is $G_\Gamma$-equivariant by construction.

We now define a map $\Phi':V\Theta'_\Gamma\to V\Theta^{\mathrm{sas}}_\Gamma$. Let $v\in V\Theta'_\Gamma$. If $v$ corresponds to a coset $gG_{st}$, we let $\Phi'(v)$ be the vertex of $\Theta^{\mathrm{sas}}_\Gamma$ associated with $gZ_{st}g^{-1}$ (this does not depend on the choice of $g$ in the left coset as $Z_{st}$ is central in $G_{st}$). If $v$ corresponds to a standard subtree $T$, we let $\Phi'(v)$ be the vertex of $\Theta^{\mathrm{sas}}_\Gamma$ corresponding to the stabilizer of some rank one vertex $x\in T$ (this does not depend on the choice of $x$ by Lemma~\ref{lemma:intersection-standard-subtrees}). One readily checks that $\Phi'$ is $G_\Gamma$-equivariant, preserves adjacency,
and that $\Phi$ and $\Phi'$ are inverse of each other. The proposition follows.
\end{proof}

\begin{cor}
	\label{cor:no-flip}
There does not exist any element of $G_{\Gamma}$ which stabilizes an edge of $\Theta_\Gamma$ but flips the endpoints of this edge.
\end{cor}

\begin{proof}
In view of Proposition~\ref{prop:curve complex}, it suffices to prove the corollary for $\Theta^{\mathrm{sas}}_\Gamma$ in place of $\Theta_\Gamma$. Let $H_1$ and $H_2$ be two commuting abelian subgroups. By Lemma~\ref{lemma:commuting-sas}, up to reversing the roles of $H_1$ and $H_2$, we can asume that either $H_1$ is of type (3) and $H_2$ is of type (4), or else that they are both of type (4). 

In the first case, by considering the homomorphism $G_{\Gamma}\to\mathbb Z$ which sends each standard generator to $1$, we see that a standard abelian subgroup of type (3) is not conjugate to one of type (4). Therefore $H_1$ and $H_2$ cannot be permuted by an element of $G_\Gamma$. 

In the second case,
up to conjugation, there exist $s,t\in V\Gamma$ with $m_{st}=2$ such that $H_1=\langle s\rangle$ and $H_2=\langle t\rangle$. If $g\in G_{\Gamma}$ is such that $g G_{st}g^{-1}=G_{st}$, 
then $g$ fixes the unique rank $2$ vertex of $\del$ fixed by $G_{st}$, so $g\in G_{st}$.
Therefore conjugation by $g$ preserves both $H_1$ and $H_2$, without flipping them.
\end{proof}

\begin{prop}
	\label{prop:3-cycle}
The graphs $\Theta'_\Gamma,\Theta_\Gamma$ and $\Theta^{\mathrm{sas}}_\Gamma$ do not contain any 3-cycle.
\end{prop}

\begin{proof}
In view of Proposition~\ref{prop:curve complex}, it is enough to prove the corollary for the graph $\Theta^{\mathrm{sas}}_\Gamma$. Suppose $\Theta^{\mathrm{sas}}_\Gamma$ has a 3-cycle. Let $H_1,H_2$ and $H_3$ be the rank one standard abelian subgroups associated with the vertices of this 3-cycle. Let $H$ be the subgroup of $G_{\Gamma}$ generated by $H_1,H_2$ and $H_3$. Then $H$ is abelian. As $G_{\Gamma}$ is torsion-free, the group $H$ is free abelian. As $G_{\Gamma}$ has cohomogical dimension at most $2$ \cite{CharneyDavis}, the group $H$ has rank 2. By Lemma~\ref{lemma:commuting-sas}, the subgroup $H_{12}$ generated by $H_1$ and $H_2$ is a rank 2 standard abelian subgroup. Since $H$ is of rank 2, the group $H_3$ has a finite index subgroup contained in $H_{12}$. By the claim given in the first paragraph of the proof of Lemma~\ref{lem:sas intersection}, we have $H_3\subseteq H_{12}$. Therefore $H_1,H_2$ and $H_3$ are three pairwise distinct rank one abelian subgroups contained in $H_{12}$, contradicting Lemma~\ref{lemma:rank-one-in-rank-two}.
\end{proof}

\subsection{Vertex and edge stabilizers of the fixed set graph}

We now provide a detailed study of stabilizers of vertices and edges of the graph $\Theta'_\Gamma$ (Lemmas~\ref{lemma:stab-vertices-theta} and~\ref{lem:intersection}, respectively). This will give us a characterization of adjacency in $\Theta'_\Gamma$. This will be crucial in order to check that the fixed set graph satisfies Assumption~4 from the definition of a curve-graph-like $G_\Gamma$-graph. 

\subsubsection{Vertex stabilizers}

\begin{lemma}\label{lemma:stab-vertices-theta}
Let $v$ be a vertex of $\Theta'_\Gamma$. 
\begin{enumerate}
\item If $v$ is of type I (associated to a standard tree $F_g$, where $g$ is a conjugate of a standard generator) then there exists $n\ge 0$ such that $\Stab(v)$ is isomorphic to $F_n\times \mathbb{Z}$, where $F_n$ is a free group of rank $n$, and the $\mathbb{Z}$ factor is generated by $g$. If in addition $v$ belongs to the core subgraph $\Theta_\Gamma$ , then $n\ge 2$.  
\item If $v$ is of type II, then $\Stab(v)$ has infinite center and is virtually isomorphic to $F_n\times\mathbb{Z}$ with $n\ge 2$.
\end{enumerate}
\end{lemma} 

\begin{proof}
When $v$ is of type I, it was proved by Martin and Przytycki \cite[Lemma~4.5]{MP} that $\Stab(v)$ is isomorphic to $F_n\times\mathbb{Z}$. The additional part of the lemma when $v$ belongs to the core subgraph $\Theta_\Gamma$ follows from \cite[Lemma~9]{crisp}.

When $v$ is of type II, the conclusion follows from the fact that when $\Gamma$ is an edge, the Artin group $G_\Gamma$ has infinite center and is virtually a direct sum of a free group by $\mathbb Z$, see e.g.\ \cite[Section~2]{crisp}.
\end{proof}

\subsubsection{Edge stabilizers}

\begin{lemma}
	\label{lem:Z2}
	Let $v_1$ and $v_2$ be two adjacent vertices in $\Theta'_\Gamma$ such that $v_1$ is of type I and $v_2$ is of type II. Then $\stab(v_1)\cap \stab(v_2)$ is a standard abelian subgroup of type (2), in particular it is isomorphic to $\mathbb Z^2$.
\end{lemma}

\begin{proof}
	Up to conjugation, we can assume that there are adjacent vertices $s,t$ of $\Gamma$ with $m_{st}\ge 3$ such that $\stab(v_2)=G_{st}$, and that there exists an element $g\in G_{\Gamma}$ (conjugate to a standard generator) such that $\stab(v_1)=\stab(F_g)$. Let $x_2\in\del$ be the rank $2$ vertex corresponding to $G_{st}$. As $d_\Theta(v_1,v_2)=1$, we have $x_2\in F_g$. Let $x_1\in F_g$ be a rank 1 vertex adjacent to $x_2$. By the second assertion of Lemma~\ref{lem:type 4}, up to conjugating $g$ by an element of $G_{st}$, we can assume that $x_1$ corresponds to the subgroup $G_s$ or $G_t$ (say $G_s$). The first assertion of Lemma~\ref{lemma:fix-sets} implies that $F_g=F_s$. 
	
	We have $\langle s,Z_{st}\rangle\subseteq \Stab(v_1)\cap\Stab(v_2)$, and we will now prove that this inclusion is in fact an equality. So let $h\in\Stab(v_1)\cap\Stab(v_2)$ -- in other words $h\in G_{st}\cap \stab(F_s)$. Then $hx_1\in F_s$, so $s$ fixes $hx_1$. This means that $shG_s=hG_s$, hence $h^{-1}sh\in G_s$, i.e.\ $h$ belongs to the normalizer of $s$ in $G_{st}$. Lemma~\ref{lem:normalizer} thus implies that $h\in \langle s,Z_{st}\rangle$, which concludes our proof.
\end{proof}

We are now in position to provide a full classification of edge stabilizers of the modified fixed set graph $\Theta'_\Gamma$.

\begin{lemma}
	\label{lem:intersection}
	 Let $v_1$ and $v_2$ be two distinct vertices of $\Theta'_\Gamma$. 
\begin{enumerate}
	\item If $d_{\Theta'_\Gamma}(v_1,v_2)=1$, then $\stab(v_1)\cap \stab(v_2)$ is a rank 2 standard abelian subgroup.
	\item If $d_{\Theta'_\Gamma}(v_1,v_2)=2$, then $\stab(v_1)\cap \stab(v_2)$ is a rank 1 standard abelian subgroup.
	\item If $d_{\Theta'_\Gamma}(v_1,v_2)\ge 3$, then $\stab(v_1)\cap \stab(v_2)$ is the trivial subgroup.
\end{enumerate}
\end{lemma}

\begin{proof}
We will classify all possibilities for $\stab(v_1)\cap\stab(v_2)$; the conclusion will follow from our classification. Suppose $\stab(v_1)\cap\stab(v_2)$ contains a nontrivial element $g$ (of infinite order since $G_{\Gamma}$ is torsion-free). In the sequel, when $v_i$ is a type II vertex of $\Theta'_\Gamma$, we will denote by $x_i\in\del$ the corresponding rank $2$ vertex.

\medskip
\noindent
\emph{Case 1: Both $v_1$ and $v_2$ are of type II}. Then the fixed point set $F_g$ of $g$ in $\del$ contains at least two distinct points, namely $x_1$ and $x_2$. Lemma~\ref{lemma:fix-sets} thus implies that $g$ is conjugate to a standard generator, and $F_g$ is a standard tree. Then $F_g$ gives rise to a vertex in $\Theta'_\Gamma$ which is adjacent to both $x_1$ and $x_2$. We have $d_{\Theta'_\Gamma}(v_1,v_2)=2$: indeed, this distance cannot be equal to $1$ because $\Theta'_\Gamma$ does not contain any $3$-cycle (Proposition~\ref{prop:3-cycle}). It follows from the bipartite structure of $F_g$ that the geodesic segment $\overline{x_1x_2}\subseteq F_g$ contains a rank 1 vertex $x$, and we have $\stab(v_1)\cap \stab(v_2)\subseteq\Stab(x)$. Conversely, for every nontrivial element $h\in \stab(x)$, we have $x\in F_{g}\cap F_h$, and Lemma~\ref{lemma:intersection-standard-subtrees} yields $F_{g}=F_h$. Thus $h$ fixes $x_1$ and $x_2$. This shows that $\stab(x)\subseteq \stab(v_1)\cap \stab(v_2)$. In summary $\stab(v_1)\cap \stab(v_2)=\stab(x)$ is a standard abelian subgroup isomorphic to $\mathbb Z$. 

\medskip
\noindent
\emph{Case 2: The vertex $v_1$ is of type I, and $v_2$ is of type II}. Let $F_1$ be the standard tree associated with $v_1$, and let $h$ be a conjugate of a standard generator such that $F_1=F_h$. If $x_2\in F_1$, then $d_{\Theta'_\Gamma}(v_1,v_2)=1$ and $\stab(v_1)\cap \stab(v_2)$ is a standard abelian subgroup isomorphic to $\mathbb Z^2$ by Lemma~\ref{lem:Z2}.

Suppose $x_2\notin F_1$. Let $x$ be the point in $F_1$ closest to $x_2$ (such a point is unique because $\del$ is $\CAT$). Since $gx_2=x_2$ and $gF_1=F_1$, we have $gx=x$. Hence $g$ fixes the geodesic segment between $x_2$ and $x$ pointwise (and this segment is not reduced to a point because $x_2\notin F_1$). Using Lemma~\ref{lemma:fix-sets}, we get the existence of a conjugate $\ell$ of a standard generator such that $g=\ell^m$ for some $m\neq 0$. Let $v'$ be the vertex of $\Theta'_\Gamma$ associated with $F_\ell$. As $x_2\in F_\ell$, we have $d_{\Theta'_\Gamma}(v',v_2)= 1$.
In addition, as $g$ preserves $F_1=F_h$, the elements $g$ and $h$ commute by Lemma~\ref{lemma:stab-vertices-theta}, so $h$ and $\ell$ virtually commute. It follows from \cite[Theorem~3]{godelle} that $h$ and $\ell$ actually commute, from which we deduce that $d_{\Theta'_\Gamma}(v',v_1)=1$ (notice that $v'=v_1$ is ruled out as we are assuming that $v_2\notin F_1$). Therefore $d_{\Theta'_\Gamma}(v_1,v_2)= 2$. We then see that $\stab(v_1)\cap \stab(v_2)$ is equal to $\langle \ell\rangle$, which is a standard abelian subgroup isomorphic to $\mathbb Z$.

\medskip
\noindent
\emph{Case 3: Both $v_1$ and $v_2$ are of type I}. For every $i\in\{1,2\}$, let $F_i$ be the standard tree associated to $v_i$, and let $h_i$ be a conjugate of a standard generator such that $F_i=F_{h_i}$. 

We first assume that $F_1\cap F_2\neq\emptyset$. Then $F_1\cap F_2$ is reduced to a single rank $2$ vertex $x$ (Lemma~\ref{lemma:intersection-standard-subtrees}).
Up to conjugation we assume that $\stab(x)=G_{st}$ where $s,t\in V\Gamma$ satisfy $m_{st}<\infty$. For every $i\in\{1,2\}$, let $y_i\in F_i$ be a rank 1 vertex adjacent to $x$. Then for every $i\in\{1,2\}$, there exist $r_i\in\{s,t\}$ and $g_i\in G_{st}$ such that $y_i$ corresponds to the coset $g_iG_{r_i}$. Lemma~\ref{lemma:intersection-standard-subtrees} shows that $F_i=F_{g_ir_ig^{-1}_i}$. Hence we assume $h_i=g_ir_ig^{-1}_i$. 

If $m_{st}=2$, then $h_1$ and $h_2$ commute. Hence $d_{\Theta'_\Gamma}(v_1,v_2)=1$. Moreover, we have $\stab(v_1)\cap\stab(v_2)=\stab(F_1)\cap\stab(F_2)=\stab(x)=G_{st}$. 

If $m_{st}\ge 3$, then the vertex $v\in V\Theta'_\Gamma$ associated to $x$ is adjacent to both $v_1$ and $v_2$, so $d_{\Theta'_\Gamma}(v_1,v_2)=2$. Then $$\stab(v_1)\cap \stab(v_2)=(\stab(v_1)\cap \stab(v))\cap (\stab(v_2)\cap \stab(v))=\langle h_1, Z_{st}\rangle\cap \langle h_2, Z_{st}\rangle$$ where the last equality follows from the argument in Lemma~\ref{lem:Z2}. As $F_1\neq F_2$, the second statement from Lemma~\ref{lem:normalizer} then implies that $\stab(v_1)\cap \stab(v_2)=\langle Z_{st}\rangle$.

We now assume that $F_1\cap F_2=\emptyset$. We first claim that $g$ acts elliptically on $\del$. Otherwise $g$ has an axis $\ell_g\subseteq \del$. By Lemma~\ref{lemma:stab-vertices-theta}, the element $g$ commutes with both $h_1$ and $h_2$. We now argue as in the last paragraph of the proof of \cite[Lemma 11]{crisp}. Namely, the axis $\ell_g$ has a \emph{parallel set} \cite[Theorem~II.6.8]{BH}, which in our case is of the form $P_g=\ell_g\times T$, where $T$ is a tree. The centralizer $C_G(g)$ of $g$ is contained in $\stab(P_g)$. As $G_{\Gamma}$ is of hyperbolic type, the tree $T$ is bounded, so the factor action $C_G(g)\actson T$ (well-defined by \cite[Theorem~II.6.8(5)]{BH}) fixes a point $u\in T$. For every $i\in\{1,2\}$, we have $h_i\in C_G(g)$, and $h_i$ acts elliptically on $\del$. Therefore $h_1$ and $h_2$ both fix the line $\ell=\ell_g\times\{u\}$ pointwise. Hence $\ell\subseteq F_{h_1}\cap F_{h_2}$, a contradiction. The claim follows.

Observe that $F_g\cap F_i\neq\emptyset$ for every $i\in\{1,2\}$: indeed $g$ fixes $F_g$ pointwise and preserves $F_i$, so $g$ fixes the nearest-point projection $x_i\in F_i$ of any point $x\in F_g$ (hence $x_i\in F_g\cap F_i$). 
Therefore $F_g$ cannot be reduced to a single rank $2$ vertex (otherwise $F_1\cap F_2\neq\emptyset$), so Lemma~\ref{lemma:fix-sets} ensures that $F_g$ is a standard subtree
(corresponding to a vertex $v\in V\Theta'_\Gamma$ of type I). As $g$ commutes with both $h_1$ and $h_2$, the vertex $v$ is adjacent to both $v_1$ and $v_2$, so 
$d_{\Theta'_\Gamma}(v_1,v_2)= 2$.

Take any nontrivial $g'\in\stab(v_1)\cap \stab(v_2)$ with $g'\neq g$. We claim that $F_{g'}=F_g$. By the previous discussion $F_{g'}$ is a standard subtree; let $v'\in V\Theta'_\Gamma$ be the corresponding vertex. The element $g'$ commutes with both $h_1$ and $h_2$. Suppose by contradiction that $F_{g'}\neq F_g$. As $\Theta'$ does not have any 3-cycle by Proposition~\ref{prop:3-cycle}, we have $d_{\Theta'}(v',v)=2$. So $g$ and $g'$ do not commute, and by \cite[Theorem~3]{godelle} they do not even have commuting powers. Let $h=h_1h_2$. On the one hand, $h$ commutes with both $g$ and $g'$, so the centralizer of $h$ is not virtually abelian. As $G_{\Gamma}$ is of hyperbolic type, $h$ has to act elliptically on $\del$ by \cite[Lemma 11]{crisp}. On the other hand, each $h_i$ stabilizes $F_g$ by commutation. Moreover, $h_1$ and $h_2$ act elliptically on $F_g$ with empty intersection of their fix point set (as $F_1\cap F_2=\emptyset$). As $F_g$ is a tree, the action of $h$ on $F_g$ is axial, a contradiction. The claim follows.

Thus any element in $\stab(v_1)\cap \stab(v_2)$ fixes $F_g$ pointwise, hence $\stab(v_1)\cap \stab(v_2)$ is equal to $\langle g\rangle$, which is a standard abelian subgroup.
\end{proof}

As a consequence, we get the following characterization of adjacency in $\Theta_\Gamma$, which verifies Assumption~4.(c) from the definition of a curve-graph-like graph  (Definition~\ref{de:curve-graph-like}).

\begin{cor}\label{cor:sas-intersection}
Let $G=G_{\Gamma}$ be a $2$-dimensional Artin group of hyperbolic type with defining graph $\Gamma$. Let $v_1$ and $v_2$ be two distinct vertices of $\Theta_\Gamma$. The following statements are equivalent.
\begin{enumerate}
\item The vertices $v_1$ and $v_2$ are adjacent in $\Theta_\Gamma$.
\item For every vertex $v$ of $\Theta_\Gamma$, if every element of $\Stab_G(v_1)\cap\Stab_G(v_2)$ has a power that fixes $v$, then $v\in\{v_1,v_2\}$.
\end{enumerate}
\end{cor}

In the following proof, we will use the the simple observation that for two vertices $v_1,v_2\in V\Theta_\Gamma$, one has  $d_{\Theta_\Gamma}(v_1,v_2)=d_{\Theta'_\Gamma}(v_1,v_2)$.

\begin{proof}
First assume that $v_1$ and $v_2$ are adjacent in $\Theta_\Gamma$, and let $v$ be a vertex such that every element of $\Stab_G(v_1)\cap\Stab_G(v_2)$ has a power that fixes $v$. By Lemma~\ref{lem:intersection}, the intersection $\Stab_G(v_1)\cap\Stab_G(v_2)$ is a standard abelian subgroup of rank $2$. In particular, for every $i\in\{1,2\}$, the intersection $\Stab_G(v_i)\cap\Stab_G(v)$ contains a free abelian subgroup of rank $2$. Lemma~\ref{lem:intersection} thus implies that $d_{\Theta_\Gamma}(v,v_1)\le 1$ and $d_{\Theta_\Gamma}(v,v_2)\le 1$. As $\Theta_\Gamma$ contains no cycle of length $3$ (Proposition~\ref{prop:3-cycle}), it follows that $v=v_1$ or $v=v_2$.

Now assume that $v_1$ and $v_2$ are not adjacent in $\Theta_\Gamma$. If $d_{\Theta_\Gamma}(v_1,v_2)\ge 3$, then by Lemma~\ref{lem:intersection}, the intersection $\Stab_G(v_1)\cap\Stab_G(v_2)$ is trivial, so it fixes every vertex of $\Theta_\Gamma$. If $d_{\Theta_\Gamma}(v_1,v_2)=2$, then $\Stab_G(v_1)\cap\Stab_G(v_2)$ is a rank $1$ standard abelian subgroup $H$, whose fixed point set in $\Theta'_\Gamma$ contains all vertices adjacent to $v_H$ (where $v_H$ is the vertex of $\Theta'_\Gamma$ associated to $H$ via Proposition~\ref{prop:curve complex}). By Lemmas~\ref{lemma:commuting-sas} and~\ref{lemma:stab-vertices-theta}, $v_H$ is adjacent to both $v_1$ and $v_2$ in $\Theta'_\Gamma$, hence $v_H\in V\Theta_\Gamma$. Moreover, $v_H$ has infinite valence in $\Theta_\Gamma$ (for example, one can consider the orbit of $v_1$ under the action of the standard abelian $\mathbb Z$-subgroup associated with $v_2$, which gives infinitely many vertices in the link of $v_H$ by Lemma~\ref{lemma:stab-vertices-theta}). Thus $\Stab_G(v_1)\cap\Stab_G(v_2)$ fixes more than two vertices.
\end{proof}

\subsection{An analogue of canonical reduction systems}\label{sec:canonical-vertices}

The results in the present section will be used to check Assumption~3 (Canonical reduction systems) from the definition of a curve-graph-like $G_\Gamma$-graph (Definition~\ref{de:curve-graph-like}). Recall that $\calp_{\flex}(V(\Theta_\Gamma))$ denotes the collection of all nonempty sets $K$ of vertices of $\Theta_\Gamma$ whose pointwise stabilizer $\pstab(K)$ is infinite.

\begin{lemma}\label{lemma:pstab-sas}
Let $K\in \calp_{\flex}(V(\Theta_\Gamma))$ with $|K|>1$. Then $\pstab(K)$ is a standard abelian subgroup. 
\end{lemma}

\begin{proof}
Write $K$ as an increasing union of finite subsets $K_n$ of cardinality at least $2$. We claim that for every $n\in\mathbb{N}$, the pointwise stabilizer $\pstab(K_n)$ is a standard abelian subgroup. Indeed, this follows from Lemma~\ref{lem:intersection} if $|K_n|=2$. As the intersection of two standard abelian subgroups is again a standard abelian subgroup (Lemma~\ref{lem:sas intersection}), the claim then follows by induction on $|K_n|$. 

As $n$ increases, the groups $\pstab(K_n)$ thus form a non-increasing chain of standard abelian subgroups. By Corollary~\ref{cor:chain}, this chain is eventually constant, with value a standard abelian subgroup equal to $\pstab(K)$.
\end{proof}

The following lemma ensures that the fixed set graph satisfies Assumption~3.(b) from the definition of a curve-graph-like $G_\Gamma$-graph.

\begin{lemma}
	\label{lem:infinite stabilizer}
Let $(K_i)_{i\in\mathbb{N}}$ be a family of elements in $\calp_{\flex}(V(\Theta_\Gamma))$ such that for every $i\in\mathbb{N}$, one has $\pstab(K_{i+1})\subseteq\pstab(K_i)$. Then there exist at most three indices $i\in\mathbb{N}$ such that $\pstab(K_{i+1})\subsetneq\pstab(K_{i})$.
\end{lemma}

\begin{proof}
Without loss of generality, we can assume that the $K_i$ are pairwise distinct. Notice that if $v$ is a vertex of $\Theta_\Gamma$ and $K$ is any set of vertices of $\Theta_\Gamma$ distinct from $\{v\}$, then $\stab(v)$ is never contained in $\stab(K)$. Therefore, for every $i\ge 2$, we have $|K_i|\ge 2$. Using Lemma~\ref{lemma:pstab-sas}, this implies that the groups $\pstab(K_i)$ form a non-increasing chain of standard abelian subgroups, so by Corollary~\ref{cor:chain} there are at most $2$ indices $i\ge 2$ such that $\pstab(K_{i+1})\subsetneq\pstab(K_{i})$.
\end{proof}

The following lemma ensures that the fixed set graph satisfies Assumption~3.(a) from the definition of a curve-graph-like $G_\Gamma$-graph.

\begin{lemma}\label{lemma:theta-preorder}
	There exists a $G_{\Gamma}$-equivariant map $f:\calp_{\flex}(V(\Theta_\Gamma))\to V(\Theta_\Gamma)$.
\end{lemma}

\begin{proof}
Let $K\in \calp_{\flex}(V(\Theta_\Gamma))$. If $|K|=1$, then $f(K)$ is defined to be the unique vertex in $K$. If $|K|>1$, then by Lemma~\ref{lemma:pstab-sas}, either $\pstab(K)$ is a standard abelian subgroup of rank 2, which corresponds to an edge in $\Theta_\Gamma$ by Proposition~\ref{prop:curve complex} and we define $f(K)$ to be a vertex of this edge; or $\pstab(K)$ is a standard abelian subgroup of rank 1, which corresponds to a vertex in $\Theta_\Gamma$ and this vertex is defined to be $f(K)$. The choice of vertex in the edge in the former case can be made in a $G_{\Gamma}$-equivariant way by Corollary~\ref{cor:no-flip}. Now the proposition follows.
\end{proof}

\subsection{The fixed set graph is curve-graph-like.}

Using all the results establish in the section, we can finally show that the graph $\Theta_\Gamma$ is a curve-graph-like $G_\Gamma$-graph in the sense of Definition~\ref{de:curve-graph-like}.

\begin{prop}\label{prop:curve-graph-like}
Let $G=G_\Gamma$ be a $2$-dimensional Artin group of hyperbolic type whose defining graph $\Gamma$ is connected and not an edge labeled by $2$.
Then the fixed set graph $\Theta_\Gamma$ is a curve-graph-like $G$-graph. 
\end{prop}

\begin{proof}
Let $\calp_\ast(V(\Theta_\Gamma))=V(\Theta_\Gamma)$. Assumption~1 of Definition~\ref{de:curve-graph-like} is obviously satisfied.

We now check Assumption~2.  Equip the modified Deligne complex with a $\cat$ metric with finitely many isometry types of simplices as in Section~\ref{sec:background-deligne}, and let $K=\overline{\del}^h$ be its horofunction compactification. Write $K=K_{\bdd}\dunion K_\infty$ as in Proposition~\ref{prop:horo-partition}. Let $\Delta:=\partial_\infty\del$, which is separable and completely metrizable (see \cite[Proposition~5.31]{Vai}). By Proposition~\ref{prop:horo-partition}, there is a Borel $G$-equivariant map $K_{\bdd}\to \del$, and by Lemma~\ref{lemma:equivariant-map}, there is a Borel $G$-equivariant map $\del\to V(\Theta_\Gamma)$, showing that Assumption~2.(a) holds. By Proposition~\ref{prop:horo-partition}, there is a $G$-equivariant  homeomorphism $K_\infty\to\partial_\infty\del$, showing that Assumption~2.(b) holds. By Proposition~\ref{prop:barycenter}, there is a $G$-equivariant Borel map $(\partial_\infty\del)^{(3)}\to\del$, which combined with the Borel $G$-equivariant map $\del\to V(\Theta_\Gamma)$ provided by Lemma~\ref{lemma:equivariant-map} shows that Assumption~2.(c) holds.  Finally, the $G$-action on $\partial_\infty\del$ is Borel amenable by Theorem~\ref{theo:amenability-boundary-artin}.

Assumption~3.(a) was checked in Lemma~\ref{lemma:theta-preorder}, and Assumption~3.(b) was checked in Lemma~\ref{lem:infinite stabilizer}.
Assumption~4.(a) is provided by Lemma~\ref{lemma:stab-vertices-theta}, and Assumption~4.(b) follows from Lemmas~\ref{lemma:stab-vertices-theta} and~\ref{lem:intersection}. Finally, Assumption~4.(c) was checked in Corollary~\ref{cor:sas-intersection}.
\end{proof}

\section{Comparison between complexes, combinatorial rigidity}\label{sec:combinatorial-rigidity}

\emph{Any $2$-dimensional hyperbolic group of hyperbolic type naturally acts on three complexes: its Cayley complex, its Deligne complex, and the fixed set graph. We now explain how to compare the automorphism groups of these complexes, and  summarize several useful statements from Crisp's work \cite{crisp} that provide combinatorial rigidity theorems.} 

\subsection{Comparison maps}

For an Artin group $G=G_\Gamma$, let $\Cay_\Gamma$ be the Cayley complex of $G_\Gamma$, i.e.\ the universal cover of its presentation complex $P_\Gamma$ (for the standard Artin presentation of $G_\Gamma$). Edges of $\Cay_\Gamma$ are labeled by vertices of $\Gamma$, and have the usual orientations as edges in the Cayley graph.

The goal of this section is to compare the groups of cellular automorphisms of the Cayley complex, the Deligne complex and the fixed set graph of $G_\Gamma$. 
Whenever $X$ is a polyhedral complex with countably many cells, the group $\Aut(X)$ of all cellular automorphisms of $X$ is equipped with the compact-open topology. When $X$ is \emph{simple} in the sense that every cell is determined by its $0$-skeleton (which is the case in all the examples we consider), every element in $\Aut(X)$ is determined by its restriction to $V(X)$, and the compact-open topology coincides with the topology of pointwise convergence on $V(X)$. As $V(X)$ is countable, this topology is metrizable. When $X$ is locally compact (e.g.\ when $X=\Cay_\Gamma$), the group $\Aut(X)$ is second countable and locally compact, see e.g.\ \cite[Example~5.B.11]{CH}.

\paragraph*{Automorphisms of the Cayley complex.} 

 We identify $G_\Gamma$ with the vertex set $V(\Cay_\Gamma)$. 
A \emph{block} of $\Cay_\Gamma$ is the full subcomplex $B_{g,e}$ spanned by all vertices in a coset of the form $gG_e$, where $g\in G_\Gamma$ and $e\subset\Gamma$ is an edge; it is a \emph{large block} if $e$ has label at least $3$, and a \emph{small block} if $e$ has label $2$. 
A \emph{standard line} of $\Cay_\Gamma$ is a line which covers a circle in $P_\Gamma$ associated with a generator (corresponding to a vertex $v$ of $\Gamma$). It is \emph{semi-solid} if either $v$ has valence at least $2$, or $v$ belongs to an edge with label at least $3$. It is \emph{solid} if $v$ has valence at least $2$. A block is \emph{semi-solid} if each standard line in this block is semi-solid.

\begin{lemma}\label{lemma:automorphisms-cayley}
Let $G_\Gamma$ be a $2$-dimensional Artin group of hyperbolic type. Then every cellular automorphism of $\Cay_\Gamma$ sends semi-solid standard lines to semi-solid standard lines, solid standard lines to solid standard lines, semi-solid blocks to semi-solid blocks and large blocks to large blocks.
\end{lemma}

\begin{proof}
Let $\alpha\in\Aut(\Cay_\Gamma)$. We first observe that the following two properties of an edge $e$ of $\Cay_\Gamma$ are preserved by $\alpha$ (i.e.\ if $e$ has the property, then so does $\alpha(e)$):
\begin{enumerate}
\item the label $a$ of $e$ (a vertex of $\Gamma$) belongs to an edge of $\Gamma$ with label at least $3$;
\item the label $a$ of $e$ has valence at least $2$ in $\Gamma$.
\end{enumerate}
Indeed, the first property is characterized by the existence of a $2$-cell whose boundary has at least $6$ edges, including $e$. The second property is characterized by the existence of three pairwise distinct edges in $\Cay_\Gamma$ with the same origin as $e$, each sharing a common $2$-cell with $e$. 

We now prove that $\alpha$ preserves the set of semi-solid standard lines and the set of solid standard lines. In view of the above observations, it is enough to show that $\alpha$ sends any semi-solid standard line $\ell$ to a standard line. For this, it is enough to prove that whenever $e_1$ and $e_2$ are two consecutive edges in $\ell$, then $\alpha(e_1)$ and $\alpha(e_2)$ have the same label. From now on, we fix two such consecutive edges $e_1,e_2$, we denote by $x$ their intersection point, and let $x'=\alpha(x)$ and $e'_i=\alpha(e_i)$; we also denote by $a\in V\Gamma$ the common label of all edges in $\ell$.

We first assume that $a$ is contained in an edge $\varepsilon$ of $\Gamma$ with label at least $3$. Then we can find two $2$-cells $C_1$ and $C_2$, both corresponding to the relator associated to $\varepsilon$, with $e_1\subseteq C_1,e_2\subseteq C_2$, and such that $C_1\cap C_2$ is a path starting from $x$ with $n-1$ edges. Thus $\alpha(C_1)\cap\alpha(C_2)$ is a path starting from $x'$ with $n-1$ edges. As $n-1\ge 2$, the form of the relators of $G_\Gamma$ implies that $e'_1$ and $e'_2$ are labeled by the same vertex of $\Gamma$, as desired.

Now suppose that $a$ has valence at least $2$ in $\Gamma$, and every incident edge has label $2$ (we say that $a$ is \emph{special}). Let $\Lambda$ be the link of $x$ in $\Cay_\Gamma$. Each vertex $v\in\Gamma$ gives a pair of vertices $v^{\pm}$ of $\Lambda$. 

We claim that if $v_1$ and $v_2$ are two distinct vertices in $\Lambda$, both corresponding to edges of $\Cay_\Gamma$ whose label is a special vertex of $\Gamma$, and such that $\lk(v_1,\Lambda)=\lk(v_2,\Lambda)$ and $|\lk(v_1,\Lambda)|\ge 4$, then there exists $v\in\Gamma$ such that $\{v_1,v_2\}=\{v^-,v^+\}$. Indeed, let $a_i\in\Gamma$ be the vertex associated with $v_i$. Our assumption implies that $\lk(a_1,\Gamma)=\lk(a_2,\Gamma)$ and $|\lk(a_1,\Gamma)|\ge 2$. If $a_1\neq a_2$, then there is an induced 4-cycle in $\Gamma$ (containing $a_1$ and $a_2$) such that each edge of the 4-cycle is labeled by $2$, violating that $\Gamma$ is of hyperbolic type. Thus $a_1=a_2$ and the claim follows. 

By our assumption on $a$, the pair of vertices in $\Lambda$ determined by $e_1$ and $e_2$ satisfies the assumptions of the claim. As $\alpha$ induces an isomorphism of links at each vertex and sends edges labeled by special vertices to edges labeled by special vertices (Observations~(1) and~(2) above), the claim implies that $e'_1$ and $e'_2$ have the same label, as desired.

Now we show that  semi-solid blocks are preserved (and as small blocks are grids, small semi-solid blocks and large blocks will automatically be preserved).
 Let $a,b\in V\Gamma$ be joined by an edge, and let $n$ be the label of this edge.
Let $\ell_a$ and $\ell_b$ be standard lines labeled by $a$ and $b$ respectively such that $\ell_a\cap \ell_b=\{x\}$ is a vertex of $\Cay_\Gamma$. We claim that for any standard $b$-line $\ell$ intersecting $\ell_a$ in a point, the images $\alpha(\ell_b)$ and $\alpha(\ell)$ are standard lines with the same label. 
\begin{figure}
	\centering
	\label{figure:1}
\includegraphics[scale=0.8]{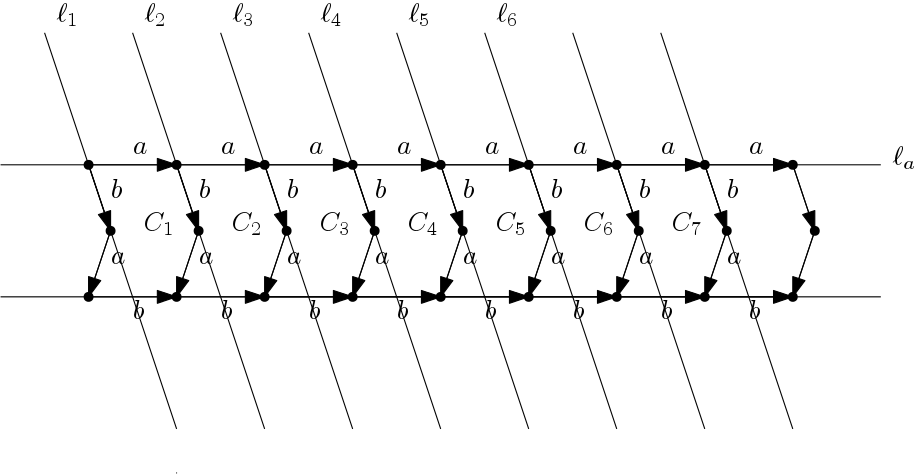}
\caption{Cellular automorphisms of $\Cay_\Gamma$ send large blocks to large blocks.}
\end{figure}
Indeed, let $\{\ell_i\}_{i\in\mathbb Z}$ be the family of all consecutive $b$-lines interesting $\ell_a$ (see Figure~\ref{figure:1}). Let $\{C_i\}_{i\in\mathbb Z}$ be 2-cells such that every $C_i$ intersects each of $\ell_a$, $\ell_i$ and $\ell_{i+1}$ in an edge and $C_{i}\cap C_{i+1}$ has $n-1$ edges (see Figure~\ref{figure:1} for the case $n=3$). As $\ell_a$ is semi-solid, $\alpha(\ell_a)$ is a standard line, thus the edges of $\alpha(C_i)$ and $\alpha(C_{i+1})$ contained in $\alpha(\ell_a)$ have the same label, say $c$. Since $\alpha(C_i)\cap\alpha(C_{i+1})$ contains at least one edge with label different from $c$, we  deduce that $\alpha(C_i)$ and $\alpha(C_{i+1})$ are contained in the same block. Hence cells $\alpha(C_i)$ (with $i\in\mathbb Z$) are all contained in the same block $B'$. Note that $B'$ is semi-solid.
 As each $\alpha(\ell_i)$ is a standard line intersecting $\alpha(C_i)$ in an edge, $\alpha(\ell_i)\subset B'$ and all lines $\alpha(\ell_i)$ have the same label, thus showing our claim. 

As the 1-skeleton of a  semi-solid block is a union of semi-solid standard lines, by repeatedly applying the claim (with possibly interchanging the role of $a$ and $b$), we deduce that $\alpha$ sends  semi-solid blocks to semi-solid blocks.
\end{proof}

\paragraph*{From $\Aut(\Cay_\Gamma)$ to $\Aut(\Theta_\Gamma)$.}

To build a homomorphism $h:\Aut(\Cay_\Gamma)\to\Aut(\Theta_\Gamma)$, we now make a few additional observations. First, notice that the $G_\Gamma$-stabilizer of a large block $B_{g,e}$ is equal to $gG_eg^{-1}$, in particular its center $gZ_eg^{-1}$ is a standard abelian subgroup (of type (3), in the sense of Definition~\ref{de:sas}), which corresponds to a vertex of $\Theta_\Gamma$. Distinct large blocks correspond to distinct vertices of $\Theta_\Gamma$. 

Second, the stabilizer of each standard line $\ell$ is a rank one standard abelian subgroup (of type (4), in the sense of Definition~\ref{de:sas}), which always corresponds to a vertex of $\Theta'_\Gamma$. Two standard lines are \emph{parallel} if they have finite Hausdorff distance. The \emph{parallel set} of a standard line $\ell$ is the union of all standard lines parallel to $\ell$. Note that two standard lines are parallel if and only if their stabilizers have finite Hausdorff distance in the Cayley graph of $G_\Gamma$; by \cite[Corollary~2.4]{mosher2011quasi}, this happens if and only if their stabilizers are commensurable. Using Lemma~\ref{lem:sas intersection}, we deduce that two standard lines are parallel if and only if they have the same stabilizer.

In summary, denoting by $\calb_{\parallel}(\Cay_\Gamma)$ the set of all large blocks and all parallel sets of standard lines in $\Cay_\Gamma$, we have constructed a $G_\Gamma$-equivariant bijection $$\pi:\calb_{\parallel}(\Cay_\Gamma)\to V(\Theta'_\Gamma)$$ sending a large block to the center of its stabilizer, and a parallel set of standard lines to the stabilizer of any line in the set. Note that in the latter case, the collection of standard lines in the parallel set $B$ is in 1-1 correspondence with the collection of rank 1 vertices in the standard tree associated with $\pi(B)$. Thus for $B\in \calb_{\parallel}(\Cay_\Gamma)$, one has $\pi(B)\in V(\Theta_\Gamma)$ if and only if $B$ is a large block or $B$ is the parallel set of a solid standard line, see Remark~\ref{rk:theta-theta'}.\footnote{Here we would like to warn the reader that a solid standard line can be parallel to a non-solid standard line: this happens precisely when the generator associated to a valence $1$ vertex of $\Gamma$ is conjugate to a generator associated to another vertex of valence at least $2$.}

\begin{lemma}\label{lemma:cay-to-theta}
Let $G_\Gamma$ be a $2$-dimensional Artin group of hyperbolic type. Then there exists a unique homomorphism $h:\Aut(\Cay_\Gamma)\to\Aut(\Theta_\Gamma)$ such that for every $\alpha\in\Aut(\Cay_\Gamma)$ and every $B\in \calb_{\parallel}(\Cay_\Gamma)$, if $\pi(B)\in V(\Theta_\Gamma)$, then $h(\alpha)(\pi(B))=\pi(\alpha(B))$.
\end{lemma}

We call $h$ the \emph{comparison homomorphism} between $\Aut(\Cay_\Gamma)$ and $\Aut(\Theta_\Gamma)$.

\begin{proof}
The uniqueness of $h$ is clear since every vertex of $\Theta_\Gamma$ is in the image of $\pi$. For its existence, by Lemma~\ref{lemma:automorphisms-cayley}, every $\alpha\in\Aut(\Cay_\Gamma)$ sends solid standard lines to solid standard lines  (sending parallel solid standard lines to parallel solid standard lines) and large blocks to large blocks, and therefore induces a bijection from $V(\Theta_\Gamma)$ to itself. It thus suffices to observe that this bijection preserves adjacency (and non-adjacency) in $\Theta_\Gamma$.

This follows from the following observations. If $B_1$ and $B_2$ are large blocks, then $\pi(B_1)$ and $\pi(B_2)$ are never adjacent. If $B_1$ is a large block and $B_2$ is a parallel set of standard lines, then $\pi(B_1)$ and $\pi(B_2)$ are adjacent if and only if some line of $B_2$ is contained in $B_1$. Finally, if $B_1$ and $B_2$ are two parallel sets of standard lines, then $\pi(B_1)$ and $\pi(B_2)$ are adjacent if and only if there is a small block $B$ such that $B_1$ contains all horizontal standard lines of $B$ and $B_2$ contains all vertical standard lines of $B$. These properties are preserved by cellular automorphisms of $\Cay_\Gamma$, which concludes our proof. 
\end{proof}

\paragraph*{From $\Aut(\Cay_\Gamma)$ to $\Aut(\del)$.}
 From now on we will further assume that every vertex of $\Gamma$ is semi-solid (i.e.\ it either has valence at least $2$, or it has valence $1$ and is adjacent to an edge with label at least $3$). Under this extra assumption, we define a homomorphism $$h_1:\Aut(\Cay_\Gamma)\to\Aut(\del)$$ as follows. Recall that a vertex of $\del$ corresponds to either a vertex of $\Cay_\Gamma$, a standard line of $\Cay_\Gamma$, or a block of $\Cay_\Gamma$. As every vertex of $\Gamma$ is semi-solid, Lemma~\ref{lemma:automorphisms-cayley} implies that $\alpha$ maps standard lines to standard lines, and blocks to blocks; so $\alpha$ induces a permutation of $V(\del)$ which extends to a unique element in $\Aut(\del)$. One readily verifies that this map from $\Aut(\Cay_\Gamma)$ to $\Aut(\del)$ is a group homomorphism.

\begin{lemma}\label{lemma:h1}
Let $G_\Gamma$ be a $2$-dimensional Artin group of hyperbolic type, and assume that every vertex of $\Gamma$ is semi-solid.

Then $h_1:\Aut(\Cay_\Gamma)\to\Aut(\del)$ is a homeomorphism onto its image. 
\end{lemma}

\begin{proof}
Injectivity of $h_1$ is clear as we can identify vertices of $\Cay_\Gamma$ with rank 0 vertices of $\del$. As the compact-open topology on either $\Aut(\Cay_\Gamma)$ or $\Aut(\del)$ is metrizable, it suffices to show that a sequence $(\alpha_n)_{n\in\mathbb{N}}\in\Aut(\Cay_\Gamma)^{\mathbb{N}}$ is eventually identity on any finite subset of $V(\Cay_\Gamma)$ if and only if the sequence $(h_1(\alpha_n))_{n\in\mathbb{N}}$ is eventually identity on any finite subset of $V(\del)$. The if direction is clear. For the only if direction, $\alpha_n$ being eventually identity on any finite subset of $V(\Cay_\Gamma)$ implies that $h_1(\alpha_n)$ is eventually identity on any finite subset of the set of rank 0 vertices of $\del$.  Now for every vertex $v\in\del$ of rank $k\in\{1,2\}$, there exists a finite collection $C$ of elements of $G_\Gamma\cong V(\Cay_{\Gamma})$ such that $v$ is the unique vertex of $\del$ of rank $k$ whose corresponding coset of $G_\Gamma$ contains $C$.  As automorphisms of $\del$ preserve ranks of vertices by Lemma~\ref{lemma:aut deligne}, any automorphism of $\del$ fixing $C$ pointwise will also fix $v$, which implies $h(\alpha_n)$ is eventually identity on any finite subset of $V(\del)$. 
\end{proof}

\paragraph*{From $\Aut(\del)$ to $\Aut(\Theta_\Gamma)$.}
Let us call a rank 2 vertex of $\del$ \emph{large} (resp.\ \emph{small}) if it is associated with a large (resp.\ small) coset of $G_\Gamma$.
\begin{lemma}
	\label{lem:standard tree}
	Let $G_\Gamma$ be a 2-dimensional Artin group of hyperbolic type, and let $\alpha\in\Aut(\del)$. Then 
	\begin{enumerate}
		\item $\alpha$ maps standard trees to standard trees;
		\item $\alpha$ preserves the set of large rank 2 vertices of $\del$.
	\end{enumerate}
\end{lemma}

\begin{proof}
Lemma~\ref{lemma:aut deligne} ensures that $\alpha$ preserves ranks of vertices.
Let $x\in\del$ be a rank 2 vertex associated with a coset $gG_e$ with $e\subset\Gamma$ being an edge.

We claim that two rank 1 vertices $x_1,x_2$ adjacent to $x$ belong to the same standard tree if and only if $\alpha(x_1)$ and $\alpha(x_2)$ belong to the same standard tree. Indeed, if $e$ is labeled by 2, then by the discussion of links of rank 2 vertices in the proof of Lemma~\ref{lemma:aut deligne}, $x_1$ and $x_2$ are in the same standard tree if and only if they correspond to non-adjacent vertices in the link of $x$. This is preserved by $\alpha$. Suppose now that the label of $e$ is at least $3$. Then \cite[Proposition~40]{crisp} implies that $\alpha$ restricted to $gG_e$ (which is identified as the collection of rank $0$ vertices adjacent to $x$) can be extended to a cellular map on the block of $\Cay_{\Gamma}$ containing $gG_e$, thus parallelism of the standard lines inside the blocks is preserved. But recall that two standard lines are parallel if and only if they have the same stabilizer, and this is precisely the condition for the corresponding rank one vertices of $\del$ to belong to the same standard tree. The claim thus follows.

We also note that for a rank 1 vertex $x$ in a standard tree $T$, the collection of rank 2 vertices in $T$ which are adjacent to $x$ coincides with the collection of rank 2 vertices of $\del$ adjacent to $x$. As vertices in a standard tree alternate between rank 1 and rank 2, we deduce that $\alpha$ sends a standard tree to a standard tree. 

The second statement follows from the discussion of vertex links of $\del$ in the proof of Lemma~\ref{lemma:aut deligne},
\end{proof}

Given any 2-dimensional Artin group of hyperbolic type $G_\Gamma$, Lemma~\ref{lem:standard tree} allows us to define a \emph{comparison homomorphism} $$h_2:\Aut(\del)\to\Aut(\Theta_\Gamma)$$
as follows. Take $\alpha\in \Aut(\del)$. Recall that a vertex of $\Theta'_\Gamma$ corresponds to either a standard tree or a rank 2 vertex of $\del$ associated with a large block, thus $\alpha$ induces a permutation of the vertex set of $\Theta'_\Gamma$ by Lemma~\ref{lem:standard tree}, which extends to an automorphism of $\Theta'_\Gamma$ by Definition~\ref{de:theta} (note that Lemma~\ref{lemma:aut deligne} and Lemma~\ref{lem:standard tree} imply a pair of standard trees intersecting in a small rank 2 vertex are sent by $\alpha$ to another pair of standard trees intersecting in a small rank 2 vertex). This automorphism restricts to an automorphism of $\Theta_\Gamma$, which gives a map $h_2:\Aut(\del)\to\Aut(\Theta_\Gamma)$. One readily verifies that this map is a homomorphism.

\begin{lemma}
	\label{lem:homos0}
	Let $G_\Gamma$ be a 2-dimensional Artin group of hyperbolic type. Suppose $\Gamma$ does not contain isolated vertices or valence one vertices.
	
	Then the map $h_2:\Aut(\del)\to \Aut(\Theta_\Gamma)$ is a homeomorphism onto its image. Moreover, $h=h_2\circ h_1$.
\end{lemma}

\begin{proof}
 We first show that $h_2$ is injective. Let $x\in\del$ be a rank 0 vertex, we also view $x$ as a vertex in $\Cay_\Gamma$ and an element of $G_\Gamma$. Let $\{\ell_i\}_{i=1}^k$ be the collection of standard lines containing $v$, let $P_i$ be the vertex set of the parallel set of $\ell_i$ and let $v_i$ be the vertex of $\Theta_\Gamma$ associated with $\ell_i$ (note that our assumption on $\Gamma$ implies that $\Theta_\Gamma=\Theta'_\Gamma$, see  Definition~\ref{de:theta}). For every $i\in\{1,\dots,k\}$, let $a_i\in V\Gamma$ be the label of $\ell_i$. 

We claim that $\cap_{i=1}^k P_i=\{x\}$. It suffices to prove this claim in the special case where $x$ is the identity element of $G_\Gamma$, and from now on we assume that we are in this special case. We first observe that $\Stab(v_i)=P_i$. Indeed, if $g\in P_i$, then $g\langle a_i\rangle$ and $\langle a_i\rangle$ have finite Hausdorff distance, which implies that $g\langle a_i\rangle g^{-1}\cap\langle a_i\rangle$ is of finite index in both $\langle a_i\rangle$ and $g\langle a_i\rangle g^{-1}$ by \cite[Corollary 2.4]{mosher2011quasi}. Hence $g\langle a_i\rangle g^{-1}=\langle a_i\rangle$ by Lemma~\ref{lem:sas intersection}. Thus $g\in \Stab(v_i)$ and $P_i\subset \Stab(v_i)$. The containment $\Stab(v_i)\subset P_i$ follows from definition. 
By Case 3 of the proof of Lemma~\ref{lem:intersection}, if $a_i$ and $a_j$ are adjacent in $\Gamma$, then $P_i\cap P_j$ is contained in the subgroup generated by $a_i$ and $a_j$. By the the moreover part of Lemma~\ref{lem:injective} the intersection of all subgroups $G_e$ with $e$ varying over the edges of $\Gamma$ is trivial, thus
$\cap_{i=1}^k P_i=\{x\}$.

Let now $\alpha\in\Aut(\del)$ be such that $h_2(\alpha)$ is the identity. Then $\alpha(T_i)=T_i$ where $T_i\subset\del$ is the standard tree associated with $v_i$. However, $P_i$ is the collection of rank 0 vertices of $\del$ which are adjacent to some rank 1 vertex in $T_i$. Then $\alpha(P_i)=P_i$. Hence $\alpha(x)=x$.

Now we show that $h_2$ is a homeomorphism onto its image. It suffices to show that a sequence $(\alpha_n)_{n\in\mathbb{N}}\in\Aut(\del)^{\mathbb{N}}$ is eventually identity on any finite subset of the vertex set of the underlying complex if and only if the same holds for the sequence $h_2(\alpha_n)$. The only if direction follows from the fact that for any standard tree $T$, there is a vertex $v$ of $T$ such that $T$ is the unique standard tree containing $v$ (see Lemma~\ref{lemma:intersection-standard-subtrees}). For the if direction, notice that each rank 0 vertex $v$ of $\del$ (viewed as an element of $G_\Gamma$) is an intersection of finitely many parallel sets of standard lines (as discussed in the previous paragraph). These parallel sets correspond to a finite set $S$ of vertices of $\Theta_\Gamma$. Hence any automorphism of $\del$ whose $h_2$-image fixes $S$ pointwise must fix $v$. Finally, as in the proof of Lemma~\ref{lemma:h1}, every vertex of $\del$ is determined (among vertices of the same rank) by finitely many adjacent rank $0$ vertices. This completes our proof that $h_2$ is a homeomorphism onto its image.  

The fact that $h=h_2\circ h_1$ follows from construction.
\end{proof}

The following is a consequence of Lemmas~\ref{lemma:h1} and~\ref{lem:homos0}.
\begin{lemma}\label{lemma:comparison-is-iso}
	Let $G_\Gamma$ be a $2$-dimensional Artin group of hyperbolic type such that every vertex of $\Gamma$ has valence at least $2$. Assume that the comparison map $\Aut(\Cay_\Gamma)\to\Aut(\Theta_\Gamma)$ is an isomorphism.
	
	Then $h,h_1,h_2$ are isomorphisms of topological groups. 
	\qed
\end{lemma}

\paragraph*{From $\Aut(\del)$ to $\Aut(\Cay_\Gamma)$.}
Suppose $G=G_\Gamma$ is an Artin group such that each vertex of $\Gamma$ is contained in a large edge. Then there is an injective homomorphism $$h'_1:\Aut(\del)\to \Aut(\Cay_\Gamma)$$ defined as follows. An automorphism of $\del$ gives a bijection $f:G\to G$. By \cite[Proposition~40]{crisp}, the restriction of $f$ to each large block can be extended to a unique cellular automorphism onto its image. By the assumption on $\Gamma$, each standard line of $\Cay_\Gamma$ is contained in a large block, so $f$ extends to a cellular automorphism. This gives $h'_1:\Aut(\del)\to \Aut(\Cay_\Gamma)$. Any element of $\Aut(\del)$ fixing rank 0 vertices of $\del$ pointwise will also fix other vertices as they correspond to certain left cosets of $G$. Thus $h'_1$ is injective.

\subsection{Rigidity statements}

\begin{theo}[Crisp \cite{crisp}]\label{theo:rigid-crisp}
Let $G_{\Gamma_1}$ and $G_{\Gamma_2}$ be two Artin groups whose defining graphs $\Gamma_1$ and $\Gamma_2$ are both connected, triangle-free, with all labels at least $3$, and have more than 2 vertices.
Suppose $\Gamma_1$ has no separating vertices and no separating edges. 

If $\Theta_{\Gamma_1}$ and $\Theta_{\Gamma_2}$ are isomorphic as graphs, then there is a label-preserving isomorphism $\Gamma_1\to\Gamma_2$, and in particular $G_{\Gamma_1}$ and $G_{\Gamma_2}$ are isomorphic as groups. 
\end{theo}

\begin{proof}
The theorem is not stated explicitly in \cite{crisp}, but follows from \cite[p.~1435]{crisp}. We identify type II vertices of $\Theta_{\Gamma_i}$ with rank 2 vertices of $\mathbb D_{\Gamma_i}$. By \cite[Proposition~41]{crisp}, the isomorphism $\alpha:\Theta_{\Gamma_1}\to\Theta_{\Gamma_2}$ maps type II vertices to type II vertices. The argument in \cite[p.~1435]{crisp} implies that $\alpha$ induces an simplicial isometric embedding $f:\mathbb D_{\Gamma_1}\to\mathbb D_{\Gamma_2}$ such that $f$ maps rank 0 vertices to rank 0 vertices, and $f$ coincides with $\alpha$ when restricted to rank 2 vertices. Let $C$ be the image of $f$. As $C$ contains all rank 2 vertices of $\mathbb D_{\Gamma_2}$, we know $\mathbb D_{\Gamma_2}$ is contained in a finite neighborhood of $C$. If there is a point $x\in D_{\Gamma_2}\setminus C$, then let $x'$ be the nearest point projection of $x$ to $C$ ($x'$ is well-defined as $C$ is a convex subset of a $\CAT$ space). As $D_{\Gamma_2}$ has geodesic extension property, we can extend the geodesic segment from $x'$ to $x$ to form a ray, which is not contained in a finite neighborhood of $C$, contradiction. Thus $f$ is surjective, hence $f$ is an isometry. By considering the restriction of $f$ to the link of a rank 0 vertex, we obtain a graph isomorphism $\Gamma_1\to \Gamma_2$, which is label preserving as explained in \cite[p.~1435]{crisp}.
\end{proof}

For the next statement, recall from the introdcution that $\Gamma$ is \emph{star-rigid} if the only label-preserving graph automorphism of $\Gamma$ which fixes the star of some vertex pointwise is the identity.

\begin{theo}[Crisp \cite{crisp}]\label{theo:rigid-crisp-2}
Let $G=G_{\Gamma}$ be an Artin group whose defining graph $\Gamma$ is connected, triangle-free, with all labels at least $3$, and has more than 2 vertices. Assume that $\Gamma$ has no separating vertex and no separating edge. 

Then the comparison homomorphism $h:\Aut(\Cay_\Gamma)\to\Aut(\Theta_\Gamma)$ is an isomorphism. If in addition $\Gamma$ is star-rigid, then the natural map $G\to\Aut(\Theta_\Gamma)$ is injective with finite index image.
\end{theo}

\begin{proof}
We identify $G$ as the collection of rank 0 vertices in $\del$. Next we define $h'_2:\Aut(\Theta_\Gamma)\to\Aut(\del)$.
Given $\beta\in \Aut(\Theta_\Gamma)$, we first define $f_\beta:G\to G$ as follows. All blocks of $G$ are large, and they are in 1-1 correspondence with type II vertices of $\Theta_\Gamma$. 
Take $g\in G$ and let $\{B_i\}_{i=1}^k$ be the collection of blocks containing $g$. Then $\cap_{i=1}^k B_i=\{g\}$. Let $B'_i$ be the block corresponding to $\beta (w_{B_i})$ where $w_{B_i}\in V\Theta_\Gamma$ is the vertex associated with $B_i$ (note that $\beta$ maps type II vertices to type II vertices by \cite[Proposition~41]{crisp}). It was explained in the first paragraph of the proof of \cite[Theorem 11.7]{HO} that one can deduce from Crisp's work \cite{crisp} that $\cap_{i=1}^k B'_i$ is exactly one element, which is defined to be $f_\beta(g)$. By applying the same argument to $\beta^{-1}$, we know $f_\beta$ is a bijection sending blocks to blocks. As cosets of $G$ associated to rank 1 vertices of $\Theta_\Gamma$ are exactly the nonempty intersections of two blocks, these cosets are preserved by $f_\beta$. Thus $f_\beta$ extends to an automorphism of $\del$. This actually gives a homomorphism $h'_2:\Aut(\Theta_\Gamma)\to\Aut(\del)$. Recall from the paragraph below Lemma~\ref{lemma:comparison-is-iso} that we also have a map $h'_1:\Aut(\del)\to\Aut(\Theta_\Gamma)$. One readily verifies that $h'_1\circ h'_2$ and $h$ are inverses of each other, which finishes the proof.

If $\Gamma$ is star-rigid, the discussion in \cite[p.~1436]{crisp} shows that any automorphism of $\del$ fixing the star of a rank 0 vertex is identity, hence any automorphism of $\Cay_\Gamma$ fixing a small neighborhood of a vertex is identity. Hence the action of $\Aut(\Cay_\Gamma)=\Aut(\Theta_\Gamma)$ on $\Cay_\Gamma$ is proper, which implies the required statement in the theorem.
\end{proof}

\section{Artin groups and the strong ICC property}\label{sec:icc}

\emph{In order to apply Theorem~\ref{theo:bfs}, we will need to know that the automorphism group of the Cayley complex of $G_\Gamma$ (which will play the role of the ambient locally compact second countable group) is strongly ICC. This section is devoted to the proof of this fact under appropriate assumptions.}

\subsection{Actions on $\cat$ spaces and the strong ICC property}

Recall that a Polish group $\mathsf{G}$ is \emph{strongly ICC} if there is no conjugation-invariant probability measure on  $\mathsf{G}\setminus\{e\}$. We will now provide criteria to establish the strong ICC property for groups acting on $\cat$ spaces. 
Recall that an isometric action of a group $\mathsf{G}$ on a Gromov hyperbolic space $X$ is \emph{non-elementary} if $\mathsf{G}$ contains two loxodromic isometries whose fixed point sets in $\partial_\infty X$ are disjoint.  Recall also that the \emph{limit set} of $\mathsf{G}$ is defined as $\Lambda_\infty \mathsf{G}=\overline{\mathsf{G}\cdot x_0}\cap\partial_\infty X$ for any $x_0\in X$.

\begin{prop}\label{prop:icc-hyperbolic}
Let $X$ be a complete separable $\cat$ space. Let $\mathsf{G}$ be a Polish group acting continuously (for the compact-open topology on $\Isom(X)$), non-elementarily by isometries on $X$. Let  $K\unlhd \mathsf{G}$ be the normal subgroup made of all automorphisms that act trivially on $\Lambda_\infty \mathsf{G}$.

Then every conjugation-invariant probability measure on $\mathsf{G}$ is supported on $K$. In particular, if $\mathsf{G}$ acts faithfully on $\Lambda_\infty \mathsf{G}$, then $\mathsf{G}$ is strongly ICC.
\end{prop}

\begin{proof}
The proof is inspired from \cite[Lemma~2.4]{BFS}.  
Let $\mu$ be a conjugation-invariant probability measure on $\mathsf{G}$. Let $\horox$ be the horofunction compactification of $X$, a compact metrizable space. There is a continuous homomorphism $\mathsf{G}\to\Homeo(\horox)$. The set $\Prob(\horox)$, equipped with the weak-$*$ topology, is compact. As in the proof of \cite[Lemma~2.4]{BFS}, the subset $\Prob_\mu(\horox)$ made of all $\mu$-stationary probability measures on $\horox$ is a non-empty $\mathsf{G}$-invariant closed (whence compact) subset of $\Prob(\horox)$. 

Using Proposition~\ref{prop:horo-partition}, we view $\partial_\infty X$ as a subspace of $\horox$. We claim that for every loxodromic isometry $g\in \mathsf{G}$, the Dirac mass at $g^{+\infty}$ belongs to $\Prob_\mu(\horox)$. Indeed, using the non-elementarity of the $\mathsf{G}$-action on $X$, we can find a loxodromic isometry $h\in \mathsf{G}$ such that $\{g^{-\infty},g^{+\infty}\}\cap\{h^{-\infty},h^{+\infty}\}=\emptyset$. Let $\nu\in\Prob_\mu(\horox)$. Combining Lemma~\ref{lemma:ns-dynamics-on-horo} with \cite[Lemma~8.3]{FLM}, we see that every weak-$*$ accumulation point $\nu_\infty$ of the measures $h^n\nu$ is supported on $\{h^{-\infty},h^{+\infty}\}$. Now, the only weak-$*$ accumulation point of the sequence $(g^n\nu_\infty)_{n\in\mathbb{N}}$ is the Dirac mass at $g^{+\infty}$, showing our claim.

Now, for every $k\in \mathsf{G}$, the $\mu$-stationarity of the Dirac mass at $k^{+\infty}$ implies that $$\mu(\{g\in \mathsf{G}|gk^{+\infty}=k^{+\infty}\})=1.$$ Combining \cite[Proposition~7.4.6]{DSU} with separability of $X$ shows that $\Lambda_\infty \mathsf{G}$ contains a dense countable set made of endpoints of axes of loxodromic elements. It follows that $\mu$ gives full measure to the subset made of all elements that act trivially on  $\Lambda_\infty \mathsf{G}$. In other words $\mu$ is supported on $K$, as desired. 
\end{proof}

\begin{rk}
 With a bit more work, we believe that this can be extended to the case where $X$ is merely a complete separable Gromov-hyperbolic space, but we will not need it in the sequel. The difference is that the map $\theta_\infty:\horox_\infty\to\partial_\infty X$ defined by Maher and Tiozzo in \cite[Section~3.2]{MT} is no longer a homeomorphism in this case (it is not necessarily injective), and one has to run the dynamical argument in $\horox_\infty$ and $\partial_\infty X$ in parallel.
\end{rk}

We recall that a $\cat$ space $X$ has the \emph{geodesic extension property} if every geodesic segment in $X$ extends to a bi-infinite geodesic line. The following lemma will help us check the faithfulness assumption from Proposition~\ref{prop:icc-hyperbolic} in the Artin group setting.

\begin{lemma}\label{lemma:criterion-for-faithfulness}
Let $X$ be a $\cat$ space with the geodesic extension property such that $|\partial_\infty X|\ge 3$.
Then the only isometry of $X$ fixing $\partial_\infty X$ pointwise is the identity. 
\end{lemma}

\begin{proof} 
Let $f$ be an isometry that fixes $\partial_\infty X$ pointwise, and let $x\in X$; we aim to prove that $f(x)=x$. As $X$ has the geodesic extension property, there is a complete geodesic line $\ell$ containing $x$. Take $\eta\in\partial_\infty X\setminus\partial_\infty\ell$, and extend the geodesic ray  $\overline{x\eta}$ to a complete geodesic line $\ell'$. As $f$ fixes $\partial_\infty X$ pointwise and $\ell$ is the only geodesic connecting the two points in $\partial_\infty\ell$, we have $f(\ell)=\ell$. Similarly $f(\ell')=\ell'$. So $f(\ell\cap\ell')=\ell\cap\ell'$. However, $\ell\cap\ell'$ is either a point, or a segment, or a ray, so $f$ fixes $\ell\cap\ell'$ pointwise. In particular $f(x)=x$, as desired. 
\end{proof}

\subsection{Specification to Artin groups}

\begin{lemma}\label{lemma:artin-icc}
Let $G=G_\Gamma$ be a 2-dimensional Artin group of hyperbolic type. Then the following statements are equivalent.
\begin{enumerate}
\item $G$ is ICC;
\item $G$ is not virtually isomorphic to the direct sum of $\mathbb{Z}$ and a free group;
\item $|V\Gamma|\ge 3$ and $\Gamma$ is not a star with all edges labeled by $2$.
\end{enumerate}
\end{lemma}

\begin{proof}
The implication $1\Rightarrow 2$ follows from the fact that being ICC is stable under passing to a finite-index subgroup (and $\mathbb{Z}\times F_n$ is not ICC). The implication $2\Rightarrow 3$ follows from the fact that every Artin group whose defining graph is an edge is virtually isomorphic to $\mathbb{Z}\times F_n$ (see e.g.\ \cite[Section~2]{crisp}), and an Artin group associated to a star with all labels labeled by $2$ is isomorphic to $\mathbb{Z}\times F_n$. Finally, the implication $3\Rightarrow 1$ follows from \cite[Theorem~2.35]{DGO}, using that $G_\Gamma$ is acylindrically hyperbolic \cite[Theorem~A]{MP} and torsion-free.
\end{proof}

\begin{lemma}\label{lemma:strong-icc}
Let $G=G_\Gamma$ be a 2-dimensional Artin group of hyperbolic type. Assume that $\Gamma$ is connected and that every vertex of $\Gamma$ has valence at least $2$. 

Then $\Aut(\mathbb{D}_\Gamma)$ is strongly ICC. 
\end{lemma}

\begin{proof}
The assumption made on $\Gamma$ implies that the Deligne complex $\del$ does not have any free face in the sense of \cite[Definition~II.5.9]{BH}, hence has the geodesic extension property by \cite[Proposition~II.5.10]{BH}. Note that $\Lambda_\infty G=\partial_\infty\del$ as the action of $G$ on $\del$ cocompact. Moreover $\partial_\infty\del$ has at least 3 points. Lemma~\ref{lemma:criterion-for-faithfulness} therefore implies that $\Aut(\del)$ acts faithfully on $\partial_\infty\del$, and Proposition~\ref{prop:icc-hyperbolic} then implies that $\Aut(\del)$ is strongly ICC.
\end{proof}

\section{Measure equivalence rigidity and applications}\label{sec:theorems}

\emph{We now complete the proofs of our main results regarding measure equivalence rigidity in the realm of $2$-dimensional Artin groups of hyperbolic type, and give applications to orbit equivalence rigidity and lattice embeddings of Artin groups.}

\subsection{Measure equivalence rigidity}

\subsubsection{A universal factor}

\begin{theo}\label{theo:main-artin}
Let $G_1$ and $G_2$ be two $2$-dimensional Artin groups of hyperbolic type, with connected defining graphs $\Gamma_1$ and $\Gamma_2$, respectively. For every $i\in\{1,2\}$, let $\Theta_{\Gamma_i}$ be the associated fixed set graph. Assume that there exists a measure equivalence coupling $\Sigma$ between $G_1$ and $G_2$.

Then $\Theta_{\Gamma_1}$ and $\Theta_{\Gamma_2}$ are isomorphic, and there exists an almost $(G_1\times G_2)$-equivariant Borel map $$\Sigma\to\Isom(\Theta_{\Gamma_1}\to\Theta_{\Gamma_2}),$$ where the $G_1\times G_2$-action on $\Isom(\Theta_{\Gamma_1}\to\Theta_{\Gamma_2})$ is via $(g_1,g_2)\cdot f (v):=g_2 f(g_1^{-1}v)$.
\end{theo}

\begin{proof}
If $G_1$ or $G_2$ is isomorphic to $\mathbb{Z}^2$, then the conclusion holds because $\mathbb{Z}^2$ is the only amenable $2$-dimensional Artin group of hyperbolic type, and it has empty fixed set graph. In all other cases, the theorem follows from Theorem~\ref{theo:blueprint} together with the fact that the fixed set graph $\Theta_{\Gamma_i}$ is a curve-graph-like $G_i$-graph (Proposition~\ref{prop:curve-graph-like}).
\end{proof}

 The assumption that $\Gamma_1$ and $\Gamma_2$ are connected is harmless, as the behavior of measure equivalence under free products has been thoroughly studied by Alvarez and Gaboriau in \cite{AG}. In fact, with a bit more work, one can use \cite{AG} to derive that the first conclusion of our theorem (namely, if $G_1$ and $G_2$ are measure equivalent, then $\Theta_{\Gamma_1}$ and $\Theta_{\Gamma_2}$ are isomorphic) still holds without assuming connectedness of $\Gamma_1$ and $\Gamma_2$.

\subsubsection{Measure equivalence classification results}

Combining Theorem~\ref{theo:main-artin} with the rigidity statement for fixed set graphs given in Theorem~\ref{theo:rigid-crisp}, we deduce the following (note that $\Out(G_1)$ being finite implies the defining graph is connected and has no separating vertices or edges by \cite{crisp}).

\begin{cor}\label{cor:me-clttf}
Let $G_1$ and $G_2$ be two large type Artin groups whose defining graphs are triangle-free and have more than 2 vertices. Suppose $\Out(G_1)$ and $\Out(G_2)$ are finite.

Then $G_1$ and $G_2$ are measure equivalent if and only if they are isomorphic.
\qed
\end{cor}

As a consequence of Theorem~\ref{theo:main-artin}, we also get a measure equivalence classification statement among right-angled Artin groups, which has since been improved in \cite{HH} where we removed the girth assumption. 

\begin{cor}\label{cor:me-raag} 
Let $G_1$ and $G_2$ be two right-angled Artin groups whose underlying graphs $\Gamma_1$ and $\Gamma_2$ are connected and have girth at least $5$. Assume that $\Out(G_1)$ and $\Out(G_2)$ are finite.

Then $G_1$ and $G_2$ are measure equivalent if and only if they are isomorphic.
\end{cor}

\begin{proof}
The if direction is clear. For the only if direction, since $G_1$ and $G_2$ are $2$-dimensional Artin groups of hyperbolic type, it follows from Theorem~\ref{theo:main-artin} that the graphs $\Theta_{\Gamma_1}$ and $\Theta_{\Gamma_2}$ are isomorphic. Notice that in our situation, the fixed set graph is isomorphic to the extension graph introduced by Kim and Koberda in \cite{KK}. The conclusion thus follows from \cite[Theorem~1.2]{Hua2} (alternatively, this also follows from \cite[Theorem~7.1]{BKS}).
\end{proof}

\subsubsection{Measure equivalence rigidity results}

\begin{theo}\label{theo:main-furman}
Let $G=G_\Gamma$ be a $2$-dimensional Artin group of hyperbolic type. Assume that $\Gamma$ is connected and has no valence one vertex, and that the comparison map $\Aut(\Cay_\Gamma)\to\Aut(\Theta_\Gamma)$ is an isomorphism. Let $H$ be a countable group which is measure equivalent to $G_\Gamma$.

Then there is a homomorphism $H\to\Aut(\Cay_\Gamma)$ with finite kernel, and whose image is a lattice.
\end{theo}

\begin{proof}
Since the comparison map $\Aut(\Cay_\Gamma)\to\Aut(\Theta_\Gamma)$ is an isomorphism (in fact an isomorphism of topological groups by Lemma~\ref{lemma:comparison-is-iso}), the group $\Aut(\Theta_\Gamma)$ is a locally compact second countable group that contains $G_\Gamma$ as a lattice. In addition, Lemma~\ref{lemma:comparison-is-iso} implies that $\Aut(\Cay_\Gamma)$ is isomorphic (as a topological group) to $\Aut(\del)$. By Lemma~\ref{lemma:strong-icc}, the group $\Aut(\Theta_\Gamma)$ is therefore strongly ICC. Since $\Aut(\Theta_\Gamma)$ is a universal $(G_\Gamma,G_\Gamma)$-factor (Theorem~\ref{theo:main-artin}), the conclusion follows by applying Theorem~\ref{theo:bfs}. 
\end{proof}

Note that the $G_\Gamma$-action on $\Theta_\Gamma$ yields a homomorphism $G_\Gamma\to \Aut(\Theta_\Gamma)$.

\begin{theo}\label{theo:main-kida}
Let $G=G_\Gamma$ be a $2$-dimensional Artin group of hyperbolic type. Assume that $\Gamma$ is connected, and the natural map $G_\Gamma\to\Aut(\Theta_\Gamma)$ is injective and has finite-index image (in particular $\Aut(\Theta_\Gamma)$ is countable). Let $H$ be a countable group which is measure equivalent to $G_\Gamma$.

Then $H$ is almost isomorphic to $G_\Gamma$. 
\end{theo}

\begin{proof}
Since the natural map $G_\Gamma\to\Aut(\Theta_\Gamma)$ is injective, the group $G_\Gamma$ is not commensurable to the direct sum of $\mathbb{Z}$ and a free group (otherwise $\Theta_\Gamma$ is reduced to a point). By Lemma~\ref{lemma:artin-icc}, the group $G_\Gamma$ is ICC, and by Theorem~\ref{theo:main-artin}, the group $\Aut(\Theta_\Gamma)$ is a universal $(G_\Gamma,G_\Gamma)$-factor, so the conclusion follows from Corollary~\ref{cor:bfs}. 
\end{proof}

\begin{cor}\label{cor:main-large}
Let $G_\Gamma$ be a large type Artin group, such that $\Out(G_\Gamma)$ is finite. Assume that $\Gamma$ is triangle-free and has more than 2 vertices. Let $H$ be a countable group which is measure equivalent to $G_\Gamma$.
	\begin{enumerate}
		\item There is a homomorphism $H\to\Aut(\Cay_\Gamma)$ with finite kernel, with image a lattice.
		\item Assuming additionally that $\Gamma$ is star-rigid, then $H$ is almost isomorphic to $G_\Gamma$.
	\end{enumerate}
\end{cor}

\begin{proof}
Since $\Out(G_\Gamma)$ is finite, the graph $\Gamma$ is connected and has no separating vertex or edge \cite{crisp}. Our assumptions therefore imply that the comparison map $\Aut(\Cay_\Gamma)\to\Aut(\Theta_\Gamma)$ is an isomorphism (Theorem~\ref{theo:rigid-crisp-2}), so the first conclusion follows from Theorem~\ref{theo:main-furman}. The second conclusion follows from Theorem~\ref{theo:main-kida} and Theorem~\ref{theo:rigid-crisp-2}.
\end{proof}

\subsection{Orbit equivalence rigidity}\label{sec:oe}

Let $G_1$ and $G_2$ be two countable groups, and for every $i\in\{1,2\}$, let $X_i$ be a standard probability space equipped with a free measure-preserving action of $G_i$ by Borel automorphisms. 

The actions $G_1\actson X_1$ and $G_2\actson X_2$ are \emph{virtually conjugate} if there exist short exact sequences $1\to F_i\to G_i\to Q_i\to 1$ with $F_i$ finite for every $i\in\{1,2\}$, finite-index subgroups $Q_1^0\subseteq Q_1$ and $Q_2^0\subseteq Q_2$, and conjugate actions $Q_1^0\actson X'_1$ and $Q_2^0\actson X'_2$ such that for every $i\in\{1,2\}$, the $Q_i$-action on $X_i/F_i$ is induced from the $Q_i^0$-action on $X'_i$ in the sense of \cite[Definition~2.1]{Kid-oe}.

The actions $G_1\actson X_1$ and $G_2\actson X_2$ are \emph{stably orbit equivalent} if there exist Borel subsets $Y_1\subseteq X_1$ and $Y_2\subseteq X_2$ of positive measure and a measure-scaling isomorphism $f:Y_1\to Y_2$ such that for a.e.\ $y\in Y_1$, one has $$f((G_1\cdot y)\cap Y_1)=(G_2\cdot f(y))\cap Y_2.$$ Virtually conjugate actions are always stably orbit equivalent.

It turns out that two countable groups are measure equivalent if and only if they admit essentially free ergodic measure-preserving actions by Borel automorphisms on standard probability spaces which are stably orbit equivalent (see \cite{Fur-oe} and \cite[Theorem~2.1]{Gab2}), so all our theorems on measure equivalence rigidity can be rephrased in this language.

Let $G$ be a countable group, and let $X$ be a standard probability space. An ergodic measure-preserving free $G$-action on $X$ by Borel automorphisms is \emph{OE-superrigid} if for every countable group $H$, and every ergodic measure-preserving free $H$-action on a standard probability space $Y$, if the actions $G\actson X$ and $H\actson Y$ are stably orbit equivalent, then they are virtually conjugate.

Again, there is a close relationship between ME-superrigidity of a countable group and OE-superrigidity of its ergodic probability measure-preserving free actions, given in \cite[Lemma~4.18]{Fur-survey}, and as in Theorem~\ref{theo:main-kida} we get the following.

\begin{theo}\label{theo:oe-superrigid}
Let $G=G_\Gamma$ be a $2$-dimensional Artin group of hyperbolic type with defining graph $\Gamma$. Assume that $\Gamma$ is connected, and that the natural map $G\to\Aut(\Theta_\Gamma)$ is injective and has finite-index image.

Then every ergodic measure-preserving essentially free action of $G$ on a standard probability space by Borel automorphisms is OE-superrigid.
\qed
\end{theo}

Notice that as in Corollary~\ref{cor:main-large}, this applies in particular to all Artin groups $G_\Gamma$ of large type, with $\Out(G_\Gamma)$ finite, whose defining graph is triangle-free and star-rigid. 

\begin{rk}\label{rk:aperiodic}
 There are possible variations on Theorem~\ref{theo:oe-superrigid}. Let $G=G_\Gamma$ be as in Theorem~\ref{theo:oe-superrigid}. Let $X$ be a standard probability space, and let $G\actson X$ be a free, measure-preserving action by Borel automorphisms which is further assumed to be \emph{aperiodic}, i.e.\ every finite index subgroup of $G$ acts ergodically on $X$. Let $H\actson Y$ be a free, ergodic, measure-preserving action on a standard probability space. If the actions $G\actson X$ and $H\actson Y$ are orbit equivalent, then they are actually conjugate. See e.g.\ \cite[Theorem~1.2(ii)]{Kid-oe}.
\end{rk}

\subsection{Lattice embeddings}

Our next two theorems describe all possible lattice embeddings of $2$-dimensional Artin groups of hyperbolic type, under different assumptions on $\Gamma$.

\begin{theo}\label{theo:lattice-embeddings}
Let $G=G_\Gamma$ be a $2$-dimensional Artin group of hyperbolic type, and let $G'$ be a finite-index subgroup of $G$. Assume that $\Gamma$ is connected and has no valence $1$ vertex, and that the comparison map $\Aut(\Cay_\Gamma)\to\Aut(\Theta_\Gamma)$ is an isomorphism. Let $\mathsf{H}$ be a locally compact second countable group, and let $f:G'\to \mathsf{H}$ be an injective homomorphism, such that $f(G')$ is a lattice in $\mathsf{H}$. 

Then there exists a continuous homomorphism $\psi:\mathsf{H}\to \Aut(\Cay_\Gamma)$ with compact kernel such that $\psi\circ f$ coincides with the natural map $G'\to\Aut(\Cay_\Gamma)$. 
\end{theo}

\begin{proof}
The proof is essentially the same as Kida's proof in the mapping class group setting \cite[Theorem~8.1]{Kid}, or Furman's proof for higher-rank lattices \cite{Fur}.

Using the lattice embedding, we view $\mathsf{H}$ as a self-coupling of $G'$ (with the actions by left and right multiplication, through $f$). By Proposition~\ref{prop:curve-graph-like} and Remark~\ref{rk:curve-graph-finite-index}, the fixed set graph $\Theta_\Gamma$ is a curve-graph-like $G'$-graph. It thus follows from Theorem~\ref{theo:blueprint} that there is an almost $(G'\times G')$-equivariant Borel map $\mathsf{H}\to\Aut(\Theta_\Gamma)$. Since the comparison map $\Aut(\Cay_\Gamma)\to\Aut(\Theta_\Gamma)$ is an isomorphism (in fact, an isomorphism of topological groups by Lemma~\ref{lemma:comparison-is-iso}), we thus have an almost $(G'\times G')$-equivariant Borel map $\varphi:\mathsf{H}\to\Aut(\Cay_\Gamma)$. Let $X=\mathsf{H}/G'$, where we mod out by the action of the rightmost factor. The map 
\begin{displaymath}
\begin{array}{cccc}
\Phi:& \mathsf{H}\times \mathsf{H} & \to & \Aut(\Cay_\Gamma) \\
& (h_1,h_2) & \mapsto & \varphi(h_1^{-1})^{-1}\varphi(h_1^{-1}h_2)\varphi(h_2)^{-1} 
\end{array}
\end{displaymath}
\noindent has the following invariance properties: for all $g\in G'$ and a.e.\ $(h_1,h_2)\in \mathsf{H}\times \mathsf{H}$, we have $\Phi(h_1g,h_2)=\Phi(h_1,h_2g)=\Phi(h_1,h_2)$, and $\Phi(gh_1,gh_2)=g\Phi(h_1,h_2)g^{-1}$. In other words $\Phi$ induces a Borel map $\overline{\Phi}:X\times X\to\Aut(\Cay_\Gamma)$, and the $G'$-invariant probability measure on $X\times X$ pushes forward to a probability measure on $\Aut(\Cay_\Gamma)$ which is invariant under the action of $G'$ by conjugation.  Now $\Aut(\Cay_\Gamma)$ is isomorphic (as a topological group) to $\Aut(\del)$ (Lemma~\ref{lemma:comparison-is-iso}), so by Lemma~\ref{lemma:strong-icc} it is strongly ICC. In particular, as $G'$ is a lattice in $\Aut(\Cay_\Gamma)$, it follows from \cite[Lemma~A.6]{BFS} that the only probability measure on $\Aut(\Cay_\Gamma)$ which is invariant under conjugation by $G'$ is the Dirac mass at the identity. This implies that $\Phi$ is almost everywhere constant, with value the neutral element of $\Aut(\Cay_\Gamma)$. In other words, we have proved that for a.e.\ $(h_1,h_2)\in \mathsf{H}\times \mathsf{H}$, one has $\varphi(h_1h_2)=\varphi(h_1)\varphi(h_2)$.
 
Now, using \cite[Theorems~B.2 and~B.3]{Zim}, there exists a continuous homomorphism $\psi:\mathsf{H}\to\Aut(\Cay_\Gamma)$ which coincides almost everywhere with $\varphi$. Now, for all $g\in G'$ and a.e.\ $h\in \mathsf{H}$, we have $$\psi(f(g))\psi(h)=\psi(f(g)h)=\varphi(f(g)h)=g\varphi(h)=g\psi(h),$$ showing that $\psi\circ f$ coincides with the natural map $G'\to\Aut(\Cay_\Gamma)$. Finally, the preimage $\varphi^{-1}(\{\mathrm{id}\})$ is a Borel fundamental domain for the leftmost action of $G$ on $\mathsf{H}$, so $\psi^{-1}(\{\mathrm{id}\})$ is a closed subgroup of finite Haar measure, whence compact. 
\end{proof}

We recall that the \emph{FC-center} of a group $G$ is the characteristic subgroup made of all elements whose conjugacy class is finite. We mention that the subgroup $F$ arising in the following statement is a finite normal subgroup of $\Aut(\Theta_\Gamma)$; in some concrete situations where $\Aut(\Theta_\Gamma)$ coincides with $\Aut(G_\Gamma)$, one could actually take advantage of the work of Crisp to show that $F$ is trivial.

\begin{theo}\label{theo:lattice-embeddings-2}
Let $G=G_\Gamma$ be a $2$-dimensional Artin group of hyperbolic type with connected defining graph, and let $G'$ be a finite-index subgroup of $G$. Assume that the natural homomorphism $G\to\Aut(\Theta_\Gamma)$ is injective and has finite index image, and let $F\unlhd\Aut(\Theta_\Gamma)$ be the FC-center of $\Aut(\Theta_\Gamma)$. Let $\mathsf{H}$ be a locally compact second countable group (equipped with its right Haar measure), and let $f:G'\to \mathsf{H}$ be an injective homomorphism, such that $f(G')$ is a lattice in $\mathsf{H}$. 

Then there exists a continuous homomorphism $g:\mathsf{H}\to \Aut(\Theta_\Gamma)/F$ with compact kernel such that $g\circ f$ coincides with the natural inclusion $G'\to\Aut(\Theta_\Gamma)/F$.
\end{theo}

\begin{proof}
The group $\Aut(\Theta_\Gamma)$ is a universal $(G',G')$-factor, and therefore so is $\Aut(\Theta_\Gamma)/F$. Notice that $F$ is finite because $G_\Gamma$ is ICC (Lemma~\ref{lemma:artin-icc}) and contained as a finite index subgroup of $\Aut(\Theta_\Gamma)$. Therefore $\Aut(\Theta_\Gamma)/F$ is ICC. The proof is then the same as the proof of Theorem~\ref{theo:lattice-embeddings}.
\end{proof}

\section{Rigidity for von Neumann algebras}\label{sec:von-neumann}

\emph{In this section, we prove that cross-product von Neumann algebras of free, ergodic, probability measure-preserving actions of many Artin groups are Cartan rigid, and derive a $W^*$-superrigidity theorem.} \\

Given a countable group $G$ and a standard probability space $(X,\mu)$ equipped with an ergodic measure-preserving free $G$-action by Borel automorphisms, a celebrated construction of Murray and von Neumann \cite{MvN} associates a cross-product von Neumann algebra $L^\infty(X)\rtimes G$ to the $G$-action on $X$, in which $L^\infty(X)$ sits as a \emph{Cartan subalgebra}, i.e.\ a maximal abelian subalgebra whose normalizer generates the whole algebra.

A countable group $G$ is \emph{Cartan rigid} if for every action $G\actson X$ as above, the subalgebra $L^\infty(X)$ is, up to unitary conjugacy, the unique Cartan subalgebra of $L^\infty(X)\rtimes G$. In a groundbreaking work \cite{PV2}, Popa and Vaes proved that finitely generated non-abelian free groups are Cartan rigid. Here we will use a theorem of Ioana \cite{Ioa} regarding more generally Cartan rigidity of amalgamated free products to show the following theorem.

\begin{theo}\label{theo:cartan-rigid}
Let $G=G_\Gamma$ be a $2$-dimensional Artin group of hyperbolic type, and assume that $\Gamma$ not a clique and that $G_\Gamma$ is not isomorphic to the direct sum of $\mathbb{Z}$ and a free group (i.e.\ $\Gamma$ is also not a star with all edges labeled by $2$).

Then $G$ is Cartan rigid.
\end{theo}

\begin{proof}
As $\Gamma$ is not a clique, we can find a vertex $v\in V\Gamma$ whose star is a proper subgraph of $\Gamma$. Then $G_\Gamma$ splits as an amalgamated free product $G_\Gamma=G_{\st(v)}\ast_{G_{\lk(v)}}G_{\Gamma\setminus\{v\}}$: indeed it follows from Lemma~\ref{lem:injective} that $G_{\st(v)}$ and $G_{\Gamma\setminus\{v\}}$ naturally embed as subgroups of $G\Gamma$ intersecting in $G_{\lk(v)}$, and the amalgamated free product structure follows from the presentation of $G_\Gamma$. Our choice of $v$ ensures that the edge group has infinite index in both vertex groups. By \cite[Theorem~1.1]{Ioa2}, it is enough to find a finite subset $F\subseteq G$ such that $\bigcap_{g\in F}gG_{\lk(v)}g^{-1}=\{1\}$. This is the contents of our next lemma.
\end{proof}

\begin{lemma}\label{lemma:cartan}
Let $G=G_\Gamma$ be a $2$-dimensional Artin group of hyperbolic type. Let $v\in V\Gamma$ be a vertex, and assume that $\st(v)\subsetneq\Gamma$.

Then there exists a finite subset $F\subseteq G$ such that $\bigcap_{g\in F}gG_{\lk(v)}g^{-1}=\{1\}$.
\end{lemma}  

\begin{proof}
We consider the modified Deligne complex $\del$ equipped with the $\mathrm{CAT}(0)$ metric that was recalled in Section~\ref{sec:background-deligne}. Let $C\subset\del$ be the subcomplex spanned by all rank $0$ vertices associated to elements of $G_{\lk(v)}$. Then $C$ is a convex subcomplex of $\del$. Let $w\in V\Gamma\setminus\st(v)$. Consider the element $g=vw\in G$ (identifying vertices of $\Gamma$ with the associated generators of $G_\Gamma$). We start by understanding the translates $g^k C$ for large $k$. For this purpose, we first construct an explicit axis of $g$ with respect to the action $G\actson \del$ as follows. For an integer $k$, consider the edge path $\ell_k$ passing through the vertices of $\del$ in the following order: $(vw)^k$, $(vw)^{k}G_v$, $(vw)^{k}v$, $(vw)^{k}vG_w$, $(vw)^{k+1}$ (here $(vw)^k$ denotes the rank 0 vertex of $\del$ represented by $(vw)^k$, $(vw)^{k}G_v$ denotes the rank 1 vertex of $\del$ represented by the coset $(vw)^{k}G_v$ etc). Let $\ell=\cup_{k=-\infty}^{\infty}\ell_k$. As $v$ and $w$ are not adjacent in $\Gamma$, the discussion on angular metric of links in Section~\ref{sec:background-deligne} implies that consecutive edges in $\ell$ form an angle at least equal to $\pi$ at their common vertex. Thus $\ell$ is a local geodesic (see \cite[p.\ 60, I.5.7 and p.\ 103, I.7.16]{BH}), hence is a global geodesic by \cite[p.\ 201, II.4.14]{BH}. Thus $\ell$ is an axis of $g=vw$. Note that $\ell\cap C=\{x_0\}$ where $x_0$ is the rank 0 vertex corresponding to the identity element of $G$. As the distance function $d(x,C)$ with $x$ moving along $\ell$ is convex (see \cite[p.\ 178, II.2.5]{BH}), we deduce that $\partial_\infty \ell\cap\partial_\infty C=\emptyset$.

Let $\pi: \del\to \ell$ be the nearest point projection. We claim $\pi(C)$ is bounded. As $\partial_\infty C$ is a closed subset of $\partial_\infty X$ (with respect to the cone topology) and $\partial_\infty \ell\cap\partial_\infty C=\emptyset$, we know from the Gromov hyperbolicity of $\del$ that there is $K>0$ such that diam$(\pi(r))\le K$ for any each geodesic ray $r\subset C$ starting from $x_0$. As $C$ is cocompact and $C$ contains an axis of some loxodromic element (when $C$ is not bounded), it follows that $C$ is contained in a finite neighborhood of the union of all geodesic lines in $C$. By applying the thin triangle property to each ideal triangle with three vertices being $x_0$ and $\partial \ell'$ where $\ell'\subset C$ is a geodesic line, we deduce that $C$ is contained in a finite neighborhood of the unions geodesic rays in $C$ emanating from $x_0$. Thus the claim follows.

For each integer $k$, let $C_k\subset C$ be the collection of points $x$ such that $d(x,g^k C)=d(C,g^k C)$. We claim there exists $k$ such that $C_k$ is bounded. By \cite[Lemma 2.10 and Remark 2.13]{Hua3}, $C_k$ is a nonempty convex subset, moreover $C_k\subset C$ and a subset $C'_k\subset g_kC$ cobound a convex subset isometric to $C_k\times [0,d(C,g^k C)]$. From the previous paragraph we know that $\lim_{k\to\infty}d(C,g^k C)=\infty$. Thus $C_k$ being unbounded for all $k$ will violate the hyperbolicity of $\del$.

As $C$ is $G_{\lk(v)}$-invariant and $g^k C$ is $g^kG_{\lk(v)}g^{-k}$-invariant, we deduce that $C_k$ is preserved by $G_{\lk(v)}\cap g^kG_{\lk(v)}g^{-k}$. As $C_k$ is bounded, let $x$ be the circumcenter of $C_k$. Then $G_{\lk(v)}\cap g^kG_{\lk(v)}g^{-k}$ fixes $x$. The stabilizer of $x$ is either trivial (in which case we are done), or a conjugate of $G_s$ for some vertex $s\in\Gamma$, or a conjugate of $G_{\overline{st}}$ for some edge $\overline{st}\subset\Gamma$ (in which case $x$ is a rank 2 vertex). It remains to find a finite subset $F'\subseteq G$ such that $\bigcap_{h\in F'}h\stab(x)h^{-1}=\{1\}$. If $\stab(x)=G_s$, then let $h$ be a standard generator which does not commute with $s$ (such $h$ exists as $G_\Gamma$ is not isomorphic to $F_k\times \mathbb Z$). It follows from \cite[Corollary~3.8]{godelle} that $G_s\cap h^{-1}G_s h=\{1\}$.
If $x$ is a rank 2 vertex, as $\Gamma$ is not a clique, there exists $h\in G$ such that $hx\neq x$. Then $\stab_G(x)\cap h\stab_G(x)h^{-1}$ fixes every point in the geodesic segment $s$ from $x$ to $hx$. As there are points in $s$ which are not rank 2 vertices, we are reduced to the previous case.
\end{proof}

Theorem~\ref{theo:cartan-rigid} has applications to rigidity for cross-product von Neumann algebras associated to ergodic actions of Artin groups. 

\begin{cor}\label{cor:von-neumann-strong}
Let $G=G_\Gamma$ be a $2$-dimensional Artin group of hyperbolic type. Assume that its defining graph $\Gamma$ is connected and is not a clique. Assume that the natural map $G\to\Aut(\Theta_\Gamma)$ is injective and has finite-index image.

Let $H$ be a countable group, and let $G\actson X$ and $H\actson Y$ be free, ergodic, measure-preserving actions by Borel automorphisms on standard probability spaces.  

If $L^\infty(X)\rtimes G$ is isomorphic to $L^\infty(Y)\rtimes H$, then $G\actson X$ and $H\actson Y$ are virtually conjugate (in particular $G$ and $H$ are almost isomorphic).
\end{cor}

\begin{proof}
We know from Theorem~\ref{theo:cartan-rigid} that $L^\infty(X)$ is, up to unitary conjugacy, the unique Cartan subalgebra of $L^\infty(X)\rtimes G$. It thus follows from a theorem of Singer \cite{Sin} that $L^\infty(X)\rtimes G$ recovers the orbit equivalence class of the action $G\actson X$, i.e.\ $G\actson X$ and $H\actson Y$ are orbit equivalent. The conclusion thus follows from Theorem~\ref{theo:oe-superrigid}.
\end{proof}

\begin{rk}
As in Remark~\ref{rk:aperiodic}, if the $G$-action on $X$ is aperiodic, then we reach the stronger conclusion that the actions $G\actson X$ and $H\actson Y$ are conjugate.
\end{rk}

 Recall that a subgroup $H$ of an Artin group $G=G_\Gamma$ is a \emph{parabolic subgroup} if there exists an induced subgraph $\Lambda\subseteq\Gamma$ such that $H$ is conjugate to $G_\Lambda$, viewed as a subgroup of $G_\Gamma$ via Lemma~\ref{lem:injective}. The following theorem shows that Cartan rigidity holds for more general classes of Artin groups, beyond the case of $2$-dimensional Artin groups of hyperbolic type.

\begin{theo}\label{theo:cartan-rigid-2}
Let $G=G_\Gamma$ be an Artin group. Assume that $\Gamma$ not a clique, and further assume that
\begin{enumerate}
\item the intersection of any two parabolic subgroups is again a parabolic subgroup;
\item every non-increasing chain of parabolic subgroups is stationary;
\item no proper nontrivial parabolic subgroup is normal.
\end{enumerate}
Then $G$ is Cartan rigid.
\end{theo}

The assumptions of Theorem~\ref{theo:cartan-rigid-2} hold for all large type Artin groups whose underlying graph is not a clique by recent work of Cumplido, Martin and Vaskou \cite{CMV}. These are all $2$-dimensional, but include examples that are not of hyperbolic type (when $\Gamma$ contains a triangle with all edges labeled $3$).

\begin{proof}
Following the proof of Theorem~\ref{theo:cartan-rigid} and keeping its notations, it is enough to find a finite subset $F\subseteq G$ such that $\bigcap_{g\in F}gG_{\lk(v)}g^{-1}=\{1\}$. Enumerate $G$ as $G=\{g_n\}_{i\in\mathbb{N}}$. For every $n\in\mathbb{N}$, the intersection $\bigcap_{i\le n}g_iG_{\lk(v)}g_i^{-1}$ is a parabolic subgroup, and the chain condition ensures that this intersection is constant for $i$ large enough, with value a proper normal parabolic subgroup $P$. It follows that $P=\{1\}$, which concludes the proof.
\end{proof}

\section{Cocycle rigidity from higher-rank lattices to Artin groups}\label{sec:cocycle}

\emph{This final short section is independent from our main measure equivalence rigidity theorems. We exploit our techniques to derive a superrigidity theorem for cocycles from higher-rank lattices to Artin groups.}\\

Using the action on the $\cat$ version of the modified Deligne complex, one can show that if  $\Lambda$ is a lattice in a higher rank simple algebraic group over a local field, and if $G$ is a $2$-dimensional Artin group of hyperbolic type, then every homomorphism from $\Lambda$ to $G$ has finite image -- for example, this can be proved using Haettel's theorem stating that higher-rank lattices cannot act nonelementarily on hyperbolic spaces \cite{Hae}. Theorem~\ref{theo:cocycle} below extends this statement about homomorphisms to a cocycle superrigidity result.

In the following statement, a \emph{cocycle} $c:\Lambda\times\Omega\to G$ is a Borel map such that for all $\lambda_1,\lambda_2\in \Lambda$ and a.e.\ $\omega\in\Omega$, one has $c(\lambda_1\lambda_2,\omega)=c(\lambda_1,\lambda_2\omega)c(\lambda_2,\omega)$. Two cocycles $c$ and $c'$ are \emph{cohomologous} if there exists a Borel map $f:\Omega\to G$ such that for all $\lambda\in \Lambda$ and a.e.\ $\omega\in\Omega$, one has $c'(\lambda,\omega)=f(\lambda\omega)c(\lambda,\omega)f(\omega)^{-1}$.

\begin{theo}\label{theo:cocycle}
Let $G=G_\Gamma$ be a $2$-dimensional Artin group of hyperbolic type. Let $\mathsf{H}$ be a product of connected higher rank simple algebraic groups over local fields. Let $\Lambda$ be either $\mathsf{H}$, or a lattice in $\mathsf{H}$. Let $\Omega$ be a standard probability space equipped with an ergodic measure-preserving $\Lambda$-action by Borel automorphisms.

Then every cocycle $\Lambda\times\Omega\to G_\Gamma$ is cohomologous to the trivial cocycle.
\end{theo}

\begin{proof}
Using a theorem that was announced by Bader and Furman that deals with the case of free products (\cite{BF}, see also \cite[Theorem~5.1]{GHL}), we can assume that $\Gamma$ is connected. We can further assume that $\Gamma$ is not an edge labeled by $2$, as otherwise $G_\Gamma$ is isomorphic to $\mathbb{Z}^2$ and the conclusion follows from \cite[Theorem~9.1.1]{Zim-book}. 

Following \cite[Definition~3.1]{GHL}, given a countable discrete group $\Upsilon$, we say that the pair $(\Lambda,\Upsilon)$ is \emph{cocycle-rigid} if for every standard probability space $\Omega$ equipped with an ergodic measure-preserving $\Lambda$-action, every cocycle $\Lambda\times\Omega\to \Upsilon$ is cohomologous to a cocycle that takes its values in a finite subgroup of $\Upsilon$. Our goal is to prove that the pair $(\Lambda,G_\Gamma)$ is cocycle-rigid (recall indeed that $G_\Gamma$ is torsion-free so the only finite subgroup of $G_\Gamma$ is the trivial subgroup).

Recall that the fixed set graph $\Theta_\Gamma$ is a curve-graph-like $G_\Gamma$-graph (Proposition~\ref{prop:curve-graph-like}). Assertion~2 in the definition of a curve-graph-like $G_\Gamma$-graph (Definition~\ref{de:curve-graph-like}) precisely says that $G_\Gamma$ is \emph{geometrically rigid} with respect to $\calp_\ast(V(\Theta_\Gamma))$ in the sense of \cite[Definition~4.9]{GHL}. It thus follows from \cite[Theorem~2]{GHL} that every cocycle $\Lambda\times\Omega\to G_\Gamma$ is cohomologous to a cocycle that takes its values in a subgroup of $G_\Gamma$ that virtually fixes a vertex of $\Theta_\Gamma$. But vertex stabilizers of $\Theta_\Gamma$ are virtually isomorphic to $F\times\mathbb{Z}$ where $F$ is a free group (Lemma~\ref{lemma:stab-vertices-theta}). Since $(\Lambda,F\times\mathbb{Z})$ is cocycle-rigid, and since cocycle-rigidity, as a property of the target group, is stable under passing to a finite-index overgroup \cite[Proposition~3.9]{GHL}, the conclusion follows.
\end{proof}

\footnotesize

\bibliographystyle{alpha}
\bibliography{ME-Artin-bib}

\begin{flushleft}
Camille Horbez\\
Universit\'e Paris-Saclay, CNRS,  Laboratoire de math\'ematiques d'Orsay, 91405, Orsay, France \\
\emph{e-mail:}\texttt{camille.horbez@universite-paris-saclay.fr}\\[8mm]
\end{flushleft}

\begin{flushleft}
Jingyin Huang\\
Department of Mathematics\\
The Ohio State University, 100 Math Tower\\
231 W 18th Ave, Columbus, OH 43210, U.S.\\
\emph{e-mail:~}\texttt{huang.929@osu.edu}\\
\end{flushleft}

\end{document}